\def\versiondate{25 Oct. 2018}
\input math.macros
\input Ref.macros

\checkdefinedreferencetrue
\continuousfigurenumberingtrue
\theoremcountingtrue
\sectionnumberstrue
\forwardreferencetrue
\tocgenerationtrue
\citationgenerationtrue
\nobracketcittrue
\hyperstrue
\initialeqmacro

\input\jobname.key
\bibsty{myapalike}

\def\gh{G}  
\def\gp{\Gamma}  
\def\gpe{\gamma}  
\def\verts{{\ss V}}

\def\edges{{\ss E}}

\def\tr{{\rm Tr}}  

\def\bfz{{\bf 0}}

\def\fo{{\frak F}}  
\def\wsf{{\ss WSF}}
\def\fsf{{\ss FSF}}
\def\UST{{\ss UST}}
\def\ust{{\ss UST}}
\def\bdv{\partial_{\verts}}

\def\Aut{{\rm Aut}} 
\def\ip#1{(\changecomma #1)}
\def\bigip#1{\bigl(\changecomma #1\bigr)}
\def\Bigip#1{\Bigl(\changecomma #1\Bigr)}

\def\changecomma#1,{#1,\,}
\def\bigchangecomma#1,{#1,\;}
\def\leftchangecomma#1,{#1,\ }

\def\GG{{\cal G}}
\def\PGW{{\ss PGW}}
\def\AGW{{\ss AGW}}
\def\UGW{{\ss UGW}}

\def\rtd{\mu}  
\def\cd{\Rightarrow}  

\def\frac#1#2{{#1 \over #2}}

\def\mathrm#1{{\rm #1}}
\def\A{{\cal A}}

\def\marks{{\Xi}}  
\def\mk{\xi}  
\def\mkmp{\psi}  

\def\fdom{{\ss FC}}  
\def\Ph{{\widehat{\rtd}}}

\def\ev#1{{\cal #1}} 
\def\traj{\theta}    
\def\invar{{\cal I}}  
\def\SS{{\cal S}}   
\def\rln{{\cal R}}   
\def\rlnmap{\phi}
\def\assumptions{Let\/ $\P$ be a percolation on a
unimodular random network}
\def\genassumptions{Let\/ $\P$ be an insertion-tolerant percolation on a
unimodular random network}
\def\ct{n}  
\def\paths{{\cal P}_*}   
\def\fin{{\ss FG_2}}  
\def\gtwo{{{\cal G}_{**}}}  
\def\uni{{\cal U}}  
\def\projection{{\ss prj}}  
\def\dmn{\GG_\#}  
\def\itimes{\boxtimes}  
\def\dist{{\rm dist}}  
\def\bn{{\rm B}}  
\def\alg{{\ss Alg}}  
\def\affalg{\overline{\ss Alg}}  
\def\two{\{0, 1\}}
\def\expdeg{\overline{\rm deg}}  
\def\Giso{{\cal H}}
\def\Gisop{{\cal H}^{{\rm min}}}

\def\ins{\Pi}   
\def\isoe{\iota_\edge}

\def\pc{p_{\rm c}}

\def\pu{p_{\rm u}}

\def\UST{{\ss UST}}
\def\MST{{\ss MST}}
\def\WUSF{{\ss WUSF}}
\def\WMSF{{\ss WMSF}}
\def\FUSF{{\ss FUSF}}
\def\FMSF{{\ss FMSF}}
\def\verts{{\ss V}}
\def\vertex{{\ss V}}
\def\vertices{{\ss V}}
\def\edges{{\ss E}}
\def\edge{{\ss E}}
\def\Stab{{\ss Stab}}
\def\bp{o}
\def\lift{\lambda}  
\def\dual#1{#1^\dagger}
\def\collec{{\cal C}}  
\def\K{{\cal K}}  
\def\gve{{\cal G}_{*-}}  
\def\type{\tau}  
\def\HH{{\cal H}}   
\def\DD{{\scr D}}   
\def\dmn{{\scr D}}   
\def\spm{\sigma}  

\def\BLPSusf{\ref b.BLPS:usf/, hereinafter referred to as %
\def\BLPSusf{BLPS (\refbyear {BLPS:usf})}\BLPSusf}

\def\BLPSgip{\ref b.BLPS:gip/, hereinafter referred to as %
\def\BLPSgip{BLPS (\refbyear {BLPS:gip})}\BLPSgip}

\def\LPSmsf{\ref b.LPS:msf/}

\ifproofmode \relax \else\head{{\it Electron. J. Probab.}, {\bf 12}, Paper
54 (2007), 1454--1508.}
{Version of \versiondate}\fi
\vglue20pt

\def\firstheader{\eightpoint\ss\underbar{\raise2pt\line 
{{\it Electron.\ J.\ Probab.}, {\bf 12}, Paper 54 (2007), 1454--1508. 
     \hfil \frenchspacing Version of \versiondate}}}

\beginniceheadline

\title{Processes on Unimodular Random Networks}

\author{David Aldous and Russell Lyons}

\abstract{
We investigate unimodular random networks.
Our motivations include 
their characterization via reversibility of an associated random walk
and their similarities to unimodular quasi-transitive graphs.
We extend various theorems concerning random walks, percolation, spanning
forests, and amenability from the known context of unimodular
quasi-transitive graphs to the more general context of unimodular random
networks.
We give properties of a trace associated to unimodular random networks
with applications to stochastic comparison of
continuous-time random walk.
}

\bottomII{Primary
60C05. 
Secondary
60K99,  
 05C80. 
}
{Amenability, equivalence relations, infinite graphs, percolation,
quasi-transitive, random walks, transitivity, weak convergence,
reversibility, trace, stochastic comparison, spanning forests,
sofic groups.}
{Research partially supported by NSF grants DMS-0203062,
DMS-0231224, DMS-0406017, and DMS-0705518.}

\articletoc

\bsection{Introduction}{s.intro}

In the setting of infinite discrete graphs,
the property of being a Cayley graph of a group is a strong form of 
``spatial homogeneity":
many results not true for arbitrary graphs are true under this strong property.
As we shall soon explain, weaker regularity properties sufficient for many
results have been studied.  
In this paper, we turn to random graphs, investigating
a notion of
``statistical homogeneity" or ``spatial stationarity"
that we call
a {\it unimodular random rooted network}.
The root is merely a distinguished vertex of the network and the
probability measure is on a certain space of rooted networks.
In a precise sense, the root is ``equally likely" to be any vertex of the
network, even though we consider infinite networks.
We shall show that many results known for
deterministic graphs under previously-studied regularity conditions 
do indeed extend to
unimodular random rooted networks.

Thus, a probabilistic motivation for our investigations is the study of
stochastic processes under unimodularity.
A second motivation is combinatorial: One often asks for asymptotics of
enumeration or optimization problems on finite networks as the size of the
networks tend to infinity.
One can sometimes answer such questions with the aid of a suitable limiting
infinite object. A survey of this approach is given by \ref b.AS:obj/.
We call ``random weak limit" the type of limit one considers; it is the
limiting ``view" from a uniformly chosen vertex of the finite networks.
What limiting objects can arise this way?
It has been observed before that
the probabilistic objects of interest,
unimodular random rooted networks, contain all
the combinatorial objects of interest, random
weak limits of finite networks.
One open question is whether these two classes in fact coincide.
An affirmative answer would have many powerful consequences, as we shall
explain.

To motivate this by
analogy, recall a simple fact about stationary sequences
$\Seq{Y_i}_{i \in \Z}$
of random variables.
For each $n \geq 1$, let 
$\Seq{Y_{n,i}}_{1 \leq i \leq n}$
be arbitrary. Center it at a uniform index 
$U_n \in \{1,2,\ldots,n\}$ 
to get a bi-infinite sequence
$\Seq{Y_{n,U_n+i}}_{i \in \Z}$,
interpreted arbitrarily outside its natural range.  
If there is a weak limit 
$\Seq{Y_i}_{i \in \Z}$ as $n \to \infty$ 
of these randomly centered sequences, then the limit is
stationary, and conversely any stationary sequence can be obtained
trivially as such a limit.

By analogy, then,
given a finite graph, take a uniform random vertex as root.
Such a randomly rooted graph automatically has a certain property
(in short, if mass is redistributed in the graph, then the expected mass
that leaves the root is equal to the expected mass the arrives at the root)
and in \ref s.notation/,
we abstract this property as unimodularity.
It is then immediate that any infinite random rooted graph that is a limit
(in an appropriate sense that we call ``random weak limit")
of uniformly randomly rooted finite graphs will be unimodular, whereas
the above question asks
conversely whether any unimodular random rooted graph arises as a random
weak limit of some sequence of randomly rooted finite graphs.

Additional motivation for the definition arises from random walk considerations.
Given any random rooted graph, simple random walk induces a Markov chain on
rooted graphs.
Unimodularity of a probability measure $\rtd$ on rooted graphs
is equivalent to the property that a reversible stationary
distribution for this chain (properly interpreted)
is given by the root-degree biasing of $\rtd$,
just as on finite graphs, 
a stationary distribution for simple random walk is proportional to the
vertex degrees; see \ref s.extreme/.

Let us return now to the case of deterministic graphs.
An apparently minor relaxation of the Cayley graph property
is the ``transitive" property 
(that there is an automorphism taking any vertex to any other vertex).
By analogy with the shift-invariant interpretation of stationary sequences, one might expect every
transitive graph to fit into our set-up.
But this is false.
Substantial research over the last ten years has shown that
the most useful regularity condition is that of a
{\it unimodular transitive graph} (or, more generally, quasi-transitive).
Intuitively, this is
an unrooted transitive graph that can be given 
a random root in such a way that each vertex is equally
likely to be the root.
This notion is, of course, precise in itself for a finite graph.
To understand how this is extended to infinite graphs, and then to
unimodular random rooted graphs, consider a finite graph $\gh$ and a
function $f(x, y)$ of ordered pairs of vertices of $\gh$.
Think of $f(x, y)$ as an amount of mass that is sent from $x$ to $y$.
Then the total mass on the graph $\gh$ before transport equals the total
after, since mass is merely redistributed on the graph.
We shall view this alternatively as saying that for a randomly
uniformly chosen vertex, the expected mass it receives is equal to the
expected mass it sends out.
This, of course, depends crucially on choosing the vertex uniformly and,
indeed, characterizes the uniform measure among all probability measures on
the vertices.

Consider now an infinite transitive graph, $\gh$.
Since all vertices ``look the same", we could just fix one, $\bp$,
rather than try to choose one uniformly.
However, a mass transport function $f$ will not conserve the mass at $\bp$
without some assumption on $f$ to make up for the fact that $\bp$ is fixed.
Although it seems special at first, it turns out that a very useful
assumption is that $f$ is invariant under the diagonal action of
the automorphism group of $\gh$.
(For a finite graph that happened to have no automorphisms other than the
identity, this would be no restriction at all.)
This is still not enough to guarantee ``conservation of mass", i.e., that
$$
\sum_x f(\bp, x) = \sum_x f(x, \bp)
\,,
\label e.mtp1
$$
but it turns out that \ref e.mtp1/ does hold when the
automorphism group of $\gh$ is unimodular. Here, ``unimodular" is used in
its original sense that the group
admits a non-trivial Borel measure that is invariant under both left and
right multiplication by group elements.
We call $\gh$ itself unimodular in that case; see Sections
\briefref s.notation/ and \briefref s.related/ for more
on this concept.
The statement that \ref e.mtp1/ holds under these assumptions is called the
Mass-Transport Principle for $\gh$.
If $\gh$ is quasi-transitive, rather than transitive, we still have a
version of \ref e.mtp1/, but we can no longer consider only one fixed
vertex $\bp$.
Instead, each orbit of the automorphism group must have a representative
vertex.
Furthermore, it must be weighted ``proportionally to its frequency" among
vertices; see \ref t.fixedgraph/. 
This principle was introduced to the study of percolation by \ref
b.Hag:deptree/, then developed and exploited heavily by \BLPSgip.
Another way of stating it is that \ref e.mtp1/ holds in expectation 
when $\bp$ is chosen randomly by an appropriate probability measure.
If we think of $\bp$ as the root, then we arrive at the notion of
random rooted graphs, and the corresponding statement that \ref e.mtp1/
holds in expectation is a general form of the Mass-Transport Principle.
This general form was called the ``Intrinsic Mass-Transport Principle" by
\ref b.BS:rdl/.
We shall call a probability measure on rooted graphs unimodular precisely
when this general form of the Mass-Transport Principle holds.
We develop this in \ref s.notation/.

Thus, we can extend many results known
for unimodular quasi-transitive graphs to our new setting of
unimodular random rooted graphs, as noted by \ref b.BS:rdl/.
As a bonus, our set-up allows the treatment of
quasi-transitive graphs to be precisely parallel to that of
transitive graphs, with no additional notation or thought needed, which had
not always been the case previously.

To state results in their natural generality, as well as for technical
convenience, we shall work in the setting of networks,
which are just graphs with ``marks" (labels) on edges and vertices.
Mainly, this paper is organized to progress from the most general to the
most specific models.
An exception is made in \ref s.related/, where we discuss random networks
on fixed underlying graphs.
This will not only help to understand and motivate the general setting, but
also will be useful in deriving consequences of our general results.

\ref s.extreme/ elaborates the comment above about reversible stationary
distributions for random walk,  discussing extremality and invariant
$\sigma$-fields, speed of random walk, and continuous-time random walk and
their explosions.
\ref s.trace/ discusses a trace associated to unimodular random networks
and comparison of return probabilities of different continuous-time random
walks, which partially answers a question of Fontes and Mathieu.
We then write out the extensions to unimodular random rooted graphs of results
known for fixed graphs in the context of
percolation (\ref s.extend/),
spanning forests (\ref s.span/)
and amenability (\ref s.amen/).
These extensions are in most (though not all) cases straightforward.
Nevertheless, we think it is useful to list these extensions in the order
they need to be proved so that others need not check the entire (sometimes
long) proofs or chains of theorems from a variety of papers.
Furthermore, we were required to find several essentially new results along
the way.

In order to appreciate the scope of our results, we list many examples of
unimodular probability measures in \ref s.ex/.
In particular, there is a significant and important overlap between our
theory and the theory of graphings of measure-preserving equivalence
relations. This overlap
is well known among a few specialists, but
deserves to be made more explicit. We do that here in \ref x.equivalence/.

Among several open problems, we spotlight a special case of
\ref q.coupling/: Suppose we
are given a partial order on the mark space and two unimodular probability
measures, one stochastically dominating the other. That is, there is
a monotone coupling of the two unimodular distributions that puts the
networks on the same graphs, but has higher marks for the second network
than for the first. Does this
imply the existence of a {\it unimodular\/} monotone coupling?
A positive answer would be of great benefit in a variety of ways.

Another especially important open question is \ref q.TII/, whether every
unimodular probability measure is a limit of uniformly rooted finite
networks.
For example, in the case that the random rooted infinite network is just a
Cayley graph (rooted, say, at the identity) with the edges marked by the
generators, a positive answer to this question on finite approximation would
answer a question of \ref b.Weiss:sofic/, by showing that all finitely
generated groups are ``sofic", although this is contrary to the belief
expressed by \ref b.Weiss:sofic/.
(Sofic groups were introduced, with a different definition, by \ref
b.Gromov:endo/; see \ref b.Elek:sofic/ for a proof that the definitions are
equivalent.)
This would establish several conjectures, since they are known to hold for
sofic groups: the direct finiteness
conjecture of \ref b.Kaplansky/ on group algebras
(see \ref b.Elek:sofic/), a
conjecture of \ref b.Gottschalk/ on ``surjunctive" groups in topological
dynamics
(see \ref b.Gromov:endo/), the Determinant
Conjecture on Fuglede-Kadison determinants
(see \ref b.Elek:hyper/), and Connes's (\refbyear{Connes})
Embedding Conjecture
for group von Neumann algebras (see \ref b.Elek:hyper/).
The Determinant Conjecture in turn implies the Approximation Conjecture of
\ref b.Schick/ and the Conjecture of Homotopy Invariance of $L^2$-Torsion
due to \ref b.Luck:3man/; see Chap.~13 of \ref b.Luck:book/ for these
implications and more information.
\ref b.Weiss:sofic/ gave another proof
of Gottschalk's conjecture for sofic groups.
One may easily extend that proof to show a form of Gottschalk's conjecture 
for all quasi-transitive unimodular
graphs that are limits of finite graphs, but there are easy
counterexamples for general transitive graphs.

Further discussion of the
question on approximation by finite networks is given in
\ref s.finite/.
A positive answer would
provide solid support for the intuition that the root of a
unimodular random rooted network is equally likely to be any vertex.
\ref s.finite/ also contains some variations that would result from a
positive answer and some
additional consequences for deterministic graphs.

The notion of weak convergence of rooted locally finite graphs or networks  
(needed to make sense of convergence of randomly rooted
finite graphs to a limit infinite graph)
has arisen before in several different contexts.
Of course, the special case where the limit network is a Cayley diagram
was introduced by \ref b.Gromov:endo/ and \ref b.Weiss:sofic/.
In the other cases, the limits provide examples of unimodular
random rooted graphs.
\ref b.Aldous:fringe/
gives many examples of models of random finite trees
which have an infinite-tree limit 
(and one such example, the limit of uniform random labeled trees
being what is now called the 
Poisson-Galton-Watson tree,
${\PGW}^\infty(1)$, goes back to 
\ref b.Grimmett:ltrees/).
The idea that random weak limits of finite planar graphs of
uniformly bounded degree provide an interesting class of infinite planar
graphs was developed by \ref b.BS:rdl/, who showed that random walk on
almost any such limit graph is recurrent.
(Thus, such graphs do not include regular trees or hyperbolic graphs,
other than trivial examples like $\Z$.)
A specialization to random weak limits of
plane triangulations was studied in more detail in
interesting recent work of \ref b.AngelSchramm/ and \ref b.Angel/.

\ref x.PWIT/ describes an infinite-degree tree, arising as a limit
of weighted finite complete graphs. 
This example provides an interface between our setting and 
related ideas of 
``local weak convergence"
and ``the objective method in the probabilistic analysis of algorithms".
A prototype is that the distribution of $n$ random points in a square of area $n$
converges in a natural sense as $n \to \infty$ to the distribution
of a Poisson point process on the plane of unit intensity.
One can ask whether solutions of combinatorial optimization
problems over the $n$ random points 
(minimum spanning tree, minimum matching, traveling salesman problem)
converge to limits that are the solutions of analogous
optimization problems over the Poisson point process in the whole plane.
\ref x.PWIT/ can be regarded as a mean-field analogue of random
points in the plane, and $n \to \infty$ limits of solutions of combinatorial optimization problems
within this model have been studied using the non-rigorous
cavity method from statistical physics.  
\ref b.APercus:SU/ illustrate what can be done by non-rigorous means,
while \ref b.AS:obj/ survey introductory rigorous theory.

The reader may find it helpful to keep in mind one additional example,
a unimodular version of family trees of Galton-Watson branching
processes; see also \ref x.fdd/.

\procl x.AGW \procname{Unimodular Galton-Watson}
Let $\Seq{p_k \st k \ge 0}$ be a probability distribution on $\N$.
Take two independent Galton-Watson trees with offspring distribution
$\Seq{p_k}$, each starting with one particle, the root,
and join them by a new edge
whose endpoints are their roots. Root the new tree
at the root of the first Galton-Watson
tree. This is {\bf augmented Galton-Watson}
measure, \AGW. (If $p_0 \ne 0$, then we have the
additional options to condition on either
non-extinction or extinction of the joined trees.) Now bias by the
reciprocal of the degree of the root to get {\bf unimodular Galton-Watson}
measure, $\UGW$.
In different language, \ref b.LPP:GW/ proved that this measure,
$\UGW$, is unimodular.
Note that
the mean degree of the root is
$$
\expdeg(\UGW) = \left(\sum_{k \ge 0} {p_k \over k+1} \right)^{-1}
\,.
\label e.degAGW
$$
\endprocl

\bsection{Definitions and Basics}{s.notation}

We denote a (multi-)graph $\gh$ with vertex set $\verts$ and undirected
edge set $\edges$ by $\gh = (\verts, \edges)$.
When there is more than one graph under discussion, we write $\verts(\gh)$ or
$\edges(\gh)$ to avoid ambiguity.
We denote the degree of a vertex $x$ in a graph $\gh$ by $\deg_\gh(x)$.
{\bf Simple random walk} on $\gh$ is the Markov chain whose state space is
$\verts$ and whose transition probability from $x$ to $y$ equals the number of
edges joining $x$ to $y$ divided by $\deg_\gh(x)$.

A {\bf network} is a (multi-)graph $\gh = (\vertex, \edge)$ together with a
complete separable metric space $\marks$ called the {\bf mark space} and
maps from $\vertex$ and $\edge$ to $\marks$. Images in $\marks$ are called
{\bf marks}.
Each edge is given two marks, one associated to (``at") each of its endpoints.
The only assumption on degrees is that they are finite.
We shall usually assume that $\marks$ is Baire space $\N^\N$, since every
uncountable complete separable metric space is Borel isomorphic to Baire
space by Kuratowski's theorem (Theorem 15.10 of \ref b.Royden/).
We generally omit mention of the mark maps from our notation for
networks when we do not need them.
For convenience, we consider graphs as special cases of networks in which
all marks are equal to some fixed mark.

We now define ends in graphs.
In the special case of a tree, an infinite path that starts at any vertex
and does not backtrack is called a {\bf ray}. Two rays are {\bf equivalent}
if they have
infinitely many vertices in common. An equivalence class of rays is
called an {\bf end}.
In a general infinite graph, $\gh$, an {\bf end} of $\gh$ is an equivalence
class of infinite simple paths in $\gh$, where two paths
are equivalent if for every finite $K\subset \vertex(\gh)$,
there is a connected component of $\gh\setminus K$ that intersects
both paths.

Let $\gh$ be a graph.  For a subgraph $H$, let its {\bf (internal) vertex
boundary} $\bdv H$ be the set of vertices of $H$ that are adjacent to some
vertex not in $H$.
We say that $\gh$ is {\bf (vertex) amenable}
if there exists a sequence of subsets $H_n\subset \vertex(\gh)$
with
$$
\lim_{n \to\infty} {|\bdv H_n| \over |\verts(H_n)|} = 0\,,
$$
where $|\cbuldot|$ denotes cardinality.
Such a sequence is called a {\bf F\o{}lner sequence}.
A finitely generated group is {\bf amenable} if its Cayley graph is amenable.
For example, every finitely generated abelian group is amenable.
For more on amenability of graphs and groups, see \BLPSgip.

A {\bf homomorphism} $\varphi: \gh_1 \to \gh_2$ from one
graph $\gh_1=(\verts_1,\edges_1)$ to another $\gh_2=(\verts_2,\edges_2)$
is a pair of maps $\varphi_\verts:\verts_1\to\verts_2$ and
$\varphi_\edges:\edges_1\to\edges_2$
such that $\varphi_\verts$ maps the endpoints of $e$ to the endpoints
of $\varphi_\edges(e)$ for every edge $e \in \edges_1$.
When both maps $\varphi_\verts:\verts_1\to\verts_2$ and
$\varphi_\edges:\edges_1\to\edges_2$ are bijections, then $\varphi$ is
called an {\bf isomorphism}.
When $\gh_1 = \gh_2$, an isomorphism is called
an {\bf automorphism}.
The set of all automorphisms of $\gh$ forms a group under composition,
denoted by $\Aut(\gh)$.
The action of a group $\gp$ on a graph $\gh$ by automorphisms is said to be {\bf
transitive} if there is only one $\gp$-orbit in $\verts(\gh)$ and to be {\bf
quasi-transitive} if there are only finitely many orbits in $\verts(\gh)$.
A graph $\gh$ is {\bf transitive} or {\bf quasi-transitive} according as
whether the corresponding action of $\Aut(\gh)$ is.
For example, every Cayley graph is transitive.
All the same terms are applied to networks when the maps in question preserve
the marks on vertices and edges.

A locally compact group is called {\bf unimodular} if its
left Haar measure is also right invariant.
In particular, every discrete countable group is unimodular.
We call a graph $\gh$ {\bf unimodular} if $\Aut(\gh)$ is unimodular, where
$\Aut(\gh)$ is given the weak topology generated by its action on $\gh$.
Every Cayley graph and, as \ref b.SoardiWoess/ and \ref b.Salvatori/ proved,
every quasi-transitive amenable graph is unimodular.
See \ref s.related/ and \BLPSgip\ for more details on unimodular graphs.

A {\bf rooted network} $(\gh, \bp)$ is a network $\gh$ with a distinguished
vertex $\bp$ of $\gh$, called the {\bf root}.
A {\bf rooted isomorphism} of rooted networks is an isomorphism of the
underlying networks that takes the root of one to the root of the other.
We generally do not distinguish between a rooted network and its
isomorphism class.
When needed, however, we use the following notation to make these distinctions:
$\gh$ will denote a graph, $\overline \gh$ will denote a network with underlying
graph $\gh$, and $[\overline \gh, \bp]$ will denote the class of rooted networks
that are rooted-isomorphic to $(\overline \gh, \bp)$.
We shall use the following notion introduced (in slightly different
generalities) by \ref b.BS:rdl/ and \ref b.AS:obj/.
%
Let $\GG_*$ denote the set of rooted isomorphism classes of rooted
{\it connected\/} locally finite networks.
Define a metric on $\GG_*$ by letting the distance between $(G_1, o_1)$ and
$(G_2, o_2)$ be $1/(1+\alpha)$, where $\alpha$ is
the supremum of those $r > 0$ such that
there is some rooted isomorphism of
the balls of (graph-distance) radius $\flr{r}$ around the roots of $G_i$ 
such that each pair of corresponding marks has distance less than $1/r$.
It is clear that $\GG_*$ is separable and complete in this metric.
For probability measures $\rtd$, $\rtd_n$ on $\GG_*$, we write $\rtd_n \cd
\rtd$ when $\rtd_n$ converges weakly with respect to this metric.

For a probability measure $\rtd$ on rooted networks, write $\expdeg(\rtd)$
for the {\it expectation\/} of the degree of the root with respect to
$\rtd$.
In the theory of measured equivalence relations (\ref x.equivalence/),
this is twice the {\bf
cost} of the graphing associated to $\rtd$.
Also, by the {\bf degree} of $\rtd$ we mean the {\it distribution\/} of
the degree of the root under $\rtd$.

For a locally finite connected rooted network, there is a canonical choice
of a rooted network in its rooted-isomorphism class.
More specifically, there is a continuous map $f$ from $\GG_*$ to the space of
networks on $\N$ rooted at 0 such that $f\big([\gh, \bp]\big) \in [\gh,
\bp]$ for all $[\gh, \bp] \in \GG_*$.
To specify this, consider the following total ordering on rooted networks
with vertex set $\N$ and root $0$.
First, total order $\N \times \N$ by the lexicographic order:
$(i_1, j_1) \prec (i_2, j_2)$ if
either $i_1 < i_2$ or $i_1 = i_2$ and $j_1 < j_2$.
Second, the lexicographic order on Baire space $\marks$ is also a total order.
We consider networks on $\N$ rooted at $0$.
Define a total order on such networks as follows.
Regard the edges as oriented for purposes of identifying the edges with $\N
\times \N$; the mark at $i$ of an edge between $i$ and $j$ will be considered
as the mark of the oriented edge $(i, j)$.
Suppose we are given a pair of networks on $\N$ rooted at $0$.
If they do not have the same edge sets, then the network that contains the
smallest edge in their symmetric difference is deemed to be the smaller
network.
If they do have the same edge sets, but not all the vertex marks are the same,
then the network that contains the vertex with the smaller mark on the
least vertex where they differ is deemed the smaller network.
If the networks have the same edge sets and the same vertex marks, but not all
the edge marks are the same, then
the network that contains the oriented edge with the smaller mark on the
least oriented edge where they differ is deemed the smaller network.
Otherwise, the networks are identical.

We claim that the rooted-isomorphism class of each locally finite connected
network contains a unique smallest rooted network on $\N$ in the above
ordering.
This is its {\bf canonical} representative.
To prove our claim, 
given a locally finite, connected, rooted network $G$ and $r \ge 1$, let
$\Giso_r$ be the class of networks on $\N$ with root $0$ that are
rooted-isomorphic to $G$ and whose vertices within distance $r$ of $0$ form
an interval, $[0, N_r]$.
Let $\Gisop_r$ be the subset of $\Giso_r$ such that the network induced on
$[0, N_r]$ is minimal for $\prec$ (there are only finitely many
possibilities for the induced network, so there is a unique minimum induced
network).
Then $\Gisop_r \supseteq \Gisop_{r+1}$ for all $r$ by the definition of
$\prec$.
Hence, there is a unique element $H \in \bigcap_{r=1}^\infty \Gisop_r$: the
network of $H$ induced on $[0, N_r]$ is determined by $\Gisop_r$.
This network $H$ is the desired canonical representative of $G$.

For a (possibly disconnected)
network $\gh$ and a vertex $x \in \verts(\gh)$, write $\gh_x$ for the
connected component of $x$ in $\gh$.
If $\gh$ is a network with probability distribution $\rtd$ on its
vertices, then $\rtd$ induces naturally a distribution on $\GG_*$,
which we also denote by $\rtd$; namely, the probability of $(\gh_x, x)$ is
$\rtd(x)$.
More precisely, $\rtd\big([\gh_x, x]\big) := \sum \big\{ \rtd(y) \st y \in
\verts(\gh), \ (\gh_y, y) \in [\gh_x, x] \big\}$.
For a finite network $\gh$, let $U(\gh)$ denote the distribution on $\GG_*$
obtained this way by choosing a uniform random vertex of $\gh$ as root.
Suppose that $\gh_n$ are finite networks and that $\rtd$ is a
probability measure on $\GG_*$.
We say the {\bf random weak limit} of $\gh_n$ is $\rtd$ if $U(\gh_n) \cd \rtd$.
If $\rtd\Big(\big\{[\gh, \bp]\big\}\Big) = 1$ for a fixed transitive
network $\gh$ (and (any) $\bp \in \vertex(\gh)$), then we say that the {\bf
random weak limit} of $\gh_n$ is $\gh$.

As usual, call a collection $\collec$
of probability measures on $\GG_*$
{\bf tight} if for each $\epsilon > 0$, there is a
compact set $\K \subset \GG_*$ such that $\rtd(\K) > 1 - \epsilon$ for all
$\rtd \in \collec$.
Because $\GG_*$ is complete, any tight collection
has a subsequence that possesses a weak limit.

The class of probability measures
$\rtd$ that arise as random weak limits of finite networks is contained in
the class of unimodular $\rtd$, which we now define.
Similarly to the space $\GG_*$, we define the space $\gtwo$ of isomorphism
classes of locally
finite connected networks with an ordered pair of distinguished vertices
and the natural topology thereon.
We shall write a function $f$ on $\gtwo$ as $f(\gh, x, y)$.

\procl d.unimodular 
Let $\rtd$ be a probability measure on $\GG_*$.
We call $\rtd$ {\bf unimodular} if it obeys the {\bf Mass-Transport
Principle}:
For all Borel
$f : \gtwo \to [0, \infty]$,
we have
$$
\int \sum_{x \in \vertex(\gh)} f(\gh, \bp, x) \,d\rtd\big([\gh, \bp]\big)
=
\int \sum_{x \in \vertex(\gh)} f(\gh, x, \bp) \,d\rtd\big([\gh, \bp]\big)
\,.
\label e.mtpgen
$$
Let $\uni$ denote the set of unimodular Borel probability measures on $\GG_*$.
\endprocl

Note that to define the sums that occur here, we choose a specific network
from its rooted-isomorphism class, but which one we choose makes no
difference when the sums are computed.
We sometimes call $f(\gh, x, y)$ the amount of ``mass" sent from $x$ to $y$.
The motivation 
for the name ``unimodular" is two fold: One is the extension of the concept
of unimodular automorphism groups of networks. The second is
that the Mass-Transport
Principle expresses the equality of two measures on $\gtwo$
associated to $\rtd$, the ``left" measure $\rtd_{\rm L}$ defined by
$$
\int_\gtwo f \,d\rtd_{\rm L}
:=
\int_{\GG_*} \sum_{x \in \vertex(\gh)} f(\gh, \bp, x) \,d\rtd\big([\gh,
\bp]\big)
$$ 
and the ``right" measure $\rtd_{\rm R}$ defined by
$$
\int_\gtwo f \,d\rtd_{\rm R}
:=
\int_{\GG_*} \sum_{x \in \vertex(\gh)} f(\gh, x, \bp) \,d\rtd\big([\gh,
\bp]\big)
\,.
$$
Thus, $\rtd$ is unimodular iff $\rtd_{\rm L} = \rtd_{\rm R}$, which can
also be expressed by saying that the left measure is absolutely continuous
with respect to the right measure and has Radon-Nikod\'ym derivative 1.

It is easy to see that any $\rtd$ that is a random weak limit of finite
networks is unimodular, as observed by \ref b.BS:rdl/, who introduced this
general form of the Mass-Transport Principle under the name 
``intrinsic Mass-Transport Principle".
The converse is open.

A special form of the Mass-Transport Principle was considered, in different
language, by \ref b.AS:obj/.
Namely, they defined $\rtd$ to be {\bf involution invariant} if \ref
e.mtpgen/ holds for those $f$ supported on $(\gh, x, y)$ with $x \sim y$.
In fact, the Mass-Transport Principle holds for general $f$ if it holds for
these special $f$:

\procl p.invinv
A measure is involution invariant iff it is unimodular.
\endprocl

\proof
Let $\rtd$ be involution invariant.
The idea is to send the mass from $x$ to $y$ by single steps, equally
spread among the shortest paths from $x$ to $y$.
For the proof, we may assume that $f(\gh, x, y) = 0$ unless $x$ and $y$ are
at a fixed distance, say $k$, from each other, since any $f$ is a sum of
such $f$.
Now write $L(\gh, x, y)$ for the set of paths of length $k$ from $x$ to $y$.
Let $n_j(\gh, x, y; z, w)$ be the number of paths in $L(\gh, x, y)$ such
that the $j$th edge goes from $z$ to $w$.
Define $f_j(\gh, z, w)$ for $1 \le j \le k$ and $z, w \in \verts(\gh)$
by
$$
f_j(\gh, z, w) :=
\sum_{x, y \in \verts(\gh)} {f(\gh, x, y) n_j(\gh, x, y; z, w) \over
|L(\gh, x, y)|}
\,.
$$
Then $f_j(\gh, z, w) = 0$ unless $z \sim w$. Furthermore,
$f_j(\gh, z, w) := f_j(\gh', z', w')$ if
$(\gh, z, w)$ is isomorphic to $(\gh', z', w')$.
Thus, $f_j$ is well defined and Borel on $\gtwo$,
whence involution invariance gives us
$$
\int \sum_{x \in \vertex(\gh)} f_j(\gh, \bp, x) \,d\rtd(\gh, \bp)
=
\int \sum_{x \in \vertex(\gh)} f_j(\gh, x, \bp) \,d\rtd(\gh, \bp)
\,.
$$
On the other hand,
$$
\sum_{x \in \vertex(\gh)} f(\gh, \bp, x)
=
\sum_{x \in \vertex(\gh)} f_1(\gh, \bp, x)
\,,
$$
$$
\sum_{x \in \vertex(\gh)} f(\gh, x, \bp)
=
\sum_{x \in \vertex(\gh)} f_k(\gh, x, \bp)
\,,
$$
and for $1 \le j < k$, we have
$$
\sum_{x \in \vertex(\gh)} f_j(\gh, x, \bp)
=
\sum_{x \in \vertex(\gh)} f_{j+1}(\gh, \bp, x)
\,.
$$
Combining this string of equalities yields the desired equation for
$f$.
\Qed

\comment{
If $\gh$ is a network, $x \in \verts(\gh)$, and $e \in \edges(\gh)$, then
the {\bf isomorphism class} of $(\gh, x, e)$ is the class $\{(\gpe \gh, \gpe
x, \gpe e) \st \gpe \hbox{ is an isomorphism of }\gh\}$.
Let $\gve$ be the set of such isomorphism classes with the natural
topology.
Given a rooted network $(\gh, x)$ and an edge $e$ incident to $x$, define
the involution $\iota(\gh, x, e) := (\gh, y, e)$, where $y$ is the other
endpoint of $e$.
Given a probability measure $\rtd$ on rooted networks, define the probability
measure $\widehat\rtd$ to be the law of the isomorphism class of
$(\gh, x, e)$, where $(\gh, x)$ is chosen
according to $\rtd$ and $e$ is then chosen uniformly among the edges
incident to $x$.
Thus,
$$
\widehat\rtd(\A) :=
\int_{\GG_*} {1 \over \deg_\gh x} |\{e \sim x \st (\gh, x, e) \in \A\}|
\,d\rtd\big([\gh, x]\big)
$$
for Borel sets $\A \subseteq \gve$.
(Note that we actually choose $[\gh, x]$ according to $\rtd$, then choose
$(\gh, x) \in [\gh, x]$ arbitrarily; this latter choice has no effect on
the number of $e \sim x$ such that $(\gh, x, e) \in \A$.)
Also, define $\widetilde \rtd$ to be the (non-probability) measure that is the
result of biasing $\widehat\rtd$ by the degree of the root; that is, the
Radon-Nikod\'ym derivative of $\widetilde \rtd$ with respect to $\widehat\rtd$
at the isomorphism class of $(\gh, x, e)$ is $\deg_\gh(x)$:
$$
\widetilde\rtd(\A) :=
\int_{\GG_*} |\{e \sim x \st (\gh, x, e) \in \A\}|
\,d\rtd\big([\gh, x]\big)
$$
for Borel sets $\A \subseteq \gve$.
(If the expected degree of the root is finite, one could obtain a
probability measure from $\widetilde \rtd$ by dividing by the expected degree;
but this is not always the case.)
The involution $\iota$ induces a pushforward map $\widetilde\rtd \mapsto
\widetilde\rtd \circ \iota^{-1}$.
We say that $\rtd$ is {\bf unimodular} or {\bf involution invariant} if
$\widetilde\rtd \circ \iota^{-1}= \widetilde\rtd$.
\comment{This definition works only for locally finite networks. We must resort
to the definition by \ref b.AS:obj/ in the general case. The latter seems much
harder for probabilists to comprehend.}
}%


Occasionally one uses the Mass-Transport Principle for functions $f$
that are not nonnegative. It is easy to see that this use is justified
when
$$
\int \sum_{x \in \vertex(\gh)} |f(\gh, \bp, x)| \,d\rtd(\gh, \bp)
< \infty
\,.
$$

As noted by Oded Schramm (personal communication, 2004), unimodularity can
be defined for probability measures on other structures, such as
hypergraphs, while involution invariance is limited to graphs (or networks
on graphs).

\comment{
Suppose that $\gh$ is an infinite quasi-transitive amenable connected network.
Choose one element $o_i$ from each vertex orbit.
It is shown in \BLPSgip, Proposition 3.6, that there is a probability
measure $\rho$ on the set $ \{ o_i \}$ such that for any F\o{}lner sequence
$H_n$, the relative frequency of vertices in $H_n$ that are in the same orbit
as $o_i$ converges to $\rho(o_i)$. We call this measure $\rho$ the {\bf
natural frequency distribution} of $\gh$.
}%

We shall sometimes use the following property of marks.
Intuitively, it says that each vertex has positive probability to be the
root.

\procl l.unmark \procname{Everything Shows at the Root}
Suppose that $\rtd$ is a unimodular probability measure on $\GG_*$.
Let $\mk_0$ be a fixed mark and $\marks_0$ be a fixed 
Borel set of marks.
If the mark of the root is a.s.\ $\mk_0$, then the mark of every vertex is
a.s.\ $\mk_0$.
If every edge incident to the root a.s.\
has its edge mark at the root in $\marks_0$,
then all edge marks a.s.\ belong to $\marks_0$.
\endprocl

\proof
In the first case, each vertex
sends unit mass to each vertex with a mark different from $\mk_0$.
The expected mass received at the root is
zero. Hence the expected mass sent is 0.
The second case is a consequence of the first, where we put the mark $\mk_0$
at a vertex when all the edge marks at that vertex lie in $\marks_0$.
\Qed

When we discuss percolation in \ref s.extend/, we shall find it crucial
that we have a unimodular coupling of the various measures (given by the
standard coupling of Bernoulli percolation in this case).
It would also be very useful to have unimodular couplings in more general
settings.
We now discuss what we mean.

Suppose that $\rln \subseteq \marks \times \marks$ is a closed set, which we
think of as a binary relation such as the lexicographic order on Baire space.
Given two measures $\rtd_1, \rtd_2 \in \uni$, say that
$\rtd_1$ is {\bf $\rln$-related} to $\rtd_2$ if there is a probability measure
$\nu$, called an {\bf $\rln$-coupling of $\rtd_1$ to $\rtd_2$},
on rooted networks with mark space
$\marks \times \marks$ such that $\nu$ is concentrated on networks all of
whose marks lie in $\rln$ and whose marginal given by taking the $i$th
coordinate of each mark is $\rtd_i$ for $i = 1, 2$.
In particular, $\rtd_1$ and $\rtd_2$ can be coupled to have the same
underlying rooted graphs.
\comment{
Call $\rln$ {\bf extra separable} if there is a sequence $\Seq{\mk_l}$ and two
Borel maps $\rlnmap_i : \marks \to 2^{\Seq{\mk_l}}$ ($i = 1, 2$)
such that for all $(\mk, \mk') \in \rln$, for all $\zeta_1 \in
\rlnmap_1(\mk)$, and for all $\zeta_2 \in \rlnmap_2(\mk')$, we have $(\zeta_1,
\zeta_2) \in \rln$.
For example, this is the case with the relation $ \le $ on $\R$, where, e.g.,
we may take $\Seq{\mk_l} := \Q$, $\rlnmap_1(\mk) := (-\infty, \mk) \cap \Q$,
and $\rlnmap_2(\mk) := (\mk, \infty) \cap \Q$.
}%

It would be very useful to have a positive answer to the following question.
Some uses are apparent in \ref s.trace/ and in \ref s.finite/,
while others appear in \ref b.Lyons:est/ and are hinted at elsewhere.

\procl q.coupling \procname{Unimodular Coupling}
Let $\rln \subseteq \marks \times \marks$ be a closed set.
If $\rtd_1, \rtd_2 \in \uni$ and $\rtd_1$ is
$\rln$-related to $\rtd_2$, is there then a unimodular $\rln$-coupling of
$\rtd_1$ to $\rtd_2$?
\endprocl

The case where $\rtd_i$ are amenable is established affirmatively in \ref
p.amen-average/.
However, the case where $\rtd_i$ are supported by a fixed underlying
non-amenable Cayley graph is open even when the marks take only two values.
Here is a family of examples to illustrate what we do not know:

\procl q.end-couple
Let $T$ be the Cayley graph of $\Z_2 * \Z_2 * \Z_2$ with respect to the
generators $a$, $b$, $c$, which are all involutions.
We label the edges with the generators.
Fix three Borel symmetric
functions $f_a$, $f_b$, $f_c$ from $[0, 1]^2$ to $[0, 1]$.
Also, fix an end $\xi$ of $T$.
Let $U(e)$ be i.i.d.\ Uniform$[0, 1]$ random variables indexed by the edges
$e$ of $T$.
For each edge $e$, let $I_e$ be the two edges adjacent to $e$ that lead
farther from $\xi$ and let $J_e$ be the two other edges that are adjacent to
$e$.
Let $L(e)$ denote the Cayley label of $e$, i.e., $a$, $b$, or $c$.
For an edge $e$ and a pair of edges $\{e_1, e_2\}$, write $f\big(e, \{e_1,
e_2\}\big) := f_{L(e)}\big(U(e_1), U(e_2)\big)$.
Define $X(e) := f(e, I_e)$ and $Y(e) := \max \big\{f(e, I_e), f(e, J_e)\big\}$.
Let $\nu$ be the law of $(X, Y)$.
Let $\rtd_1$ be the law of $X$ and $\rtd_2$ be the law of $Y$.
We use the same notation for the measures in $\uni$ given by rooting $T$ at
the vertex corresponding to the identity of the group.
Let $\rln$ be $\le$ on $[0, 1] \times [0, 1]$.
Since $X(e) \le Y(e)$ for all $e$, $\nu$ is an $\rln$-coupling of $\rtd_1$ to
$\rtd_2$.
In addition, $\rtd_2$ is clearly $\Aut(T)$-invariant (recall that the edges
are labeled), while
the same holds for $\rtd_1$ since it is an i.i.d.\ measure.
Thus, $\rtd_i$ are both unimodular for $i = 1, 2$.
On the other hand, $\nu$ is not $\Aut(T)$-invariant except in the trivial
case that the functions $f_a$, $f_b$, and $f_c$ are all constant.
Is there an invariant $\rln$-coupling of $\rtd_1$ to $\rtd_2$?
In other words, is there a unimodular $\rln$-coupling of $\rtd_1$
to $\rtd_2$?
\endprocl

Another example concerns monotone coupling of the wired and free uniform
spanning forests (whose definitions are given below in \ref s.span/).
This question was raised in \BLPSusf; a partial answer was given by \ref
b.Bowen:coupling/.
This is not the only interesting situation involving graph inclusion.
To be more precise about this relation,
for a map $\psi : \marks \to \marks$ and a network $G$, let $\psi(G)$
denote the network obtained from $G$ by replacing each mark with its image
under $\psi$.
Given a Borel subset $\marks_0 \subseteq \marks$ and a network $\gh$,
call the subnetwork consisting of those edges both of whose edge marks lie in
$\marks_0$ the {\bf $\marks_0$-open subnetwork} of $\gh$.
If $\rtd$ and $\rtd'$ are two probability measures on rooted networks, let
us say that $\rtd$ is {\bf edge dominated} by $\rtd'$ if there exists a
measure $\nu$ on $\GG_*$, a Borel subset $\marks_0 \subseteq \marks$, and
Borel functions $\psi, \psi': \marks \to \marks$ such that if $(\gh', \bp)$
denotes a network with law $\nu$ and $(\gh, \bp)$ the component of $\bp$ in
the $\marks_0$-open subnetwork,
then $\big(\psi(\gh), \bp\big)$ has law $\rtd$ and
$\big(\psi'(\gh'), \bp\big)$ has law $\rtd'$.
If the measure $\nu$ can be chosen to be unimodular, then we say that
$\rtd$ is {\bf unimodularly edge dominated} by $\rtd'$.
As a special case of \ref q.coupling/, we do not know whether the existence
of such a measure $\nu$ that is not unimodular implies the existence of
$\nu$ that is unimodular when $\rtd$ and $\rtd'$ are both unimodular
themselves.

\bsection{Fixed Underlying Graphs}{s.related}

Before we study general unimodular probability measures, it is useful to
examine the relationship between unimodularity in the classical sense for
graphs and unimodularity in the sense investigated here for random rooted
network classes.

Given a graph $\gh$ and a vertex $x \in \verts(\gh)$, write $\Stab(x) := \{
\gpe \in \Aut(\gh) \st \gpe x = x\}$ for
the {\bf stabilizer subgroup} of $x$.
Also, write $[x] := \Aut(\gh) x$ for the orbit of $x$.
Recall the following principle from \BLPSgip:

\proclaim Mass-Transport Principle. If $\gh= (\vertex, \edge)$ is any graph,
$f : \vertex \times \vertex \to [0,
\infty]$ is invariant under the diagonal action of $\Aut(\gh)$, and $o, o' \in
\vertex$, then
$$
\sum_{z \in [o']} f(o, z) |\Stab(o')| = \sum_{y \in [o]} f(y, o') |\Stab(y)| \,.
$$

Here, $|\cbuldot|$ denotes Haar measure on $\Aut(\gh)$, although we
continue to use this notation for cardinality as well.
Since $\Stab(x)$ is compact and open, $0 < |\Stab(x)| < \infty$.
As shown in \ref b.Schlichting/ and \ref b.Trofimov/,
$$
|\Stab(x)y| / |\Stab(y)x| = |\Stab(x)| / |\Stab(y)|
\,.
\label e.schl-trof
$$
It follows easily that $\gh$ is unimodular iff
$$
|\Stab(x)y| = |\Stab(y)x|
\label e.unimod-schl-trof
$$
whenever $x$ and $y$ are in the same orbit.

\procl t.fixedgraph \procname{Unimodular Fixed Graphs}
Let $G$ be a fixed connected graph. Then
$G$ has a random root that gives a unimodular
measure iff $G$ is a unimodular graph with
$$
c := \sum_i |\Stab(o_i)|^{-1} <\infty
\,,
\label e.alm-trn
$$
where $\{ o_i\} $ is a complete orbit
section.
In this case, there is only one such measure $\rtd$ on random rooted graphs
from $G$ and it satisfies
$$
\rtd([G, x]) = c^{-1} |\Stab(x)|^{-1}
\label e.mustab
$$
for every $x \in \verts(G)$.
\endprocl

Of course, a similar statement holds for fixed networks.
An example of a graph satisfying \ref e.alm-trn/, but that is not
quasi-transitive, is obtained from the random weak limit of balls in a
3-regular tree.
That is, let $\vertex := \N \times \N$.
Join $(m, n)$ by edges to each of $(2m, n-1)$ and $(2m+1, n-1)$ for $n \ge
1$.
The result is a tree with only one end and $\big|\Stab\big((m, n)\big)\big|
= 2^n$.

\proof
Suppose first that $G$ is unimodular and that $c <\infty$.
Define $\rtd$ by
$$
\all i \rtd([G, o_i]) := c^{-1} |\Stab(o_i)|^{-1}
\,.
$$
To show that $\rtd$ is unimodular, let
$f : \gtwo \to [0, \infty]$ be Borel.
Since we are concerned only with the graph $\gh$, we shall write $f$
instead as a function
$f : \vertex \times \vertex \to [0,
\infty]$ that is $\Aut(\gh)$-invariant.
Then
$$\eqaln{
\int \sum_x f(o, x) \,d\rtd(G, o)
&=
c^{-1} \sum_i \sum_x f(o_i, x) |\Stab(o_i)|^{-1}
\cr&=
c^{-1} \sum_i |\Stab(o_i)|^{-1}
\sum_j |\Stab(o_j)|^{-1} \sum_{x \in [o_j]} f(o_i, x) |\Stab(o_j)|
\cr&=
c^{-1} \sum_i |\Stab(o_i)|^{-1}
\sum_j |\Stab(o_j)|^{-1} \sum_{y \in [o_i]} f(y, o_j) |\Stab(y)|
\cr&\hskip1in\hbox{[by the Mass-Transport Principle for $\gh$]}
\cr&=
c^{-1} \sum_i |\Stab(o_i)|^{-1}
\sum_j |\Stab(o_j)|^{-1} \sum_{y \in [o_i]} f(y, o_j) |\Stab(o_i)|
\cr&\hskip1in\hbox{[by unimodularity of $\gh$]}
\cr&=
c^{-1} \sum_j \sum_y f(y, o_j) |\Stab(o_j)|^{-1}
\cr&=
\int \sum_y f(y, o) \,d\rtd(G, o)
\,.
}$$
Since $\rtd$ satisfies the Mass-Transport Principle, it is unimodular.

Conversely, suppose that $\rtd$ is a unimodular probability measure on rooted
versions of $\gh$. To see that $\gh$ is unimodular, consider any two
vertices $u, v$.
Define
$$
\rtd\big([x]\big) :=
\rtd\big([\gh, x]\big)
\,.
$$
We first show that $\rtd\big([u]\big) > 0$.
Every graph isomorphic to $\gh$ has a well-defined notion of vertices of
type $[u]$.
Let each vertex $x$ send mass 1 to each vertex of type $[u]$ that is
nearest to $x$.
This is a Borel function on $\gtwo$ if we transport no mass on graphs
that are not isomorphic to $\gh$.
The expected mass sent is positive, whence so is the expected mass
received.
Since only vertices of type $[u]$ receive mass, it follows that
$\rtd\big([u]\big) > 0$, as desired.

Let $f(x, y) := \I{\gp_{u, x} v}(y)$, where $\gp_{u, x} := \{ \gpe
\in\Aut(\gh) \st \gpe u = x \}$.
Note that $y \in \gp_{u, x} v$ iff $x \in \gp_{v, y} u$.
It is straightforward to check that $f$ is diagonally invariant under
$\Aut(\gh)$.
Note that
$$
|\Stab(x) y| \I{[x]}(o)
=
|\gp_{x, \bp} y|
$$
for all $x, y, \bp \in \verts(\gh)$.
Therefore, we have
$$\eqaln{
|\Stab(u) v| \rtd([u])
&=
\int |\gp_{u, o} v| \,d\rtd(\gh, o)
=
\int \sum_x \I{\gp_{u, o} v}(x) \,d\rtd(\gh, o)
\cr&=
\int \sum_x f(o, x) \,d\rtd(\gh, o)
=
\int \sum_x f(x, o) \,d\rtd(\gh, o)
\cr&\hskip1in\hbox{[by the Mass-Transport Principle for $\rtd$]}
\cr&=
\int \sum_x \I{\gp_{u, x} v}(o) \,d\rtd(\gh, o)
=
\int \sum_x \I{\gp_{v, o} u}(x) \,d\rtd(\gh, o)
\cr&=
\int |\gp_{v, o} u| \,d\rtd(\gh, o)
=
|\Stab(v) u| \rtd([v])
\,.
}$$
That is,
$$
|\Stab(u) v| \rtd([u])
=
|\Stab(v) u| \rtd([v])
\,.
\label e.rtdstab
$$
If $u$ and $v$ are in the same orbit, then $[u] = [v]$, so $\rtd([u]) =
\rtd([v])$. Since $\rtd([u]) > 0$, we obtain \ref e.unimod-schl-trof/.
That is, $\gh$ is unimodular. In general,
comparison of \ref e.rtdstab/ with \ref e.schl-trof/ shows \ref e.mustab/.
\Qed

Automorphism invariance for random unrooted networks on fixed underlying
graphs is also closely tied to unimodularity of random rooted networks.
Here, we shall need to distinguish between graphs, networks, and
isomorphism classes of rooted networks.
Recall that
$\overline \gh$ denotes a network whose
underlying graph is $\gh$ and $[\overline \gh, \bp]$ denotes an
equivalence class of networks $\overline \gh$ on $\gh$ with root $\bp$.

Let $\gh$ be a fixed connected unimodular graph satisfying \ref e.alm-trn/.
Fix a complete orbit section $\{\bp_i\}$ of $\vertex(\gh)$.
For a graph $\gh'$ and $x \in \vertex(\gh')$, $z \in
\vertex(\gh)$, let $\Phi(x, z)$ be the set of rooted isomorphisms, if any,
from $(\gh', x)$ to $(\gh, z)$.
When non-empty, this set carries a natural probability measure,
$\lift'_{(\gh', x; z)}$ arising from the Haar probability measure on
$\Stab(z)$.
When $\Phi(x, z) = \emptyset$, let $\lift'_{(\gh', x, z)} := 0$.
Define
$$
\lift_{(\gh', x)} := \sum_{z \in \vertex(\gh)} \lift'_{(\gh', x; z)}
\,.
$$
This is the analogue for isomorphisms from $\gh'$ to $\gh$
of Haar measure on $\Aut(\gh)$.
In particular, any $\gpe \in \Aut(\gh)$ pushes forward $\lift'_{(\gh', x,
z)}$ to $\lift'_{(\gh', x, \gpe z)}$.

For a graph $\gh'$ isomorphic to $\gh$ and $x \in \vertex(\gh')$, let
$\type(\gh', x) := \bp_i$ for the unique $\bp_i$ for which
$\Phi(x, \bp_i) \ne \emptyset$.
Note that $\lift_{(\gh', x)} = \lift_{(\gh', y)}$ when $\type(\gh', x) =
\type(\gh', y)$.

Every probability measure $\rtd$ on $\GG_*$ that is
concentrated on network classes whose underlying graph is $\gh$ induces a
probability measure $\lift_\rtd$ on unrooted networks on $\gh$:
$$
\lift_\rtd(\A)
:=
\int \int_{\textstyle \Phi\big(\bp, \type(\gh', \bp)\big)}
\I{\A}(\phi \overline \gh')  \,d\lift_{(\gh', \bp)}(\phi)
\,d\rtd([\overline \gh', \bp])
$$
for Borel sets $\A$ of networks on $\gh$.
It is easy to see that this is well defined (the choice of $(\overline
\gh', \bp)$ in its equivalence class not mattering).

\procl t.uni-vs-invar \procname{Invariance and Unimodularity}
Let $\gh$ be a fixed connected unimodular graph satisfying \ref e.alm-trn/.
Let $\nu$ be an $\Aut(\gh)$-invariant probability measure on unrooted
networks whose underlying graph is $\gh$.
Then randomly rooting the network as in \ref e.mustab/ gives a measure $\rtd
\in \uni$.
Conversely, let $\rtd \in \uni$ be supported on networks whose underlying
graph is $\gh$.
Then $\lift_\rtd$ is $\Aut(\gh)$-invariant.
\endprocl

\proof
The first part of the theorem is proved just as is the first part of \ref
t.fixedgraph/,
so we turn to the second part.
Let $\gpe_0 \in \Aut(\gh)$
and $F$ be a bounded Borel-measurable function of networks on $\gh$.
Invariance of $\lift_\rtd$ means that $\int F(\overline
\gh)\,d\lift_\rtd(\overline \gh) = \int F(\gpe_0 \overline
\gh)\,d\lift_\rtd(\overline \gh)$.
To prove that this holds, let
$$
f(\overline \gh', x, y) :=
\int_{\textstyle \Phi\big(x, \type(\gh', y)\big)
\atop \textstyle \cap \Phi\big(y, \gpe_0 \type(\gh, y)\big)}
F(\phi \overline \gh')\,d\lift_{(\gh', y)}(\phi)
\,.
$$
It is straightforward to check that
$f$ is well defined and Borel on $\gtwo$.
Therefore, unimodularity of $\rtd$ gives
$$\eqaln{
\int F(\gpe_0 \overline \gh)\,d\lift_\rtd(\overline\gh)
&=
\int \int_{\textstyle \Phi\big(\bp, \type(\gh', \bp)\big)}
F(\gpe_0 \phi \overline \gh') \,d\lift_{(\gh', \bp)}(\phi)
\,d\rtd([\overline \gh', \bp])
\cr&=
\int \int_{\textstyle \Phi\big(\bp, \gpe_0 \type(\gh', \bp)\big)}
F(\phi \overline \gh') \,d\lift_{(\gh', \bp)}(\phi)
\,d\rtd([\overline \gh', \bp])
\cr&=
\int \sum_{x \in \vertex(\gh')}
\int_{\textstyle \Phi\big(x, \type(\gh', \bp)\big)
\atop \textstyle \cap \Phi\big(\bp, \gpe_0 \type(\gh', \bp)\big)}
F(\phi \overline \gh') \,d\lift_{(\gh', \bp)}(\phi)
\,d\rtd([\overline \gh', \bp])
\cr&=
\int \sum_{x \in \vertex(\gh')}
f(\overline \gh', x, \bp)
\,d\rtd([\overline \gh', \bp])
\cr&=
\int \sum_{x \in \vertex(\gh')}
f(\overline \gh', \bp, x)
\,d\rtd([\overline \gh', \bp])
\cr&=
\int \sum_{x \in \vertex(\gh')}
\int_{\textstyle \Phi\big(\bp, \type(\gh', x)\big)
\atop \textstyle \cap \Phi\big(x, \gpe_0 \type(\gh', x)\big)}
F(\phi \overline \gh') \,d\lift_{(\gh', x)}(\phi)
\,d\rtd([\overline \gh', \bp])
\cr&=
\int \sum_{x \st \type(\gh', x) = \type(\gh', \bp)}
\int_{\textstyle \Phi\big(\bp, \type(\gh', x)\big)
\atop \textstyle \cap \Phi\big(x, \gpe_0 \type(\gh', x)\big)}
F(\phi \overline \gh') \,d\lift_{(\gh', x)}(\phi)
\,d\rtd([\overline \gh', \bp])
\cr&=
\int \sum_{x \st \type(\gh', x) = \type(\gh', \bp)}
\int_{\textstyle \Phi\big(\bp, \type(\gh', \bp)\big)
\atop \textstyle \cap \Phi\big(x, \gpe_0 \type(\gh', \bp)\big)}
F(\phi \overline \gh') \,d\lift_{(\gh', \bp)}(\phi)
\,d\rtd([\overline \gh', \bp])
\cr&=
\int \sum_{x \in \vertex(\gh')}
\int_{\textstyle \Phi\big(\bp, \type(\gh', \bp)\big)
\atop \textstyle \cap \Phi\big(x, \gpe_0 \type(\gh', \bp)\big)}
F(\phi \overline \gh') \,d\lift_{(\gh', \bp)}(\phi)
\,d\rtd([\overline \gh', \bp])
\cr&=
\int
\int_{\textstyle \Phi\big(\bp, \type(\gh', \bp)\big)}
F(\phi \overline \gh') \,d\lift_{(\gh', \bp)}(\phi)
\,d\rtd([\overline \gh', \bp])
\cr&=
\int F(\overline \gh)\,d\lift_\rtd(\gh)
\,.
\Qed
}$$

\procl r.anygroup
As this section shows, unimodular quasi-transitive graphs are special cases
of unimodular rooted networks.
However, sometimes one is interested in random networks on a graph $\gh$
that are not necessarily invariant under the full group $\Aut(\gh)$, but
only under some subgroup, $\gp \subset \Aut(\gh)$.
This is common when $\gh$ is a Cayley graph of $\gp$.
In this case, we could mark the edges by the generators they represent;
that is, if $x, y \in \gp$ and $y = x a$ with $a$ one of the generators
used to form $\gh$, then we can mark the edge $[x, y]$ at $x$ by $a$.
This makes the full automorphism group of the network
$\overline \gh$ equal to $\gp$, rather than to $\Aut(\gh)$.
The theory here then goes through with only a complication of notation.
However, given any graph $\gh$ and any closed subgroup $\gp \subset
\Aut(\gh)$ that acts quasi-transitively on $\gh$, we do not know whether it
is possible to mark the edges and vertices of $\gh$ to get a network whose
automorphism group is equal to $\gp$.
Yet, the theory for quasi-transitive subgroups is the same; see \BLPSgip.
\endprocl

\bsection{Random Walks and Extremality}{s.extreme}

Random walks on networks,
besides being of intrinsic interest,
form an important tool for studying networks.
A random walk is most useful when it has a stationary measure, in other
words, when the distribution of $(\gh, w_0)$ is the same as the
distribution of $(\gh, w_1)$, where $w_0$ is the initial location of the
random walk and $w_1$ is the next location of the random walk.

Consider simple random walk on a random graph chosen by a unimodular
probability measure $\rtd$ on rooted graphs, where we start the random walk
at the root.
Just as for finite graphs, we do not expect $\rtd$ to be stationary for the
random walk; rather, we get a stationary measure by biasing $\rtd$ by the
degree of the root.
The fact that this measure is stationary follows from the definition of
involution invariance;
in fact, the definition is precisely that 
the distribution of the isomorphism class of $(\gh, w_0,
w_1)$ is the same as that of $(\gh, w_1, w_0)$ when $(\gh, w_0)$ has the
distribution $\rtd$ biased by the degree of the root and $w_1$ is a uniform
random neighbor of the root.
This implies that simple random walk is
reversible, i.e., that the distribution of $\big((\gh, w_0), (\gh,
w_1)\big)$ is the same as the distribution of $\big((\gh, w_1), (\gh,
w_0)\big)$, where again 
$(\gh, w_0)$ has distribution $\rtd$ biased by the degree
of the root and $w_1$ is a uniform random neighbor of the 
root.\ftnote{*}{Note that the degree times counting measure is
reversible on every graph, regardless of
unimodularity of the measure on rooted graphs.}
If $\expdeg(\rtd) < \infty$, then we can normalize the biased measure
to obtain a probability measure.

In particular, recall from
\ref x.AGW/ the definition of the augmented Galton-Watson measure \AGW.
In \ref b.LPP:GW/, it was remarked in reference to the stationarity of $\AGW$
for simple random walk that ``unlike the situation for finite
graphs, there is no biasing in favor of vertices of large degree".
However, we now see that contrary to this remark,
the situations of finite graphs and $\AGW$ are, in fact, parallel.
That is because the biasing by the degree has already been made part of the
probability measure $\AGW$. The correct comparison of the uniform measure
on vertices of finite graphs is to the
unimodular Galton-Watson probability measure on trees, $\UGW$, because it
is for this measure that ``all vertices are equally likely to be the root".

More generally, we can consider
stationarity of random walk in a random environment with random scenery.
Here, if the graph underlies a network, the marks are not restricted to play a
passive role, but may, in fact, determine the transition probabilities (as
in \ref s.trace/) and provide a scenery for the random walk.
That is, a Borel function $p : \gtwo \to [0, 1]$, written as
$p : (\gh, x, y) \mapsto p_\gh(x, y)$, such that $\sum_{y \in \verts}
p_\gh(x, y) = 1$ for all vertices $x$
is called an {\bf environment}.
A Borel map $\nu : \GG_* \to (0, \infty)$, written
$\nu : (\gh, x) \mapsto \nu_{\gh}(x)$,
is called an {\bf initial bias}.
It is called {\bf $p$-stationary} if for all $\gh$, the measure
$\nu_\gh$ is stationary for
the random walk on $\gh$ given by the environment $p_\gh$.
Write $\paths$ for the set of (equivalence classes of)
pairs $\big((\gh, w_0), \Seq{w_n \st n \ge
0}\big)$ with $(\gh, w_0) \in \GG_*$ and $w_n \in \vertex(\gh)$.
Let $\Ph$ denote the distribution on $\paths$ of
the trajectory of the Markov chain
determined by the environment starting at $\bp$ with initial distribution
equal to $\rtd$ biased by $\nu_{\gh}(\bp)$.
That is, if $\traj_{(\gh, \bp)}$ denotes the probability measure on $\paths$
determined by the environment on $\gh$ with initial vertex
$w_0 = \bp$, then for all events $\ev B$, we have
$$
\Ph(\ev B)
:=
\int_{\GG_*} \traj_{(\gh, \bp)}(\ev B) \nu_{\gh}(\bp)
\,d\rtd(\gh, \bp)
\,.
$$
Let $\invar$ denote the $\sigma$-field of events (in the Borel
$\sigma$-field of $\GG_*$) that are invariant under non-rooted
isomorphisms.
To avoid possible later confusion,
note that this does not depend on the measure $\rtd$, so that even if there
are no non-trivial non-rooted isomorphisms $\rtd$-a.s., the $\sigma$-field
$\invar$ is still not equal (mod 0) to the $\sigma$-field of
$\rtd$-measurable sets.
It is easy to see that
for any $\rtd \in \uni$ and $\ev A \in \invar$ with $\rtd(\ev A) > 0$, the
probability measure $\rtd(\;\cbuldot \mid \ev A)$ is also unimodular.
Define the {\bf shift} $\SS : \paths \to \paths$ by
$$
\SS \big((\gh, w_0), \Seq{w_n}\big) := \big((\gh, w_1), \Seq{w_{n+1}}\big)
\,.
$$
The following extends Theorem 3.1 of \ref b.LS:rwrers/;
the proof is essentially the same.


\procl t.rwrers \procname{Random Walk in a Random Environment and Random
Scenery}
Let $\rtd$ be a unimodular probability measure on $\GG_*$.
Let $p_{\cbuldot}(\cbuldot)$ be an environment and
$\nu_{\cbuldot}(\cbuldot)$ be an initial bias that is $p$-stationary.
Let $\Ph$ be the corresponding measure on trajectories.
Then $\Ph$ is stationary for the shift.
If $p$ is also reversible with respect to $\nu_{\cbuldot}(\cbuldot)$, then
$\Ph$ is reversible, in other words, for all events $\ev A$, $\ev B$, we
have 
$$
\Ph[(\gh, w_0) \in \ev A, (\gh, w_1) \in \ev B]
=
\Ph[(\gh, w_1) \in \ev A, (\gh, w_0) \in \ev B]
\,.
$$
If
$$
\int \nu_{\gh}(\bp) \,d\rtd(\gh, \bp) = 1
\,,
\label e.proba
$$
then $\Ph$ is a probability measure.
\endprocl

\proof
The reversibility was not mentioned in prior work, so we give that proof
here.
Assuming that $p$ is $\nu$-reversible, we have 
$$\eqaln{
\Ph[(\gh, w_0) \in \ev A, (\gh, w_1) \in \ev B]
&=
\Ebig{\sum_{x \in \vertex(\gh)} \I{\ev A}(\gh, \bp) \nu_\gh(\bp) p_\gh(\bp,
x) \I{\ev B}(\gh, x)}
\cr&=
\Ebig{\sum_{x \in \vertex(\gh)} \I{\ev A}(\gh, \bp) \nu_\gh(x) p_\gh(x,
\bp) \I{\ev B}(\gh, x)}
\,.
\cr
}$$
The Mass-Transport Principle now gives that this 
$$
=
\Ebig{\sum_{x \in \vertex(\gh)} \I{\ev A}(\gh, x) \nu_\gh(\bp) p_\gh(\bp,
x) \I{\ev B}(\gh, \bp)}
=
\Ph[(\gh, w_1) \in \ev A, (\gh, w_0) \in \ev B]
\,. \Qed
$$

\procl r.choice
This theorem is made more useful by noticing that for any $\rtd \in \uni$,
there is a choice of $p_{\cbuldot}(\cbuldot)$ and
$\nu_\cbuldot(\cbuldot)$ that satisfies all the hypotheses, including
\ref e.proba/.
For example, if $F_\gh(x)$
denotes $\sum_{y \sim x} 1/\deg_\gh(y)$,
then let $p_\gh(x, y) := 1/[F_\gh(x) \deg_\gh(y)]$ and
$\nu_\gh(x) := Z^{-1} F_\gh(x)/\deg_\gh(x)$,
where 
$$
Z := \int F_\gh(\bp)/\deg_\gh(\bp) \,d\rtd(\gh, \bp)
\,.
$$
It is clear that $p$ is an environment.
Since $F_\gh(\bp) \le \sum_{y \sim x} 1 = \deg_\gh(\bp)$, we also have that
$Z < \infty$, so that $\nu$ is a $p$-stationary initial bias and $p$ is
$\nu$-reversible.
\endprocl

Given a network with positive edge
weights and a time $t > 0$, form the {\bf transition operator} $P_t$ for
continuous-time random walk whose rates are the edge weights; in the case
of unbounded weights (or degrees), we take the minimal process, which dies
after an explosion.
That is, if the entries of a matrix $A$ indexed by the vertices are equal
off the diagonal to the negative of the edge weights and the diagonal
entries are chosen to make the row sums zero, then $P_t := e^{-A t}$; in
the case of unbounded weights, we take the self-adjoint extension of $A$
corresponding to the minimal process.
The matrix $A$ is called 
the {\bf Laplacian} of the network; it is
the negative of the
{\bf infinitesimal generator} of the random walk.

\procl c.continuous-time-stationary
Suppose that $\rtd \in \uni$ is carried by networks with non-negative edge
weights such that the corresponding continuous-time Markov chain has no
explosions a.s.
Then $\rtd$ is stationary and reversible.
\endprocl

\proof
Fix $t > 0$ and let $p_\gh(x, y) := P_t(x, y)$.
It is well known that $p$ is reversible with respect to
the uniform measure $\nu_\gh \equiv 1$.
Thus, \ref t.rwrers/ applies.
\Qed

We can also obtain a sufficient condition for lack of explosions:

\procl c.no-explode
Suppose that $\rtd \in \uni$ is carried by networks with non-negative edge
weights $c_\gh(e)$ such that $Z := \E[\sum_{x \sim \bp} c_\gh(\bp, x)] <
\infty$.
Then the corresponding continuous-time Markov chain has no explosions.
\endprocl

\proof
In this case, consider the discrete-time Markov chain corresponding to
these weights.
It has a stationary probability measure arising from the choice
$\nu_\gh(x) := \sum_{y \sim \bp} c_\gh(x, y)/Z$.
It is well known that explosions occur iff 
$$
\sum_{n \ge 0} \nu_\gh(w_n)^{-1} < \infty
$$
with positive probability.
However, stationarity guarantees that this sum is infinite a.s.\
(by the Poincar\'e recurrence theorem).
\Qed

\procl r.poss-explode
It is possible for explosions to occur:
For example, consider the uniform spanning tree $T$ in $\Z^2$ (see
\BLPSusf). The only fact we use about $T$ is that it has one end a.s.\
and has an invariant distribution.
Let $c_\gh(e) := 0$ for $e \notin T$ and
$c_\gh(e) := 2^{f(e)}$ when $e \in T$
and $f(e)$ is the number of vertices in the finite component of $T
\setminus e$.
Then it is easy to verify that the corresponding continuous-time Markov
chain explodes a.s.

Furthermore, explosions may occur
on a fixed transitive graph that is not unimodular, 
even if the condition in \ref c.no-explode/ is satisfied.
To see this, let
$\xi$ be a fixed end of a regular tree $T$ of degree 3. 
Thus, for every vertex $x$ in $T$, there is a unique ray $x_\xi := 
\Seq{x_0=x, x_1, x_2, \ldots}$ starting at $x$ such that $x_\xi$ and
$y_\xi$ differ by only finitely many vertices for any pair $x$, $y$.
Call $x_1$ the {\bf $\xi$-parent} of $x$, call $x$ a {\bf $\xi$-child} of
$x_1$, and call
$x_2$ the {\bf $\xi$-grandparent} of $x$. Let $\gh$ be the graph obtained
from $T$ by adding the edges $(x, x_2)$ between each $x$ and its
$\xi$-grandparent.
Then $\gh$ is a transitive graph, first mentioned by \ref b.Trofimov/.
In fact, every automorphism of $\gh$ fixes $\xi$.
Now consider the following random weights on $\gh$.
Put weight 0 on every edge in $\gh$ that is not in $T$.
For each vertex of $\gh$, declare open the edge to precisely
one of its two $\xi$-children, chosen uniformly and
independently for different vertices.
The open components are rays.
Let the weight of every edge that is not open also be 0.
If an edge $(x, y)$ between a vertex $x$ and its $\xi$-parent $y$ is open and
$y$ is at distance $n$ from the beginning of the open ray containing $(x, y)$,
then let the weight of the edge be $(3/2)^n$.
Since this event has probability $1/2^{n+1}$, the condition of \ref
c.no-explode/ is clearly satisfied.
It is also clear that the Markov chain explodes a.s.
\endprocl

The class $\uni$ of unimodular probability measures on $\GG_*$ is clearly
convex.
An element of $\uni$ is called {\bf extremal} if it cannot be written as a
convex combination of other elements of $\uni$.
We shall show that the extremal measures are those for which $\invar$
contains only sets of measure 0 or 1.
Intuitively,
they are the extremal measures for unrooted networks since
the distribution of the root is forced given the distribution of the
unrooted network.
For example, one may show that $\UGW$ is extremal when conditioned on
non-extinction.
First, we show the following ergodicity property, analogous to Theorem 5.1
of \ref b.LS:indist/.
Recall that a $\sigma$-field is called {\bf $\rtd$-trivial} if all its
elements have measure 0 or 1 with respect to $\rtd$.

\procl t.erg \procname{Ergodicity}
Let $\rtd$ be a unimodular probability measure on $\GG_*$.
Let $p_{\cbuldot}(\cbuldot)$ be an environment that satisfies
$$
\all{\gh} \all{x, y \in \vertex(\gh)} \quad x \sim y \Longrightarrow
p_\gh(x, y) > 0
\label e.getsall
$$
and
$\nu_{\cbuldot}(\cbuldot)$ be an initial bias that is $p$-stationary and
satisfies \ref e.proba/.
Let $\Ph$ be the corresponding probability measure on trajectories.
If $\invar$ is $\rtd$-trivial, then
every event that is shift invariant is $\Ph$-trivial.
More generally, the events $\ev B$ in the $\Ph$-completion of the
shift-invariant $\sigma$-field are those of the form 
$$
\ev B
=
\big\{ \big((\gh, \bp), w\big) \in \paths \st (\gh, \bp) \in \ev A \big\}
\xor \ev C
\label e.shiftev
$$
for some $\ev A \in \invar$ and some event $\ev C$ with $\Ph(\ev C) = 0$.
\endprocl

\proof
Let $\ev B$ be a shift-invariant event.
As in the proof of Theorem 5.1 of \ref b.LS:indist/, we have $\traj_{(\gh,
\bp)}(\ev B) \in \{ 0, 1 \}$ $\rtd$-a.s.
The set $\ev A$ of $(\gh, \bp)$ where this probability equals 1 is in
$\invar$ by \ref e.getsall/, and a little thought reveals that
\ref e.shiftev/ holds for some $\ev C$ with $\Ph(\ev C) = 0$.
If $\invar$ is $\rtd$-trivial, then
$\rtd(\ev A) \in \{0, 1\}$, whence $\Ph(\ev B) \in \{0, 1\}$ as desired.
Conversely, every event $\ev B$ of the form \ref e.shiftev/ is clearly 
in the $\Ph$-completion of the shift-invariant $\sigma$-field.
\Qed

We may regard the space $\paths$ as the space of sequences of rooted
networks, where all roots belong to the same network.
Thus, $\paths$ is the natural trajectory space for the Markov chain with
the transition probability from $(\gh, x)$ to $(\gh, y)$ given by $p_\gh(x,
y)$.
With this interpretation, \ref t.erg/ says that this Markov chain is
ergodic when $\invar$ is $\rtd$-trivial.
The next theorem says that this latter condition is, in turn, equivalent to
extremality of $\rtd$.

\procl t.findroot \procname{Extremality}
A unimodular probability measure $\rtd$ on $\GG_*$ is
extremal iff $\invar$ is $\rtd$-trivial.
\endprocl

\comment{
For the entire proof, choose an environment and stationary initial bias
that satisfy \ref e.proba/ and \ref e.getsall/.
Also, let $\kappa$ be the law of $(\gh, \bp, \bp')$ on $\gtwo$ that couples
$\rtd$ and $\rtd'$.

Let $\A$ be an event of $\GG_*$.
Let $\alpha$ be the function on $\paths$ that gives
the asymptotic frequency of visits to $\A$:
$$
\alpha\big(((\gh, w_0), \Seq{w_n})\big) :=
\liminf_{N \to\infty} {1 \over N} |\{n \le N \st (\gh, w_n) \in \A\}|
\,.
$$
Then $\alpha$ is a shift-invariant function.
For any $t \in [0, 1]$, let $\ev B_t := \{(\gh, \bp) \st \traj_{(\gh,
\bp)}[\alpha \le t] \in \{0, 1\}\}$.
Then $\rtd(\ev B_t), \rtd'(\ev B_t) \in \{ 0, 1 \}$,
as we saw in the proof of \ref t.erg/.
Therefore, there exists a Borel function $f : \GG_* \to [0, 1]$ such
that $\rtd\big[\traj_{(\gh, \bp)}[\alpha = f(\gh, \bp)] = 1 \big] = 1$.
Since a random walk has positive probability of going from $\bp$ to $\bp'$,
it follows that
$$
\kappa[\big[\traj_{(\gh, \bp)}[\alpha = f(\gh, \bp) = f(\gh, \bp')] = 1 \big]
= 1
\,.
\label e.same
$$
Let $\nu \rtd$ stand for the measure $d(\nu \rtd)(\gh, \bp) = \nu_\gh(\bp)
d\rtd(\gh, \bp)$, and likewise for $\nu\rtd'$.
By \ref t.rwrers/, these are stationary measures for simple random walk,
whence the ergodic theorem gives
$$\eqaln{
\int \I\A \nu_\gh(\bp) \,d\rtd(\gh, \bp)
&=
\int f(\gh, \bp) \nu_\gh(\bp) \,d\rtd(\gh, \bp)
}$$
}%

\proof
Let $\A \in \invar$. If $\A$ is
not $\rtd$-trivial, then we may write $\rtd$ as a convex combination of $\rtd$
conditioned on $\A$ and $\rtd$ conditioned on the complement of $\A$. 
Each of these two new probability measures is unimodular,
yet distinct, so $\rtd$ is not extremal.

Conversely, suppose that $\invar$ is $\rtd$-trivial.
Choose an environment and stationary initial bias
that satisfy \ref e.proba/ and \ref e.getsall/, as in \ref r.choice/.
Let $\A$ be an event of $\GG_*$.
Let $\alpha$ be the function on $\paths$ that gives
the frequency of visits to $\A$:
$$
\alpha\big((\gh, w_0), \Seq{w_n})\big) :=
\liminf_{N \to\infty} {1 \over N} |\{n \le N \st (\gh, w_n) \in \A\}|
\,.
$$
\ref t.rwrers/ allows us to apply the ergodic theorem to deduce that
$\int \alpha \,d\Ph = (\nu\rtd)(\A)$,
where $\nu \rtd$ stands for the measure $d(\nu \rtd)(\gh, \bp) = \nu_\gh(\bp)
d\rtd(\gh, \bp)$.
On the other hand, $\alpha$ is a shift-invariant function, which,
according to \ref t.erg/, means that $\alpha$ is a constant $\Ph$-a.s.
Thus, we conclude that
$\alpha = (\nu\rtd)(\A)$ $\Ph$-a.s.
Consider any non-trivial
convex combination of two unimodular probability
measures, $\rtd_1$ and $\rtd_2$, that gives $\rtd$.
Then $\Ph$ is a (possibly different) convex combination of $\Ph_1$ and $\Ph_2$.
The above applies to each of $\Ph_i$ ($i=1, 2$) and the associated
probability measures $a_i \nu \rtd_i$, where
$a_i := \left(\int \nu_\gh(\bp) \,d\rtd_i(\gh, \bp)\right)^{-1}$.
Therefore, we obtain that $a_1(\nu\rtd_1)(\A) = (\nu\rtd)(\A) =
a_2(\nu\rtd_2)(\A)$.
Since this holds for all $\A$, we obtain $a_1 (\nu\rtd_1) = a_2
(\nu\rtd_2)$.
Since $\rtd_1$ and $\rtd_2$ are probability measures, this is the same as
$\rtd_1 = \rtd_2$, whence $\rtd$ is extremal.
\Qed

We define the {\bf speed} of a path $\Seq{w_n}$ in a graph $\gh$ to be
$\lim_{n \to\infty} \dist_\gh(w_0, w_n)/n$ when this limit exists, where
$\dist_\gh$ indicates the distance in the graph $\gh$.

The following extends Lemma 4.2 of \ref b.BLS:pert/.

\procl p.speedexists \procname{Speed Exists}
Let $\rtd$ be a unimodular probability measure on $\GG_*$ with an
environment and stationary initial distribution $\nu_{\cbuldot}(\cbuldot)$ with
$\int \nu_{\gh}(\bp) \,d\rtd(\gh, \bp) = 1$, so that the
associated random walk distribution $\Ph$ is a probability measure.
Then the speed of random walk exists $\Ph$-a.s.\ and is equal $\Ph$-a.s.\
to an $\invar$-measurable function.
The same holds for simple random walk when $\expdeg( \rtd) < \infty$.
\endprocl

\proof
Let $f_n\big((\gh, \bp), w\big) := \dist_\gh\big(w(0), w(n)\big)$. Clearly
$$
f_{n+m}\big((\gh, \bp), w\big)
\le f_n\big((\gh, \bp), w\big) + f_m\big(\SS^n \big((\gh, \bp), w\big)\big)
\,,
$$
so that the Subadditive Ergodic Theorem ensures
that the speed $\lim_{n \to\infty} f_n\big( (\gh, \bp), w\big)/n$ exists
$\Ph$-a.s.
Since the speed is shift invariant, \ref t.erg/ shows that the speed is
equal $\Ph$-a.s.\ to an $\invar$-measurable function.
The same holds for simple random walk since it has an equivalent stationary
probability measure (degree biasing) when $\expdeg(\rtd) < \infty$.
\Qed

In the case of simple random walk on trees, we can actually calculate the
speed\ftnote{*}{The publisher inadvertently changed the following proposition
to a theorem in the published version. Also, the published version had an
incorrect formula for \ref e.genspeed/.}:

\procl p.treespeed \procname{Speed on Trees}
Let $\rtd \in \uni$ be concentrated on infinite trees.
If $\rtd$ is extremal and $\expdeg(\rtd) < \infty$, then the speed of
simple random walk is $\Ph$-a.s.\ 
$$
1 - {2 \over \overline{\deg}(\mu)}
\,.
\label e.genspeed
$$
\endprocl

\proof Given a rooted tree $(\gh, \bp)$ and $x \in \vertex(\gh)$, write
$|x|$ for the distance in $\gh$ between $\bp$ and $x$.
The speed of a path $\Seq{w_n}$ is the limit
$$
\lim_{n \to \infty} {1 \over n} |w_n|
= \lim_{n \to \infty} {1 \over n} \sum_{k = 0}^{n-1} \bigl(|w_{k+1}|
       - |w_{k}|\bigr)\,.
$$
Now the strong law of large numbers for martingale differences (\ref
b.Feller:book/, p.~243) gives
$$
\lim_{n \to \infty} {1 \over n} \sum_{k = 0}^{n-1} \bigl(|w_{k+1}|
       - |w_{k}|\bigr)
= \lim_{n \to \infty} {1 \over n} \sum_{k = 0}^{n-1} \E\bigl[|w_{k+1}|
       - |w_{k}| \bigm| \Seq{w_i \st i \le k}\bigr] \quad\hbox{a.s.} 
$$
Provided $w_k \ne o$, the $k$th term on the right equals
$$
{\deg_\gh w_k - 2 \over \deg_\gh w_k }\,. 
$$
Since $\gh$ is a.s.\ infinite,
$w_k = \bp$ for only a set of $k$ of density 0 a.s., whence the speed equals
$$
\lim_{n \to \infty} {1 \over n} \sum_{k = 0}^{n-1} {\deg_\gh w_k - 2 \over
\deg_\gh w_k }\,. 
$$
Since this is the limit of averages of an ergodic stationary sequence for
the measure $d\sigma(\gh, \bp) = \deg_\gh (\bp)\, d\rtd(\gh,
\bp)/\overline{\deg}(\rtd)$,
the ergodic theorem
tells us that it converges a.s.\ to the $\sigma$-mean of an element of the
sequence,
$$
\int {\deg_\gh(\bp) - 2 \over \deg_\gh(\bp)} d\sigma(\gh, \bp)
\,,
$$
which is the same as \ref e.genspeed/.
\Qed

When we study percolation, the following consequence will be useful.

\procl p.transtrees \procname{Comparison of Transience on Trees}
Suppose $\rtd \in \uni$ is concentrated on networks whose underlying graphs
are trees that are transient for simple random walk.
Suppose that the mark space is $(0, \infty)$, that marks $\mkmp(\cbuldot,
\cbuldot)$ on edges are the same at both endpoints, that the environment is
$p_\gh(x, y) := \mkmp(x, y)/\nu_\gh(x)$, where $\nu_\gh(x) := \sum_{y \sim
x} \mkmp(x, y)$, and that $\int \nu_{\gh}(\bp) \,d\rtd(\gh, \bp) = 1$, so
that the associated random walk distribution $\Ph$ is a probability
measure.
Then random walk is also transient with respect to the environment
$p_\cbuldot(\cbuldot)$ $\rtd$-a.s.
\endprocl

\proof
Let $\ev A$ be the set of $p_\cbuldot(\cbuldot)$-recurrent networks.
Suppose that $\rtd(\ev A) > 0$.
By conditioning on $\ev A$, we may assume without loss of
generality that $\rtd(\ev A) = 1$.
By \ref t.deg2/ and the recurrence of simple random walk on trees with at
most two ends, we have $\expdeg(\rtd) > 2$.
Let $\epsilon > 0$ be sufficiently small that 
$$
\int |\{x \sim \bp \st \mkmp(\bp, x) \ge \epsilon\}| \,d\rtd(\gh, \bp) > 2
\,.
$$
Since finite trees have average degree strictly less than 2, is follows
that the subnetwork $(\gh_\epsilon, \bp)$, defined to be the connected
component of $\bp$ formed by the edges with marks at least $\epsilon$,
is infinite with positive probability.
Let $\rtd'$ be the law of $(\gh_\epsilon, \bp)$ when $(\gh, \bp)$ has the
law $\rtd$, and conditioned on the event $\ev B$ that 
$(\gh_\epsilon, \bp)$ is infinite.
Then $\rtd' \in \uni$ and $\expdeg(\rtd') > 2$.
By \ref p.treespeed/, simple random walk has positive speed $\rtd'$-a.s.,
so, in particular, is transient a.s.
Now simple random walk is the walk corresponding to all edge weights in
$\gh_\epsilon$ being, say, $\epsilon$.
Rayleigh's monotonicity principle (\ref b.DoyleSnell/ or \ref b.LP:book/)
now implies that random walk is transient $\Ph$-a.s.\ on $\ev B$.
Thus, $\ev A \cap \ev B = \emptyset$.
This contradicts our initial assumption that $\rtd(\ev A) = 1$.
\Qed

The converse of \ref p.transtrees/ is not true, as there are transient
reversible random walks on 1-ended trees (see \ref x.halfplane/ for an
example of such graphs; weights can be defined appropriately).
Also, \ref p.transtrees/ does not extend to arbitrary networks, as one may
construct an invariant network on $\Z^3$ that gives a recurrent random
walk.

Given two probability measures $\rtd$ and $\rtd'$
on rooted networks and one of the standard
notions of product networks, one can define the {\bf independent product}
$\rtd \itimes \rtd'$ of
the two measures by choosing a network from each measure independently and
taking their product, rooted at the ordered pair of the original roots.

\procl p.product \procname{Product Networks}
Let $\rtd$ and $\rtd'$ be two unimodular probability measures on $\GG_*$.
Then their independent product $\rtd \itimes \rtd'$ is also unimodular.
If $\rtd$ and $\rtd'$ are both extremal, then so is $\rtd \itimes \rtd'$.
\endprocl

\proof
Let $\gh_n$ and $\gh'_n$ be finite connected networks whose random weak
limits are $\rtd$ and $\rtd'$, respectively.
Then $\gh_n \times \gh'_n$
clearly has random weak limit $\rtd \itimes \rtd'$, whence the product
is unimodular.
Now suppose that both $\rtd$ and $\rtd'$ are extremal.
Let $\ev A \in \invar$.
Then $\ev A_{(\gh', \bp')} :=
\{(\gh, \bp) \st (\gh, \bp) \times (\gh', \bp') \in \ev A \} \in
\invar$ since $\Aut(\gh) \times \Aut(\gh') \subseteq \Aut(\gh \times
\gh')$.
Therefore, $\rtd \ev A_{(\gh', \bp')} \in \{0, 1\}$.
On the other hand, $\ev A_{(\gh', \bp')} = \ev A_{(\gh', \bp'')}$ for all
$\bp'' \in \verts(\gh')$ because $\ev A \in \invar$.
Therefore,
$\ev B := \{(\gh', \bp') \st \rtd \ev A_{(\gh', \bp')}  = 1\} \in \invar$,
whence $\rtd' \ev B \in \{0, 1\}$.
Hence Fubini's theorem tells us that $(\rtd \itimes \rtd')(\ev A) \in
\{0, 1\}$, as desired.
\Qed

\procl r.weakmixing
Another type of product that can sometimes be defined does not always
produce an extremal measure from two extremal measures.
That is, suppose that
$\rtd$ and $\rtd'$ are two extremal unimodular probability measures on
$\GG_*$ that, for simplicity, we assume are concentrated on networks with a
fixed underlying transitive graph, $\gh$.
Let $\rtd''$ be the measure on networks given by taking a fixed root,
$\bp$, and choosing the marks as $(\mkmp, \mkmp')$, where $\mkmp$ gives a
network with law $\rtd$, $\mkmp'$ gives a network with law $\rtd'$, and
$\mkmp$, $\mkmp'$ are independent.
Then it may be that $\rtd''$ is not extremal.
For an example, consider the following.
Fix an irrational number, $\alpha$.
Given $x \in [0, 1]$, form the network $\gh_x$ on the
integer lattice graph by marking each integer $n$ with the indicator
that the fractional part of $n\alpha + x$ lies in $[0, 1/2]$.
Let $\rtd$ be the law of $(\gh_x, 0)$ when $x$ is chosen uniformly.
Then $\rtd$ is unimodular and extremal (by ergodicity of Lebesgue measure
with respect to rotation by $\alpha$), but if $\rtd' = \rtd$ and $\rtd''$
is the associated measure above, then $\rtd''$ is not
extremal since when the marks come from $x, y \in [0, 1]$,
the fractional part of $x - y$ is $\invar$-measurable.
\endprocl

It may be useful to keep in mind the vast difference between stationarity
and reversibility in this context.
For example, let $T$ be a 3-regular tree and $\zeta$ be an end of $T$.
Mark each edge by two independent random variables, one that is uniform on
$[0, 1]$ and the other uniform on $[1, 2]$, with the latter one
at its endpoint closer to $\zeta$ and with all these random variables
mutually independent for different edges.
Then simple random walk is stationary in this scenery, but not
reversible, even though $T$ is a Cayley graph.

\bsection{Trace and Stochastic Comparison}{s.trace}

There is a natural trace associated to every measure in $\uni$.
This trace is useful for making various comparisons.
We illustrate this by extending results of \ref b.PittSC/ and
\ref b.FontesMathieu/ on return
probabilities of continuous-time random walks.

Suppose that $\rtd$ is a unimodular probability measure on $\GG_*$.
Consider the Hilbert space $H := \int^\oplus \ell^2\big(\vertex(\gh)\big)
\,d\rtd(\gh, \bp)$, a direct integral (see, e.g., \ref b.Nielsen/ or \ref
b.KR/, Chapter 14).
Here, we always choose the canonical representative for each network,
which, recall, is a network on the vertex set $\N$.
The space $H$ is defined as the set
of ($\rtd$-equivalence classes of) $\rtd$-measurable functions $f$ defined on
canonical rooted networks $(\gh, \bp)$ that satisfy $f(\gh, \bp) \in
\ell^2(\N)$ and $\int \|f(\gh, \bp)\|^2 \,d\rtd(\gh, \bp) < \infty$.
We write $f = \int^\oplus f(\gh, \bp) \,d\rtd(\gh, \bp)$.
The inner product is given by $(f, g) := \int \big(f(\gh, \bp), g(\gh,
\bp)\big) \,d\rtd(\gh, \bp)$.
Let $T : (\gh, \bp) \mapsto T_{\gh, \bp}$ be a measurable assignment of
bounded linear operators on $\ell^2\big(\vertex(\gh)\big) = \ell^2(\N)$
with finite
supremum of the norms $\|T_{\gh, \bp}\|$.
Then $T$ induces a bounded linear operator $T := T^\rtd := \int^\oplus T_{\gh,
\bp} \,d\rtd(\gh, \bp)$ on $H$ via
$$
T^\rtd : \int^\oplus f(\gh, \bp)
\,d\rtd(\gh, \bp) \mapsto \int^\oplus T_{\gh, \bp} f(\gh, \bp)
\,d\rtd(\gh, \bp)
\,.
$$
The norm $\|T^\rtd\|$ of $T^\rtd$ is the $\rtd$-essential supremum of
$\|T_{\gh, \bp}\|$.
Identify each $x \in \vertex(\gh)$ with the vector $\II x \in
\ell^2\big(\vertex(\gh)\big)$.
Let $\alg = \alg(\rtd)$ be the von Neumann algebra of ($\mu$-equivalence
classes of) such maps $T$ that are equivariant in
the sense that for all network isomorphisms $\phi : \gh_1 \to \gh_2$, all
$\bp_1, x, y \in \verts(\gh_1)$ and all $\bp_2 \in \verts(\gh_2)$,
we have $(T_{\gh_1, \bp_1} x, y) = (T_{\gh_2, \bp_2} \phi x, \phi y)$.
For $T \in \alg$, we have in particular that $T_{\gh, \bp}$ depends on
$\gh$ but not on the root $\bp$, so we shall simplify our notation and
write $T_\gh$ in place of $T_{\gh, \bp}$.
For simplicity, we shall even write $T$ for $T_\gh$ when no confusion can
arise.
Recall that if $S$ and $T$ are self-adjoint operators on a Hilbert space
$H$, we write $S \le T$ if $\ip{S u, u} \le \ip{T u, u}$ for all $u \in H$.
We claim that
$$
\tr(T) := \tr_\rtd(T)
:= \Ebig{(T_\gh \bp, \bp)} := \int (T_\gh \bp, \bp)
\,d\rtd(\gh, \bp)
$$
is a {\bf trace} on $\alg$, i.e., $\tr(\cbuldot)$ is linear,
$\tr(T) \ge 0$ for $T \ge 0$, and $\tr(S T) = \tr(T
S)$ for $S, T \in \alg$.
Linearity of $\tr$ is obvious.
Also, the second property is obvious since the integrand is nonnegative for $T
\ge 0$. The third property follows from the Mass-Transport Principle:
We have
$$\eqaln{
\Ebig{(S T \bp, \bp)}
&=
\Ebig{(T \bp, S^* \bp)}
=
\Ebigg{\sum_{x \in \vertex(\gh)} (T \bp, x) (x, S^* \bp)}
\cr&=
\Ebigg{\sum_{x \in \vertex(\gh)} (T \bp, x) (S x, \bp)}
=
\Ebigg{\sum_{x \in \vertex(\gh)} (T x, \bp) (S \bp, x)}
\cr&=
\Ebig{(T S \bp, \bp)}
\,.
}$$
In order to justify this use of the Mass-Transport Principle, we check
absolute integrability:
$$\eqaln{
\Ebigg{\sum_{x \in \vertex(\gh)} |(T \bp, x) (x, S^* \bp)|}
&\le
\bigg(\Ebigg{\sum_{x \in \vertex(\gh)} |(T \bp, x)|^2}
\Ebigg{\sum_{x \in \vertex(\gh)} |(x, S^* \bp)|^2} \bigg)^{1/2}
\cr&=
\bigg(\Ebig{\|T \bp\|^2} \Ebig{\|S^* \bp\|^2} \bigg)^{1/2}
\le
\|T\| \cdot \|S\|
<\infty
\,.
}$$

A general property of traces that are finite and normal, as ours is,
is that if $S \le T$,
then $\tr f(S) \le \tr f(T)$ for any increasing function $f : \R \to \R$.
One proof is as follows.
First, if $f \ge 0$ and $T$ is self-adjoint, then $f(T) \ge 0$.
Second, if $S, T \ge 0$, then $\tr (S T) = \tr \big(S^{1/2} T S^{1/2}\big)
\ge 0$ since $S^{1/2} T S^{1/2} = (T^{1/2} S^{1/2})^* (T^{1/2}S^{1/2})  \ge 0$.
Third, if $f$ is an increasing polynomial and $S \le T$, then
$$
{d \over d z} \Big(\tr f\big(S + z(T-S)\big)\Big)
=
\tr \Big(f'\big(S + z(T-S)\big)(T-S)\Big)
\ge 0
$$
for $z \ge 0$
since $f' \ge 0$, $S + z(T-S)$ is self-adjoint, and $0 \le T-S$.
This shows that with these restrictive hypotheses, $\tr f(S) \le \tr f(T)$.
Fourth, any monotone increasing function can be approximated by 
an increasing polynomial.
This shows the result in general.
See \ref b.BrownKosaki/, pp.~6--7 for another proof. They stated the result
only for continuous $f$ with $f(0) = 0$ because they dealt with more
general traces and operators, but such restrictions are not needed in our
situation.
In fact, their proof shows that
$\tr f(S) \le \tr f(T)$
for bounded increasing $f \colon \R \to \R$ and $0 \le S \le T$ that are
$\tr$-measurable operators affiliated to $\alg$. 
Definitions are as follows.
A closed densely defined operator is {\bf affiliated}
with $\alg$ if it commutes with all unitary operators
that commute with $\alg$.
We call an affiliated operator $T$ {\bf
$\tr$-measurable} if for all $\epsilon > 0$,
there is an orthogonal projection $E \in \alg$ whose
image lies in the domain of $T$ and $\tr(E^\perp) < \epsilon$.

Recall from \ref s.extreme/ that given
a network with positive edge
weights and a time $t > 0$, we form the transition operator $P_t$ for
continuous-time random walk whose rates are the edge weights; in the case
of unbounded weights (or degrees), we take the minimal process, which dies
after an explosion.
If $A$ is the Laplacian of the network,
then $P_t := e^{-A t}$.
The Laplacian $A$ as an operator belongs to $\alg(\rtd)$ if the sum of the
edge weights at $\bp$ is $\rtd$-essentially
uniformly bounded. In any case, $A$ is affiliated to
$\alg(\rtd)$. Also, $A$ is $\tr_\rtd$-measurable because 
if $E_n$ denotes the orthogonal projection to the space of functions that
are nonzero only on those $(\gh, \bp)$
where the sum of the edge weights at $\bp$ and $x$ is at most $n$ for every
$x \sim \bp$, then $\lim_{n \to\infty} \tr(E_n^\perp) = 0$ and $\|A
E_n\| \le n$.

\procl t.rtnpr \procname{Return Probabilities}
Let the mark space be $\R^+$ and let $\rln$ be $\le$.
Let $\rtd_i \in \uni$ have edge weights that are the same at both ends of each
edge ($i = 1, 2$).
Suppose that there is a unimodular $\rln$-coupling $\nu$ of
$\rtd_1$ to $\rtd_2$.
Let $P^{(i)}_t$ be the transition operators corresponding to the edge
weights ($i = 1, 2$).
Then 
$$
\int P^{(1)}_t(\bp, \bp) \,d\rtd_1(\gh, \bp) \ge \int P^{(2)}_t(\bp,
\bp) \,d\rtd_2(\gh, \bp)
$$
for all $t > 0$.
\endprocl

\proof
The Laplacians $A^{(i)}$ affiliated to
$\alg(\nu)$ satisfy $A^{(1)} \le A^{(2)}$, so that for all $t > 0$, we have
$-A^{(1)} t \ge -A^{(2)} t$.
Therefore
$\int P^{(1)}_t(\bp, \bp) \,d\rtd_1(\gh, \bp) = \tr_{\nu}(e^{-A^{(1)} t})
\ge \tr_{\nu}(e^{-A^{(2)} t}) = \int P^{(2)}_t(\bp, \bp)
\,d\rtd_2(\gh, \bp)$.
\Qed

\ref t.rtnpr/ extends a result of \ref b.FontesMathieu/, who proved it in the
case of $\Z^d$ for processes without explosions.
\ref b.PittSC/, Lemma 3.1, prove an analogous
comparison result for Cayley graphs, but with different assumptions on the
pairs of rates (which are deterministic for them).
The case of \ref t.rtnpr/ specialized to
finite networks was proved earlier by Benjamini and Schramm;
see Theorem 3.1 of \ref b.HeicklenHoffman/.

This theorem also gives a partial answer to a
question of Fontes and Mathieu, who asked whether the same holds when
$\rtd_i$ are supported on a single Cayley graph, are invariant under the
group action, and are $\rln$-related.
For example, the theorem shows that
it holds when the networks (which are the environments for the
random walks) are given by i.i.d.\ edge marks, since in such a case, it is
trivial that being $\rln$-related implies the existence of a unimodular
$\rln$-coupling.
Also, in the amenable case, the existence of a unimodular $\rln$-coupling
follows from the existence of an $\rln$-coupling, as is proved in \ref
p.amen-average/; in the case of fixed amenable transitive graphs, this is a
well-known averaging principle.

\procl q.rtnpr
Does \ref t.rtnpr/ hold without the assumption of a unimodular
$\rln$-coupling, but just an $\rln$-coupling?
\endprocl

This question asks whether we can compare the traces from two different von
Neumann algebras.
One situation where we can do this is as follows.

If $\rtd_1$ and $\rtd_2$ are probability measures on $\GG_*$, then a
probability measure $\nu$
on $\GG_* \times \GG_*$ whose coordinate
marginals are $\rtd_1$ and $\rtd_2$ is called a {\bf monotone graph
coupling} of $\rtd_1$ and $\rtd_2$ if
$\nu$ is concentrated on pairs of rooted
networks $\big( (\gh_1, \bp), (\gh_2, \bp) \big)$
that share the same roots and satisfy $\vertex(\gh_1) \subseteq
\vertex(\gh_2)$.
In this instance,
let $V_1$ be the inclusion of $\ell^2\big(\vertex(\gh_1)\big)$ in
$\ell^2\big(\vertex(\gh_2)\big)$.

When there is a unimodular coupling (as in \ref t.rtnpr/),
the following result is easy.
The fact that it holds more generally is useful.

\comment{
\procl p.equality \procname{Trace Comparison}
Let $\nu$ be a monotone graph coupling of two unimodular probability measures
$\rtd_1$ and $\rtd_2$.
Let $T^{(i)} \in \alg$  ($i = 1, 2$).
Suppose that 
$$
T^{(1)}_{\gh_1} \oplus \bfz_2 \le T^{(2)}_{\gh_2}
\label e.cond1
$$
for $\nu$-almost all pairs $\big( (\gh_1, \bp), (\gh_2, \bp) \big)$.
Then $\tr_{\rtd_1}(T^{(1)}) \le \tr_{\rtd_2}(T^{(2)})$ with equality iff
$$
T^{(1)}_{\gh_1} \oplus \bfz_2 = T^{(2)}_{\gh_2} \quad \nu\hbox{-a.s.}
\label e.defeq1
$$
Similarly, if
$$
T^{(1)}_{\gh_1} \le V_1^* T^{(2)}_{\gh_2} V_1
\label e.cond2
$$
for $\nu$-almost all pairs $\big( (\gh_1, \bp), (\gh_2, \bp) \big)$,
then $\tr_{\rtd_1}(T^{(1)}) \le \tr_{\rtd_2}(T^{(2)})$ with equality iff
$$
T^{(1)}_{\gh_1} = V_1^* T^{(2)}_{\gh_2} V_1 \quad \nu\hbox{-a.s.}
\label e.defeq2
$$
\endprocl
}

\procl p.equality \procname{Trace Comparison}
Let $\nu$ be a monotone graph coupling of two unimodular probability measures
$\rtd_1$ and $\rtd_2$.
Let $T^{(i)} \in \alg(\rtd_i)$ be self-adjoint with
$$
T^{(1)}_{\gh_1} \le V_1^* T^{(2)}_{\gh_2} V_1
\label e.cond1
$$
for $\nu$-almost all pairs $\big( (\gh_1, \bp), (\gh_2, \bp) \big)$.
Then $\tr_{\rtd_1}(T^{(1)}) \le \tr_{\rtd_2}(T^{(2)})$.
If in addition
for $\nu$-almost all pairs $\big( (\gh_1, \bp), (\gh_2, \bp) \big)$ and for
all $x \in \vertex(\gh_1)$ we have
$$
\deg_{\gh_1}(x) < \deg_{\gh_2}(x)
\implies
\bigip{T^{(1)}_{\gh_1} x, x} < \bigip{T^{(2)}_{\gh_2} x, x}
\,,
\label e.cond2
$$
then either
$$
\tr_{\rtd_1}(T^{(1)}) < \tr_{\rtd_2}(T^{(2)})
\label e.concl1
$$
or 
$$
\vertex(\gh_1) = \vertex(\gh_2)\,,\ \edges(\gh_1) = \edges(\gh_2)\,,
\ \hbox{ and } \ T^{(1)}_{\gh_1} = T^{(2)}_{\gh_2} \quad \nu\hbox{-a.s.}
\label e.concl2
$$
\endprocl

\proof
Suppose that \ref e.cond1/ holds.
The fact that $\tr_{\rtd_1}(T^{(1)}) \le \tr_{\rtd_2}(T^{(2)})$ is an
immediate consequence of the definition of trace and of the hypothesis:
$$\eqaln{
\tr_{\rtd_1}(T^{(1)}) 
&=
\int \bigip{T^{(1)}_{\gh} \bp, \bp} \,d\rtd_1(\gh, \bp)
=
\int \bigip{T^{(1)}_{\gh_1} \bp, \bp} \,d\nu\big( (\gh_1, \bp), (\gh_2,
\bp) \big)
\cr&\le
\int
\bigip{V_1^* T^{(2)}_{\gh_2} V_1 \bp, \bp} \,d\nu\big( (\gh_1, \bp), (\gh_2,
\bp) \big)
\cr&=
\int \bigip{T^{(2)}_{\gh_2} \bp, \bp} \,d\nu\big( (\gh_1, \bp), (\gh_2,
\bp) \big)
=
\tr_{\rtd_2}(T^{(2)})
\,.
}$$

Suppose that equality holds in this inequality, i.e., \ref e.concl1/ fails.
Then
$$
\bigip{T^{(1)}_{\gh_1} \bp, \bp}
=
\bigip{T^{(2)}_{\gh_2} \bp, \bp}
\quad \nu\hbox{-a.s.}
$$
\comment{
Let $B_n(\bp)$ be the vertex ball of radius $n$ about $\bp$.
By the Mass-Transport Principle and the preceding equation, we have
$$\eqaln{
\int \sum_{x \in B_n(\bp)} \bigip{T^{(1)}_\gh x, x} \,d\rtd_1(\gh, \bp)
&=
\int |B_n(\bp)| \bigip{T^{(1)}_\gh \bp, \bp} \,d\rtd_1(\gh, \bp)
\cr&=
\int |B_n(\bp)| \bigip{T^{(2)}_\gh \bp, \bp} \,d\rtd_2(\gh, \bp)
\cr&=
\int \sum_{x \in B_n(\bp)} \bigip{T^{(2)}_\gh x, x} \,d\rtd_2(\gh, \bp)
\,.
}$$
Comparing the first and last quantities in this chain of equalities and using
the hypothesis again yields that
}
Of course, we also have by hypothesis that $\nu$-a.s.,
$$
\bigip{T^{(1)}_{\gh_1} x, x}
\le
\bigip{T^{(2)}_{\gh_2} x, x}
\label e.ineqatx
$$
for all $x \in \vertex(\gh_1)$.
Assume now that \ref e.cond2/ holds.
We shall prove that \ref e.concl2/ holds.
First,
we claim that $\nu$-a.s.,
$$
\vertex(\gh_1) = \vertex(\gh_2),\ 
\edges(\gh_1) = \edges(\gh_2)\,,
\label e.eqgraphs
$$
and
$$
\bigip{T^{(1)}_{\gh_1} x, x}
=
\bigip{T^{(2)}_{\gh_2} x, x}
\label e.atx
$$
for all $x \in \vertex(\gh_2)$.
If not, let $k$ be the smallest integer such that with positive
$\nu$-probability, there is a vertex $x$ at $\gh_1$-distance $k$ from $\bp$
where \ref e.atx/ does not hold.
Such a $k$ exists by virtue of \ref e.cond2/.
Consider $\bigip{T^{(i)}_{\gh_i} x, x}$ as part of the mark at $x$.
According to \ref t.rwrers/, the random walk on $\gh_i$ given in \ref
r.choice/ yields a shift-stationary measure $\Ph_i$ on trajectories
$\big((\gh_i, w_0), \Seq{w_n \st n \ge 0}\big)$ for each $i$.
In particular, the distribution of $\bigip{T^{(i)}_{\gh_i} w_k, w_k}$ is
the same as that of $\bigip{T^{(i)}_{\gh_i} w_0, w_0}$.
Now the latter is the same for $i=1$ as for $i=2$ (since $w_0 = \bp$).
Note that for all $x$ at distance less than $k$ from $\bp$, we have
$\deg_{\gh_1}(x) = \deg_{\gh_2}(x)$ by \ref e.cond2/.
Thus, the walks may be coupled together up to time $k$, whence
the distribution of $\bigip{T^{(i)}_{\gh_i} w_k, w_k}$ is
not the same for $i=1$ as for $i=2$ in light of \ref
e.ineqatx/ and choice of $k$.
This is a contradiction.
It follows that $\deg_{\gh_1}(x) = \deg_{\gh_2}(x)$ for all $x \in
\vertex(\gh_1)$, whence \ref e.eqgraphs/ and \ref e.atx/ hold.

Now $T := T^{(2)}_{\gh_2} - T^{(1)}_{\gh_1} \ge 0$ is self-adjoint.
It follows that
for any $x, y \in \vertex(\gh_2)$
and any complex number $\alpha$ of modulus 1,
$$
0
\le
\bigip{T (\alpha x + y), \alpha x + y}
=
2 \Re \big\{ \alpha \ip{T x, y} \big\}
\,,
$$
whence $\bigip{T x, y} = 0$.
That is, $T = 0$ and \ref e.concl2/ holds.
\Qed

\comment{
Since we also want to handle unbounded operators, sometimes arising as
unbounded functions of bounded operators, we recall that
a closed densely defined operator is {\bf affiliated}
with $\alg$ if it commutes with all unitary operators
that commute with $\alg$; see, e.g., \ref b.KadRing1/, p.~342.
Write $\affalg$ for the set of closed densely defined operators affiliated
with $\alg$.
Also note that our trace on $\alg$ is obviously finite.
Let $T \in \affalg$ be a self-adjoint operator 
with spectral resolution $E_T$. We define the Borel measure $\spm_{\rtd, T}$
by
$$
\spm_{\rtd, T}(B) :=
\tr_\rtd\big(E_T(B)\big)
\label e.specmsr
$$
for Borel subsets $B \subseteq \R$.
We extend the trace by defining 
$$
\tr_\rtd(T) := 
\int_0^\infty \lambda \,d\spm_{\rtd, T}(\lambda) 
$$
for positive operators $T \in \affalg$ and then by linearity to all of
$\affalg$ when it makes sense.
If $T \ge 0$ and $\tr(T) < \infty$, then $\bp \in \dmn(\sqrt T)$ since
$$\eqaln{
\lim_{n \to\infty} \|E_T[0, n] \sqrt T E_T[0, n] \bp \|^2 
&=
\lim_{n \to\infty} \bigip{E_T[0, n] T E_T[0, n] \bp, \bp}
\cr&=
\lim_{n \to\infty} \int_0^n \lambda \,d\spm_{\rtd, T}(\lambda) 
=
\tr_\rtd(T) < \infty
\,.
}$$
In particular, 
$$
\|\sqrt T \bp\|^2 = \tr_\rtd(T)
\,.
$$
Similarly, if $T$ is self-adjoint with polar decomposition $T = U |T|$ and
$\tr_\rtd\big(|T|\big) < \infty$,
then 
$$
\bigip{U \sqrt{|T|} \bp, \sqrt{|T|} \bp} = \tr_\rtd(T)
\,.
$$



We now extend \ref p.equality/.

\procl p.equalityunbdd \procname{Trace Comparison of Unbounded Operators}
Let $\nu$ be a monotone graph coupling of two unimodular probability measures
$\rtd_1$ and $\rtd_2$.
Let $T^{(i)} \in \affalg(\rtd_i)$ be self-adjoint with
$\tr_{\rtd_i} (|T^{(i)}|) < \infty$ ($i = 1, 2$).
If $\|\sqrt{|T^{(1)}_{\gh_1}|} \bp\| \le \|\sqrt{|T^{(2)}_{\gh_2}|} \bp\|$
for $\nu$-almost all pairs $\big( (\gh_1, \bp), (\gh_2, \bp) \big)$,
then $x \in \dmn(\sqrt{|T^{(i)}|})$ for all $x \in \vertex(\gh_i)$ ($i = 1,
2$) and $\tr_{\rtd_1}(T^{(1)}) \le \tr_{\rtd_2}(T^{(2)})$.
If in addition
for $\nu$-almost all pairs $\big( (\gh_1, \bp), (\gh_2, \bp) \big)$ and for
all $x \in \vertex(\gh_1)$ we have
$$
\deg_{\gh_1}(x) < \deg_{\gh_2}(x)
\implies
\|\sqrt{|T^{(1)}_{\gh_1}|} x\| < \|\sqrt{|T^{(2)}_{\gh_2}|} x\|
\,,
\label e.cond3
$$
then either
\ref e.concl1/ holds or
\ref e.concl2/ holds.
\endprocl

\proof
Note first that $e_n |T^{(2)}| e_n$ is a bounded operator, which is easily
seen to be in $\alg(\rtd_1)$.
Furthermore, this is an increasing sequence of operators that converges
since $L^1(\rtd_1)$ is complete by Theorem 13 of \ref b.Segal/.
\msnote{Or we might want to not use $\le$ in the hypo and use instead
sqrts. Maybe then we don't even need the limits.}
}

\comment{
\procl r.eq-ext
An extension of the previous proposition is proved similarly, where we
consider marks in $\marks^2$ and assume that $\nu$ is concentrated on pairs
of rooted networks such that when the second coordinate is forgotten, we
obtain identical rooted networks.
\msnote{Explain what this means}
\endprocl
}

\comment{
The following extension is useful as well (see \ref b.Lyons:est/).
A continuous function $f : (0, \infty) \to \R$ is called {\bf operator
monotone on $(0, \infty)$} if for any bounded self-adjoint operators $A, B$
with spectrum in $(0, \infty)$ and $A \le B$, we have $f(A) \le f(B)$.
According to L\"owner's theorem,
such $f$ are precisely those that can be represented as
$$
f(t) =
\alpha + \beta t + \int \left[{s \over 1+s^2} - {1 \over t +s}\right]
d\kappa(s)
\label e.pick
$$
for some $\alpha \in \R$, $\beta \ge 0$, and
some positive Borel measure $\kappa$ on $\CO{0, \infty}$ with $\int
d\kappa(s)/(1+s^2) < \infty$; see \ref b.Bhatia/, pp.~138--144 and \ref
b.Donoghue/.
Given $f$ as above, $0 \le T \in \alg$ and $\rtd \in \uni$, say that the
pair $(f, T)$ is {\bf mild} with respect to $\rtd$ if 
$$
\int 
\left[
\int_0^1 \bigip{(T+s)^{-1} \bp, \bp}
\,d\kappa(s)
+
\int_1^\infty \Big( {1 \over s} - \bigip{(T+s)^{-1} \bp, \bp} \Big)
\,d\kappa(s)
\right]
\,d\rtd(\gh, \bp)
< \infty
\,.
$$
We are interested in extending the equality condition of \ref p.equality/
to unbounded functions of unbounded operators.
For simplicity, we extend $\alg$ to the set of equivariant measurable maps
$T$ such that $T_{\gh}$ is a linear operator 
defined on the linear span of the vectors $x \in \vertex(\gh)$.

\procl l.incrlimit
Let $T_n$ be bounded self-adjoint operators on a Hilbert space $\HH$ and
$T$ an unbounded self-adjoint operator on $\HH$ with domain $\DD(T)$.
Suppose that $I \le T_n \le T_{n+1}$ for all $n$ and that for all $x \in
\DD(T)$, we have $\bigip{T_n x, x} \to \bigip{T x, x}$ as $n \to\infty$.
Then $\bigip{T_n^{-1} x, x} \to \bigip{T^{-1} x, x}$ as $n \to\infty$ for
all $x \in \HH$.
\endprocl

\proof
Note that $T$ is invertible by
Theorems 13.13 and 13.31 of \ref b.Rudin:FA/.
Since the function $t \mapsto -1/t$ is operator monotone on $(0, \infty)$,
we have that the pointwise limit $A$ of $T_n^{-1}$ exists and satisfies
$T_n^{-1} \ge A \ge T^{-1}$.
Therefore
$T_n \le A^{-1} \le T$.
(The proof of operator monotonicity on p.~114 of \ref b.Bhatia/ applies as
well to unbounded self-adjoint operators.)
On the other hand, the pointwise limit of $T_n$ is $T$ on $\DD(T)$, whence
$A = T$ on $\DD(T)$.
Since $A$ is also self-adjoint, it follows that its domain is the same as
$\DD(T)$ and so $A = T$.
\Qed

\msnote{Perhaps we want to do the previous theorem for unbounded operators
instead.}

\procl c.unbdd \procname{Traces of Unbounded Operators}
Let $\nu$ be a monotone graph coupling of two unimodular probability measures
$\rtd_1$ and $\rtd_2$.
Let $T^{(i, M)} \in \alg$ for $i = 1, 2$ and $M \in \N$ be such that
$$
T^{(1, M)}_{\gh_1} \oplus \bfz_2 \le T^{(2, M)}_{\gh_2} 
\label e.mono
$$
and
$$
0 \le T^{(i, M)}_{\gh_i} \le T^{(i, M+1)}_{\gh_i}
\label e.increasing
$$
for $i = 1, 2$, all $M \in \N$ and for
$\nu$-almost all pairs $\big( (\gh_1, \bp), (\gh_2, \bp) \big)$, where
$\bfz_2$ denotes the 0 operator on
$\ell^2\big(\vertex(\gh_2) \setminus \vertex(\gh_1)\big)$.
Assume that for all $\gh$ and all $x, y \in \vertex(\gh)$, we have
$\lim_{M \to\infty} \bigip{T^{(i, M)}_{\gh} x, y} = \bigip{T^{(i)}_{\gh} x,
y}$ for some $T^{(i)} \in \alg$ that has a self-adjoint
extension ($i = 1, 2$).
Let $f$ be an operator monotone function on $(0, \infty)$ with
$$
\lim_{t \downarrow 0} t f(t) = 0
\label e.nomass
$$
and such that $(f, T^{(i, 1)})$ is mild with respect to $\rtd_i$ for $i = 1,
2$.
If
$$
\lim_{M \to\infty} \lim_{\epsilon \downarrow 0}
\left[
\tr_{\rtd_2}\big(f(T^{(2, M)} + \epsilon)\big) -
\tr_{\rtd_1}\big(f(T^{(1, M)} + \epsilon)\big)
\right]
=
0
\label e.limit
$$
and
for $\nu$-almost all pairs $\big( (\gh_1, \bp), (\gh_2, \bp) \big)$ and for
all $x \in \vertex(\gh_1)$ we have
$$
\deg_{\gh_1}(x) < \deg_{\gh_2}(x)
\implies
\forall s > 0 \quad
\Bigip{\big(T^{(1)}_{\gh_1} + s\big)^{-1} x, x} <
\Bigip{\big(T^{(2)}_{\gh_2} + s\big)^{-1} x, x}
\,,
\label e.invs
$$
then $T^{(1)}_{\gh_1} = T^{(2)}_{\gh_2}$ $\nu$-a.s.
\endprocl

\proof
Write $f$ as in \ref e.pick/.
If $\kappa = 0$, then the result is trivial, so assume $\kappa \ne 0$.
The assumption \ref e.nomass/ ensures that $\kappa\big(\{0\}\big) = 0$.
For $\epsilon > 0$, $M \in \N$ and 
$\nu$-almost each pair $\big( (\gh_1, \bp), (\gh_2, \bp) \big)$, 
we have
$$\eqaln{
\bigip{f(T^{(2, M)}_{\gh_2} + \epsilon) \bp, \bp}
&-
\bigip{f(T^{(1, M)}_{\gh_1} + \epsilon) \bp, \bp}
\cr&\hskip-1in=
\beta \bigg[\bigip{T^{(2, M)}_{\gh_2} \bp, \bp}
-
\bigip{T^{(1, M)}_{\gh_1} \bp, \bp} \bigg]
\cr&\hskip-1in\qquad+
\int
\left[
\Bigip{ \big(T^{(1, M)}_{\gh_1} + \epsilon +s \big)^{-1} \bp, \bp} -
\Bigip{ \big(T^{(2, M)}_{\gh_2} + \epsilon +s \big)^{-1} \bp, \bp}
\right]
d\kappa(s)
\cr&\hskip-1in\ge
\int
\left[
\Bigip{ \big(T^{(1, M)}_{\gh_1} + \epsilon +s \big)^{-1} \bp, \bp} -
\Bigip{ \big(T^{(2, M)}_{\gh_2} + \epsilon +s \big)^{-1} \bp, \bp}
\right]
d\kappa(s)
\,,
\label e.upb
}
$$
where the inequality is from
\ref e.mono/ and the fact that $\beta \ge 0$.
By \ref e.increasing/,
we have for every $M \ge 1$, $\epsilon > 0$, and $s > 0$,
$$
\bigip{(T^{(i, M)}_\gh+\epsilon + s)^{-1} \bp, \bp}
\le
\bigip{(T^{(i, 1)}_\gh+s)^{-1} \bp, \bp}
\,,
$$
and so for $\epsilon < 1 < s$,
$$
{1 \over s} - \bigip{(T^{(i, M)}_\gh+\epsilon + s)^{-1} \bp, \bp} 
\le
{1 \over 1+s^2} + {1 \over 1+s} - \bigip{(T^{(i, 1)}_\gh+1+s)^{-1} \bp, \bp} 
\,.
$$
Because $(f, T^{(i, 1)})$ is mild with respect to $\rtd_i$, we may
therefore apply the Lebesgue
dominated convergence theorem and \ref l.incrlimit/ to obtain that
$$\eqaln{
\lim_{M \to\infty}
\lim_{\epsilon \downarrow 0}
&\int
\int_0^\infty
\left[
{s \over 1 + s^2} -
\Bigip{ \big(T^{(i, M)}_{\gh} + \epsilon +s \big)^{-1} \bp, \bp}
\right]
d\kappa(s)
\,d\rtd_i(\gh, \bp)
\cr&=
\int
\int_0^\infty
\left[
{s \over 1 + s^2} -
\Bigip{ \big(T^{(i)}_{\gh} + \epsilon +s \big)^{-1} \bp, \bp}
\right]
d\kappa(s)
\,d\rtd_i(\gh, \bp)
\,.
}$$
Combining this with \ref e.limit/ and \ref e.upb/ yields
$$
\int
\int_0^\infty
\left[
\Bigip{ \big(T^{(1)}_{\gh_1} +s \big)^{-1} \bp, \bp} -
\Bigip{ \big(T^{(2)}_{\gh_2} +s \big)^{-1} \bp, \bp}
\right]
d\kappa(s)
\,d\nu\big( (\gh_1, \bp), (\gh_2, \bp) \big)
=
0
\,.
$$
Now the integrand is non-negative, whence
for $\kappa$-almost all $s > 0$ and
for $\nu$-almost all pairs $\big( (\gh_1, \bp), (\gh_2, \bp) \big)$,
$$
\Bigip{ \big(T^{(1)}_{\gh_1} +s \big)^{-1} \bp, \bp}
=
\Bigip{ \big(T^{(2)}_{\gh_2} +s \big)^{-1} \bp, \bp}
\,.
$$
Note that the operators $\big(T^{(i)}_{\gh_i} +s \big)^{-1}$ are
bounded with norm at most $s^{-1}$.
Therefore, we may conclude from \ref p.equality/ that these operators are
equal $\nu$-a.s.
Now Theorems 13.13 and 13.31 of \ref b.Rudin:FA/ tell us that
$T^{(i)}_{\gh} + s$ is a bijection of its domain with
$\ell^2\big(\vertex(\gh)\big)$.
Hence $T^{(1)}_{\gh_1} + s = T^{(2)}_{\gh_2} + s$ $\nu$-a.s.,
which gives us the desired result.
\Qed
}

\bsection{Percolation}{s.extend}

We now begin our collection of
extensions of results that are known for unimodular fixed graphs.
For most of the remainder of the paper, we consider graphs without marks, or,
equivalently, with constant marks, except that marks are used as explained
below to perform percolation on the given graphs.
We begin this section on percolation with some preliminary results on
expected degree.

\procl t.deg-infinite \procname{Minimal Expected Degree}
If $\rtd$ is a unimodular probability measure on $\GG_*$ concentrated on
infinite graphs, then $\expdeg(\rtd) \ge 2$.
\endprocl

This is proved exactly like Theorem 6.1 of \BLPSgip\ is proved.
In the context of equivalence relations, this is well known and was perhaps
first proved by \ref b.Levitt/.

\procl t.deg2 \procname{Degree Two}
If $\rtd$ is a unimodular probability measure on $\GG_*$ concentrated on
infinite graphs,
then $\expdeg (\rtd) = 2$ iff $\rtd$-a.s.\ $\gh$ is a tree with at most 2 ends.
\endprocl

The proof is like that of Theorem 7.2 of \BLPSgip.

\procl p.limittree \procname{Limits of Trees}
If $G_n$ are finite trees with random weak limit $\rtd$, then
$\expdeg(\rtd) \le 2$ and $\rtd$ is concentrated on trees with at most 2
ends.
\endprocl

\proof
Since $\expdeg\big(U(G_n)\big) < 2$, we have $\expdeg(\rtd) \le 2$.
The remainder follows from \ref t.deg2/.
\Qed

We now discuss what we mean by percolation on a random rooted network.
Given a probability measure $\rtd \in \uni$, we may wish to randomly
designate some of the edges of the random network ``open".
For example, in Bernoulli($p$) bond percolation, each edge is independently
open with probability $p$.
More generally, we'd like to couple together all these Bernoulli
percolation measures. We do this by using the canonical networks.
We wish the second coordinates to be uniformly distributed on $[0, 1]$ and 
independent (but the same at each endpoint of a given edge).
For $0 \le i < j$, let $U_{i, j}$ be i.i.d.\ uniform $[0, 1]$ random
variables.
Then for each canonical network $(\gh, 0) \in \GG_*$ and for each $0 \le i
< j$, change the mark at each endpoint of the edge between $i$ and $j$, if
there is an edge, by adjoining a second coordinate equal to $U_{i, j}$.
Let $\rtd^\bn$ be the law of the resulting network class when $[\gh, 0]$
has law $\rtd$.
It is clear that $\rtd^\bn$ is unimodular when $\rtd$ is.
We refer to $\rtd^\bn$ as the {\bf standard coupling of
Bernoulli percolation on} $\rtd$.
For $p \in [0, 1]$, one can then define {\bf Bernoulli($p$) bond
percolation on} $\rtd$
as the measure $\rtd^\bn_p$ that replaces the second coordinate of each edge
mark by ``open" if it is at most $p$ and by ``closed" otherwise.
In the future, we shall not be explicit about how randomness is added to
random networks. 

Every map $\mkmp : \marks \to \marks$ induces a map on $\GG_*$ by
applying $\mkmp$ to all the marks of a network.
For simplicity, we shall denote this induced map still by $\mkmp$.
Note that if $\phi : \marks \times [0, 1] \to \marks$ is the projection
onto the first coordinate, then $\rtd = \rtd^\bn_p \circ \phi^{-1}$.
We have changed the mark space, but a fixed homeomorphism would bring it back
to $\marks$.
Thus, more generally, if $\psi : \marks \to \marks$ is Borel, then we call
$\rtd$ a {\bf percolation} on $\rtd \circ \psi^{-1}$.

\procl d.ins
Let $G=\big(\verts(G),\edges(G)\big)$ be a graph.
Given a configuration $A\in \two^{\edges(G)}$ and an edge
$e\in\edges(G)$, denote $\ins_e A$ the element of $\two^{\edges(G)}$ that
agrees with $A$ off of $e$ and is $1$ on $e$.
For $\A\subset \two^{\edges(G)}$, we write $\ins_e\ev A := \{\ins_e A \st A \in
\A\}$.
For bond percolation, call an edge ``closed" if it is marked ``0" and ``open"
if it is marked ``1".
A bond percolation process $\P$ on $G$ is {\bf insertion tolerant} if
$\P(\ins_e\ev A)>0$
for every $e\in\edges(G)$ and every Borel $\A\subset \two^{\edges(G)}$
satisfying $\P(\A)>0$.
The primary subtlety in extending this notion to percolation on unimodular
random networks is that it may not be possible to pick an edge measurably from
a rooted-automorphism-invariant set of edges.
Thus, we shall make an extra assumption of distinguishability with marks.
That is, a percolation process $\P$ on a unimodular probability measure on
$\GG_*$ is {\bf insertion tolerant} if $\P$-a.s.\ there is no nontrivial
rooted isomorphism of the marked network and
for any event $\A \subseteq \{ (A,\gh)
\st A \subseteq \edges(\gh),\; \gh \in \GG_* \}$ with $\P(\A) > 0$ and any
Borel function $e : \gh \mapsto e(\gh) \in \edges(\gh)$ defined on
$\GG_*$, we have $\P(\ins_e \A) > 0$.
\endprocl

For example, Bernoulli($p$) bond percolation
is insertion tolerant when $p\in \OC{0, 1}$.

We call a connected component of open edges (and their endpoints) a {\bf
cluster}.
Given a rooted graph $(\gh, \bp)$, define
$$
\pc(\gh, \bp) := \sup \{ p \st \hbox{Bernoulli}(p)\hbox{ percolation on }
\gh \hbox{ has no infinite clusters a.s.}\}
\,.
$$
Clearly $\pc$ is $\invar$-measurable, so if $\rtd \in \uni$ is extremal,
then there is a constant $\pc(\rtd)$ such that $\pc(\gh, \bp) = \pc(\rtd)$
for $\rtd$-a.e.\ $(\gh, \bp)$.

\procl x.pclower
Even if $\rtd \in \uni$ satisfies $\expdeg(\rtd) < \infty$, it does
not necessarily follow that $\pc(\gh) > 0$ for $\rtd$-a.e.\ $(\gh, \bp)$.
For example, let $p_k := 1/[k(k+1)]$ for $k \ge 1$ and $p_0 := 0$.
Let $\UGW$ be the corresponding unimodular Galton-Watson measure (see
\ref x.AGW/).
Then $\expdeg(\UGW) = 6/(12-\pi^2)$ by \ref e.degAGW/,
but since $\sum k p_k = \infty$,
we have $\pc(\gh) = 0$ a.s.\ by \ref b.Lyons:rwpt/.
\endprocl

A more elaborate example
shows that no stochastic bound on the
degree of the root,
other than uniform boundedness,
implies $\pc(\rtd) > 0$:

\procl x.pclower2
Given $a_n > 0$, we shall construct $\rtd \in \uni$ such that
$\rtd[\deg_\gh(\bp) \ge n] < a_n$ for all large $n$
and $\pc(\gh) = 0$ for $\rtd$-a.e.\ $(\gh, \bp)$.
We may assume that $\sum a_n < \infty$.
Consider the infinite tree $T$ of degree $3$ and $\bp \in \vertex(T)$.
Let $j \geq 2$ and set $p_j := \frac{1}{2} + \frac{1}{j}$. 
Consider supercritical Bernoulli($p_j$)
bond percolation on $T$,
so that
$$
\theta_j := \P[\bp
\hbox{ belongs to an infinite component}] > 0
\,.
$$
Define $q_j < 1$ by
$p_j q_j = \frac{1}{2} + \frac{1}{2j}$,
so that the fragmentation of an infinite $p_j$-component by an independent
$q_j$ percolation process will still contain infinite components.
Let $N_j$ be the smallest integer with
$(1 - \frac{1}{j^2})^{N_j} < 1 - q_j$.
Finally, choose some sequence $1 > r_j \downarrow 0$ sufficiently fast.
We label the edges of $T$ by the following
operations, performed independently for each $j \geq 2$.

Take the infinite components of Bernoulli($p_j$) bond percolation on $T$.
``Thin" by retaining each component independently with probability $r_j$ and
deleting other components.
For each edge $e$ in the remaining components, let $L_j(e) := N_j$, while
$L_j(e) := 1$ for deleted edges $e$.

Now let $L(e) := \sup_j L_j(e)$.
Since $r_j \to 0$ fast, $L(e) < \infty$ for all $e$ a.s.
Consider the graph $G$ obtained by replacing each edge $e$ of $T$ by $L(e)$
parallel chains of length 2.
To estimate $\deg_G(\bp)$, note that the chance that $\bp$ is incident
to an edge $e$ with $L_j(e) = N_j$ (thus contributing at most $3N_j$ to the
degree) equals $r_j\theta_j$.  
Thus, by choice of $\Seq{r_j}$, we can make the root-degree distribution
have tail probabilities eventually less than $a_n$.
Now consider Bernoulli($1/j$) bond percolation on $G$.
For an edge $e$ of $T$ which is replaced by $N_j$ chains of $G$,
the chance of percolating across at least one of these $N_j$
chains is (by definition of $N_j$) larger than $q_j$.
Thus the percolation clusters on $G$ dominate the $q_j$-percolation
clusters on the retained components of the original infinite
$p_j$-percolation clusters on $T$, and as observed above must therefore
contain infinite components.

To see how to make this into a probability measure in $\uni$, note that
$L$ is an invariant random network on $T$.
Therefore, we obtain a probability measure in
$\uni$ by \ref t.uni-vs-invar/.
We may now use the edge labels to replace an edge labeled $n$ by $n$ parallel
chains of length 2, followed by
a suitable re-rooting
as in \ref x.stretched/ (below).
This gives a new probability measure in $\uni$ that has
the property desired, since 
the expected degree of the root is finite by the hypothesis $\sum_n a_n <
\infty$.
Also, the re-rooting stochastically decreases the degree since it
introduces roots of degree 2.
\endprocl

The following extends a well-known result of \ref b.HP:uniq/.
The proof is the same.

\procl t.uniq \procname{Uniqueness Monotonicity and Merging Clusters}
Let $\rtd$ be a unimodular probability measure on $\GG_*$.
Let $p_1 < p_2$ and $\P_i$ ($i=1,2$) be the corresponding Bernoulli($p_i$) bond
percolation processes on $\rtd$.
If there is a unique infinite cluster $\P_1$-a.s., then
there is a unique infinite cluster $\P_2$-a.s.  Furthermore, in the
standard coupling of Bernoulli percolation processes, if $\rtd$ is extremal,
then $\rtd$-a.s.\ for all $p_1, p_2$ satisfying $\pc(\rtd) < p_1 < p_2 \le 1$,
every infinite $p_2$-cluster contains an infinite $p_1$-cluster.
\endprocl

As a consequence, for extremal $\rtd \in \uni$,
there is a constant $\pu(\rtd)$ such that for any
$p > \pu(\rtd)$, we have $\P_{p}$-a.s., there is a
unique infinite cluster, while for any $p < \pu(\rtd)$,
we have $\P_{p}$-a.s., there is not a unique infinite cluster.

Every unimodular probability measure on $\GG_*$ can be written as a Choquet
integral of extremal measures.
In the following, we refer to these extremal measures as ``extremal
components".

\procl l.ergcomp
If\/ $\P$ is an insertion-tolerant percolation on a
unimodular random network that is concentrated on infinite graphs, then
almost every extremal component of\/ $\P$ is insertion tolerant.
\endprocl

\proof
The proof can be done precisely as that of Lemma 1 of \ref b.GKN/, but
by using the present
Theorems \briefref t.rwrers/, \briefref t.erg/, and \briefref
t.findroot/, as well as \ref r.choice/ to replace the use of a
measure-preserving transformation in \ref b.GKN/ by the shift on
trajectories of a stationary Markov chain.
The fact that $\P$ is concentrated on infinite graphs is used to deduce
that the corresponding Markov chain is not positive recurrent, whence the
asymptotic frequency of visits to any given neighborhood of the root is 0.
\Qed

\procl c.ninf \procname{Number of Infinite Clusters}
If\/ $\P$ is an insertion-tolerant percolation on a
unimodular random network, then
$\P$-almost surely, the number of infinite clusters is $0,1$ or $\infty$.
\endprocl

The proof is standard for the extremal components; cf.\ \ref b.NS:number/.

The following extends \ref b.HP:uniq/ and Proposition 3.9 of \ref
b.LS:indist/.
The proof is parallel to that of the latter.

\procl p.niso
\assumptions.
Then $\P$-a.s.\ each infinite cluster that has at least $3$ ends has
no isolated ends.
\endprocl

The following corollary is proved just like Proposition 3.10 of \ref
b.LS:indist/.
There is some overlap with Theorem 3.1 of \ref b.Paulin/.

\procl c.isolated \procname{Many Ends}
\genassumptions.
If there are infinitely many infinite clusters $\P$-a.s., then $\P$-a.s.\ every
infinite cluster has continuum many ends and no isolated end.
\endprocl

The following extends Lemma 7.4 and Remark 7.3
of \BLPSgip\ and is proved similarly.

\procl l.trim \procname{Subforests}
\assumptions.
If\/ $\P$-a.s.\ there is a component of the open subgraph
$\omega$ with at least three ends, then
there is a percolation $\fo$ on $\omega$ whose components are trees
such that
a.s.\ whenever a component $K$ of $\omega$ has at least three ends,
there is a component of $K \cap \fo$ that has infinitely many ends and has
$\pc < 1$.
\endprocl

The following extends Proposition 3.11 of \ref b.LS:indist/ and is proved
similarly (using the preceding \ref l.trim/).

\procl p.trans \procname{Transient Subtrees}
\genassumptions.
If there are $\P$-almost surely infinitely many
infinite clusters,  then $\P$-a.s.\ each infinite cluster is transient and,
in fact, contains a transient tree.
\endprocl

In order to use this, we shall use the comparison of simple to network
random walks given in \ref p.transtrees/.

\procl d.inds
A percolation process $\P$ on a unimodular probability measure on
$\GG_*$ has {\bf indistinguishable infinite clusters} if
for any event $\A \subseteq \{ \big(A,(\gh, \bp)\big)
\st A \in \two^{\vertices(\gh)} \times \two^{\edges(\gh)},\;
(\gh, \bp) \in \GG_* \}$ that is invariant under non-rooted isomorphisms,
almost surely, for all infinite clusters $C$ of the open subgraph
$\omega$, we have $(C,\omega)\in\ev A$, or for all
infinite clusters $C$, we have $(C,\omega)\notin\ev A$.
\endprocl

The following extends Theorem 3.3 of \ref b.LS:indist/ and is proved
similarly using the preceding results:
For example, instead of the use of delayed simple random walk by \ref
b.LS:indist/, we use the network random walk in \ref r.choice/. 
This is a reversible random walk corresponding to edge weights $(x, y)
\mapsto 1/[(\deg x)(\deg y)]$.
It is transient by Propositions \briefref p.trans/ and \briefref
p.transtrees/, combined with Rayleigh's monotonicity principle.

\procl t.cerg \procname{Indistinguishable Clusters}
If\/ $\P$ is an insertion-tolerant percolation on a
unimodular random network, then
$\P$ has indistinguishable infinite clusters.
\endprocl

Among the several consequences of this result is the following extension of
Theorem 4.1 of \ref b.LS:indist/, proved similarly.

\procl t.tau \procname{Uniqueness and Long-Range Order}
\genassumptions, $\rtd$.
If\/ $\P$ is extremal and there is more than one infinite cluster $\P$-a.s.,
then $\rtd$-a.s.,
$$
\inf \big \{ \P[\hbox{there is an open path from $x$ to } y \mid \gh]
\st x, y \in \vertex(\gh) \big \} = 0
\,.
$$
\endprocl

The following extends Theorem 6.12 of \ref b.LS:indist/ and is proved
similarly.

\procl t.product \procname{Uniqueness in Products}
Suppose that $\rtd$, $\rtd_1$, and $\rtd_2$ are extremal unimodular probability
measures on $\GG_*$, with $\rtd$ supported on infinite graphs and $\rtd_1$
a percolation on $\rtd_2$. Then $\pu(\rtd \itimes \rtd_1) \ge
\pu(\rtd \itimes \rtd_2)$. In particular, $\pu(\rtd) \ge
\pu(\rtd \itimes \rtd_2)$.
\endprocl

More results on percolation will be presented in \ref s.amen/.

\bsection{Spanning Forests}{s.span}

An interesting type of percolation other than Bernoulli is given by certain
random forests.
There are two classes of such random forests that have been widely studied,
the uniform ones and the minimal ones.

We first discuss the uniform case.
Given a finite connected graph, $\gh$, let $\ust(\gh)$ denote the uniform
measure on spanning trees on $\gh$.
\ref b.Pemantle:ust/
proved a conjecture of Lyons, namely,
that if an infinite connected graph $\gh$ is exhausted by a sequence of finite
connected subgraphs $\gh_n$, then the weak limit of
$\Seq{\ust(\gh_n)}$ exists.
However, it may happen that the limit measure is not supported on trees,
but on forests.
This limit measure is now called the {\bf free (uniform) spanning forest} on
$\gh$, denoted $\fsf$ or $\FUSF$.
If $\gh$ is itself a tree, then this measure is
trivial, namely, it is concentrated on $\{\gh\}$. Therefore, \ref
b.Hag:umesft/
introduced another limit that had been considered on $\Z^d$
more implicitly by \ref b.Pemantle:ust/ and explicitly by
\ref b.Hag:rcust/, namely, the weak limit of the uniform
spanning tree measures on $\gh_n^*$, where $\gh_n^*$ is the graph $\gh_n$ with
its boundary identified (``wired")
to a single vertex. As \ref b.Pemantle:ust/ showed,
this limit also always exists
on any graph and is now called the {\bf wired (uniform) spanning forest},
denoted $\wsf$ or $\WUSF$.
It is clear that both $\fsf$ and $\wsf$ are concentrated on the set of
spanning forests\ftnote{*}{By a
``spanning forest'', we mean a subgraph without cycles that contains every
vertex.} of $\gh$ that are {\bf essential}, meaning that all their
trees are infinite.
Both $\fsf$ and $\wsf$ are important in their own right; see \ref
b.Lyons:bird/ for a survey and \BLPSusf\ for a comprehensive
treatment.

In all the above, one may work more generally with a weighted graph, where
the graph has positive weights on its edges.
In that case, $\ust$ stands for the measure such that the probability of a
spanning tree is proportional to the product of the weights of its edges.
The above theorems continue to hold and we use the same notation for the
limiting measures.

Most results known about the uniform spanning forest measures hold for
general graphs. Some, however, require extra hypotheses such as
transitivity and unimodularity.
We extend some of these latter results here.

Given $\rtd$, taking the wired uniform spanning forest on each graph
gives a percolation that we denote $\WUSF(\rtd)$.
Our first result shows, among other things,
that the kind of limit considered in this paper,
i.e., random weak convergence, gives another natural way to define $\wsf$.
It might be quite useful to have an explicit description of measures on
finite graphs whose random weak limit is the {\it free\/} spanning forest.

\procl p.ust-limit \procname{UST Limits}
If $\rtd$ is a unimodular probability measure on infinite networks in
$\GG_*$, then $\expdeg\big(\WUSF(\rtd)\big) = 2$.
If $\gh_n$ are finite connected networks whose random weak limit is $\rtd$,
then $\UST(\gh_n) \cd \WUSF(\rtd)$.
More generally, if $\rtd_n$ are unimodular probability measures on $\GG_*$
with $\rtd_n \cd \rtd$, then $\WUSF(\rtd_n) \cd \WUSF(\rtd)$.
\endprocl

\proof
We begin by proving part of the third sentence, namely, 
$$
\hbox{every weak limit point of $\Seq{\WUSF(\rtd_n)}$ 
stochastically dominates $\WUSF(\rtd)$}
\,.
\eqno(*)
$$
Given a positive integer $R$, let $\UST_R(\rtd)$ be the uniform spanning
tree on the wired ball of radius $R$ about the root.
(Although $\UST_R(\rtd) \notin \uni$, this will not affect our argument.)
Identify the edges of the wired ball of radius $R$ with the edges of the
ball itself.
By definition, we have $\UST_R(\rtd) \cd \WUSF(\rtd)$ as $R \to\infty$.
Clearly, $\UST_R(\rtd_n) \cd \UST_R(\rtd)$ as $n \to\infty$.
Furthermore, the intersection of $\WUSF(\rtd_n)$ with the ball of radius $R$
stochastically dominates $\UST_R(\rtd_n)$ by a theorem of \ref
b.FedMih/.
Therefore, every weak limit point of $\Seq{\WUSF(\rtd_n)}$ stochastically
dominates $\UST_R(\rtd)$ and therefore also $\WUSF(\rtd)$.

Suppose now that $\rtd$ is concentrated on recurrent networks.
If $\rtd$ is concentrated on networks with bounded degree, then so is
$\WUSF(\rtd)$, and the latter is also concentrated on recurrent networks by
Rayleigh's monotonicity principle. By \ref p.treespeed/, the claim of the first
sentence follows.
If $\rtd$ has unbounded degree, then let $\rtd_n$ be the law of the
component of the root when all edges incident to vertices of degree larger
than $n$ are deleted. Clearly $\rtd_n \in \uni$ and $\rtd_n \cd \rtd$.
We have shown that $\expdeg(\WUSF(\rtd_n)) = 2$, so that $(*)$ and
%
Fatou's lemma yield that $\expdeg(\WUSF(\rtd)) \le 2$,
whence equality results from \ref t.deg-infinite/.

\comment{
On the other hand, from \ref p.limittree/,
every weak limit point of $\Seq{\UST(\rtd_n)}$
has expected degree at most 2, while by \ref t.deg-infinite/, we have
$\expdeg\big(\WUSF(\rtd)\big) \ge 2$.
It follows that $\UST(\rtd_n) \cd \WUSF(\rtd)$ as $n \to\infty$ and that
$\expdeg\big(\WUSF(\rtd)\big) = 2$.

By \ref t.TII/, it follows that
$\expdeg\big(\WUSF(\rtd)\big) = 2$ for all $\rtd \in \uni$.
}%

Suppose next that $\rtd$ is concentrated on transient networks.
Then the proof of Theorem 6.5 of \BLPSusf\ gives the same result.

Finally, if $\rtd$ is concentrated on neither recurrent nor transient
networks, then we may write $\rtd$ as a mixture of two unimodular measures
that are concentrated on recurrent or on transient networks and apply the
preceding.

This proves the first sentence. The second sentence is a special case of
the third, so it remains to finish the proof of the third.
By Fatou's lemma and \ref t.deg-infinite/, after what we have shown, we
know that
all weak limits of $\Seq{\WUSF(\rtd_n)}$
have expected degree 2, as does $\WUSF(\rtd)$. Since all such
weak limits lie in $\uni$, $(*)$ shows that
all weak limits are equal.
\Qed

We next show that the trees of the $\wsf$ have only one end a.s.
The first theorem of this type was proved by \ref b.Pem:ust/.
His result was completed and extended in
Theorem 10.1 of \BLPSusf, which
dealt with the transitive unimodular case.
The minor modifications needed for the quasi-transitive unimodular case
were explained by \ref b.Lyons:est/.
Another extension is given by \ref b.LMS:wsf/, who showed that for graphs with
a ``reasonable" isoperimetric profile, each tree has only one end
$\wsf$-a.s.

\procl t.1endUSF \procname{One End}
If $\rtd$ is a unimodular probability measure on $\GG_*$ that is
concentrated on transient networks with bounded degree,
then $\WUSF(\rtd)$-a.s., each tree has exactly one end.
\endprocl

\proof
The proof is essentially the same as in \BLPSusf, with the following
modifications.
In the proof of Theorem 10.3 of \BLPSusf, which is the case where there is
only one tree a.s., we replace $x$ and $y$ there by $\bp$ and $X_n$, where
$\Seq{X_n}$ is the simple random walk starting from the root; to use
stationarity, we bias the underlying network by the degree of the root.
This gives a measure equivalent to $\rtd$, so that almost sure conclusions for
it hold for $\rtd$ as well.
The stationarity and reversibility give that the probability that the
random walk from $\bp$ ever visits $X_n$ is equal to the probability that a
random walk from $X_n$ ever visits $\bp$.
By transience, this tends to 0 as $n \to\infty$, which allows the proof of
\BLPSusf\ to go through.
[Here, we needed finite
expected degree to talk about probability since we used the equivalent
probability measure of biasing by the degree.]

In the proof of Theorem 10.4 of \BLPSusf, the case when there is more than
one tree a.s.,
we need the degrees to be bounded for the displayed equality on p.~36 of
\BLPSusf\ to hold up to a constant factor.
\Qed

For our proof, we had to assume transience;
there is presumably an extension to the recurrent case, which would say
that the number of ends $\WUSF(\rtd)$-a.s.\ is the same as the number of
ends $\rtd$-a.s.\ when $\rtd$ is concentrated on recurrent networks.
Also, presumably the assumption that the degrees are bounded is not needed.
In any case, our result here goes beyond what has been done before and
gives a partial answer to Question 15.4 of \BLPSusf; removing the
assumption of bounded degrees would completely answer that question.
It also applies, e.g., to transient clusters of Bernoulli percolation; see
\ref b.GKZ/ and \ref b.BLS:pert/ for sufficient conditions for transience.

We now prove analogous results for the other model of spanning trees, the
minimal ones.
Given a finite connected graph, $\gh$, and independent uniform $[0, 1]$
random variables on its edges, the spanning tree that minimizes the sum of
the labels of its edges has a distribution denoted $\MST(\gh)$, the {\bf
minimal spanning tree} measure on $\gh$.
If $\gh$ is infinite, there are two analogous measures, as in the uniform
case. They can be defined by weak limits, but also directly (and by
pointwise limits).
Namely, given independent uniform $[0, 1]$ edge labels, remove all edges
whose label is the largest in some cycle containing that edge.
The remaining edges form the {\bf free minimal spanning forest},
$\FMSF$.
If one also removes the edges $e$ both of whose endpoints belong to infinite
paths of edges that are all labeled smaller than $e$ is, then the resulting
forest is called the {\bf wired minimal spanning forest}, $\WMSF$.

The following is analogous to \ref p.ust-limit/ above and is proved
similarly to it and part of Theorem 3.12 of \LPSmsf, using \ref t.death/
below. Parts of it were also proved by \ref b.AS:obj/.

\procl p.mst-limit \procname{MST Limits}
If $\rtd$ is a unimodular probability measure on infinite networks in
$\GG_*$, then $\expdeg\big(\WMSF(\rtd)\big) = 2$.
If $\gh_n$ are finite connected networks whose random weak limit is $\rtd$,
then $\MST(\gh_n) \cd \WMSF(\rtd)$.
More generally, if $\rtd_n$ are unimodular probability measures on $\GG_*$
with $\rtd_n \cd \rtd$, then $\WMSF(\rtd_n) \cd \WMSF(\rtd)$.
\endprocl

The following extends a result of \LPSmsf, which in turn extends a result
of \ref b.Alexander:MSF/, who proved this in fixed Euclidean lattices.
Our proof follows slightly different lines.
For information on when the hypothesis is satisfied, see \ref t.death/
below.

\procl t.msf1end \procname{One End}
If $\rtd$ is an extremal
unimodular probability measure on infinite networks in
$\GG_*$ and there is $\P_{\pc(\rtd)}$-a.s.\ no infinite cluster,
then $\WMSF(\rtd)$-a.s., each tree has exactly one end.
\endprocl

\proof
By the first part of \ref p.mst-limit/ and by \ref t.deg2/, each tree has at
most 2 ends $\WMSF(\rtd)$-a.s.
Suppose that some tree has 2 ends with positive probability.
A tree with precisely two ends has a {\bf trunk}, the unique bi-infinite
simple path in the tree.
By the definition of $\WMSF$, the labels on a trunk cannot have a maximum.
By the Mass-Transport Principle, the limsup in one direction must equal the
limsup in the other direction, since otherwise we could identify the one
edge that has label larger than the average of the two limsups and is the
last such edge in the direction from the larger limsup to the smaller
limsup.
Let $p$ be this common limsup.
By the preceding, all the labels on the trunk are strictly less than $p$.
The root belongs to the trunk with positive probability.
Assume this happens.
Then the root belongs to an infinite $p$-cluster (the one containing the
trunk).
Now invasion from the root will fill (the vertices of)
this $p$-cluster and is part of the
tree containing the root (see \LPSmsf), whence the tree contains (the
vertices of) the entire $p$-cluster of the root.
By \ref t.uniq/, the tree therefore also contains an infinite $p'$-cluster
for every $p' \in (\pc(\rtd), p)$ if $p > \pc(\rtd)$.
Let $x$ be a vertex in the tree that is in an infinite $p'$-cluster $C$ for
$p' := (\pc(\rtd) + p)/2$.
Now invasion from $x$ has a finite symmetric difference with
invasion from $o$ (see \LPSmsf),
invasion from $x$ does not leave $C$, and invasion from $\bp$ fills the
trunk.
It follows from the definition of $p$ that $p' \ge p$.
That is, $p = \pc(\rtd)$.
Therefore, there is $\P_{\pc(\rtd)}$-a.s.\ an infinite cluster.
\Qed

\bsection{Amenability and Nonamenability}{s.amen}

Recall that a graph $\gh$ is (vertex) amenable
iff there is a sequence of subsets $H_n\subset \vertex(\gh)$
with
$$
\lim_{n \to\infty} {|\bdv H_n| \over |\verts(H_n)|} = 0\,,
$$
where $|\cbuldot|$ denotes cardinality.

Amenability, originally defined for groups, now appears in several areas of
mathematics, including probability theory and ergodic theory.
Its presence provides many tools one is used to from $\Z$ actions, yet its
absence also provides a powerful threshold principle.
There are many equivalent definitions of amenability.
We choose one that is not standard, but is useful for our probabilistic
purposes.
We show that it is equivalent to other definitions.
Then we shall illustrate its uses.

\procl d.amen
Let $\projection : \marks \to \marks$ be the composition of a homeomorphism of
$\marks$ with $\marks^2$ followed by the projection onto the first coordinate.
If a rooted network $(\gh, \bp)$ is understood,
then for a subset $\marks_0 \subseteq \marks$ and vertex $x$, the {\bf
$\marks_0$-component of $x$} is the set of vertices that can be reached
from $x$ by edges both of whose marks lie in $\marks_0$.
Write $K(\marks_0)$ for the $\marks_0$-component of the root.
For a probability measure $\rtd$ on rooted graphs, denote by $\fdom(\rtd)$ the
class of percolations on $\rtd$ that have only finite components.
That is, $\fdom(\rtd)$ consists of pairs $(\nu, \marks_0)$ such that $\nu$ is
a unimodular probability measure on $\GG_*$, $\marks_0 \subseteq \marks$ is
Borel, $\rtd = \nu \circ \projection^{-1}$, and
$K(\marks_0)$ is finite $\nu$-a.s.
(By \ref l.unmark/, all $\marks_0$-components are then finite $\nu$-a.s.)
For $\marks_0 \subseteq \marks$ and $x \in \vertex(\gh)$, write
$$
\ct(x, \marks_0) := | \{y \in \vertex(\gh) \st (x, y) \in \edges(\gh),
\hbox{ some edge mark of } (x, y) \hbox{ is }\notin \marks_0 \} |
\,.
$$
For $K \subset \vertex(\gh)$, define
$$
\ct(K, \marks_0) :=
\sum_{x \in K} \ct(x, \marks_0)
\,.
$$
Define
$$
\isoe(\rtd)
:=
\inf \left\{ \int {\ct\big(K(\marks_0), \marks_0\big) \over |K(\marks_0)|}
d\rtd'(\gh, \bp)
\st (\rtd', \marks_0) \in \fdom(\rtd) \right\}
\,.
$$
Call $\rtd$ {\bf amenable} if $\isoe(\rtd) = 0$.
Define 
$$
\expdeg(\rtd', \marks_0)
:=
\int \big[\deg_\gh(\bp) - n(\bp, \marks_0)\big] \,d\rtd'(\gh, \bp)
\,,
$$
the expected degree in the $\marks_0$-component of the root,
and
$$
\alpha(\rtd)
:=
\sup \big\{ \expdeg(\rtd', \marks_0) \st (\rtd', \marks_0) \in \fdom(\rtd)
\big\}
\,.
$$
Of course, neither $\isoe(\rtd)$ nor $\alpha(\rtd)$ depends on the choice
of homeomorphism in $\projection$.
Furthermore, these quantities depend only on the probability measure on
the underlying graphs of the networks, not on the marks.
\endprocl

This definition of amenability is justified in three ways: It agrees with
the usual definition of amenability for fixed unimodular quasi-transitive
graphs by Theorems 5.1 and 5.3 of \BLPSgip; it agrees with the usual notion
of amenability for equivalence relations by \ref t.many/ below; and it
allows us to extend to non-amenable unimodular random rooted graphs many
theorems that are known for non-amenable unimodular fixed graphs, as we
shall see.

We say that a graph $\gh$ is {\bf anchored amenable}
if there is a sequence of subsets $H_n\subset \vertex(\gh)$
such that $\bigcap_n H_n \ne \emptyset$, each $H_n$ induces a connected
subgraph of $\gh$, and
$$
\lim_{n \to\infty} {|\bdv H_n| \over |\verts(H_n)|} = 0\,.
$$
The relationship between amenability of $\rtd$ and amenability or anchored
amenability of $\rtd$-a.e.\ graph is as follows.
The first clearly implies the third, which implies the second, but the third
does not imply the first.
Indeed, take a 3-regular tree and
randomly subdivide its edges by a number of vertices whose
distribution does not have a finite
exponential tail. \ref b.ChenPeres/ show that the result is anchored
amenable a.s.
However, there is an appropriate unimodular version if the subdividing
distribution has finite mean (see Example 2.4.4 of \ref
b.Kaim:harmonic/ or \ref x.stretched/ below),
and it is non-amenable by \ref c.amentrees/ below.

In order to work with this definition, we shall need some easy facts.

\procl l.bdryo
If $\rtd$ is a unimodular probability measure on $\GG_*$ and $(\rtd',
\marks_0) \in \fdom(\rtd)$, then
$$
\int {\ct\big(K(\marks_0), \marks_0\big) \over |K(\marks_0)|} d\rtd'
=
\int \ct(\bp, \marks_0) \, d\rtd'
\,.
$$
\endprocl

\proof
Let $K_x$ be the $\marks_0$-component of $x$.
Let each vertex $x$ send mass $\ct(y, \marks_0)/|K_x|$ to each $y \in K_x$.
Then the left-hand side is the expected mass sent from the root and the
right-hand side is the expected mass received by the root.
\Qed

%

\procl p.alpha
If $\rtd$ is a unimodular probability measure on $\GG_*$, then
$$
\isoe(\rtd) + \alpha(\rtd) = \expdeg(\rtd)
\,.
\label e.alpha
$$
Therefore, if\/ $\expdeg(\rtd) < \infty$, then $\rtd$ is amenable iff
$\alpha(\rtd) = \expdeg(\rtd)$.
\endprocl

\proof
This is obvious from \ref l.bdryo/.
\Qed

\procl l.disint
Let $\rtd, \nu \in \uni$ with $\nu$ a percolation on $\rtd$, that is, there
is some Borel $\psi : \marks \to \marks$ such that $\rtd = \nu \circ
\psi^{-1}$.
Let $\kappa$ be a regular conditional probability measure of $\rtd$
with respect to the $\sigma$-field generated by $\psi$, i.e., a
disintegration of $\nu$ with respect to $\psi$, with $\kappa_{(\gh, \bp)}$
being the probability measure on the fiber over $(\gh, \bp)$.
Let $h : \gtwo \to [0, 1]$ be Borel and symmetric: $h(G, x, y) = h(G, y,
x)$.
Define $k(G, x, y) := \int h(G, x, y) \,d\kappa_{(\gh, x)}$.
Then there is a symmetric Borel $\lambda$ such that for $\rtd$-a.e.\ $(\gh,
\bp)$ and all $x \in \vertex(\gh)$, we have $k(\gh, \bp, x) = \lambda(\gh,
\bp, x)$.
\endprocl


\proof
It suffices to show that for all $f : \gtwo \to [0, 1]$, we have 
$$
\int \sum_{x \in \vertex(\gh)} k(\gh, \bp, x) f(\gh, \bp, x) \,d\rtd(\gh,
\bp)
=
\int \sum_{x \in \vertex(\gh)} k(\gh, x, \bp) f(\gh, \bp, x) \,d\rtd(\gh,
\bp)
\,,
$$
since this shows symmetry of $k$ a.e.\ with respect to the left measure
$\rtd_{\rm L}$.
To see that this equation holds, observe that 
$$\eqaln{
\int \sum_{x \in \vertex(\gh)} k(\gh, \bp, x) f(\gh, \bp, x) &\,d\rtd(\gh,
\bp)
=
\int \int \sum_{x \in \vertex(\gh)} h(\gh, \bp, x) f(\gh, \bp, x)
\,d\kappa_{(\gh, \bp)} \,d\rtd(\gh, \bp)
\cr&=
\int \int \sum_{x \in \vertex(\gh)} h(\gh, x, \bp) f(\gh, x, \bp)
\,d\kappa_{(\gh, \bp)} \,d\rtd(\gh, \bp)
\cr&\hskip1in\hbox{[by the Mass-Transport Principle for $\nu$]}
\cr&=
\int \int \sum_{x \in \vertex(\gh)} h(\gh, \bp, x) f(\gh, \bp, x)
\,d\kappa_{(\gh, x)} \,d\rtd(\gh, \bp)
\cr&\hskip1in\hbox{[by the Mass-Transport Principle for $\rtd$]}
\cr&=
\int \int \sum_{x \in \vertex(\gh)} h(\gh, x, \bp) f(\gh, \bp, x)
\,d\kappa_{(\gh, x)} \,d\rtd(\gh, \bp)
\cr&\hskip1in\hbox{[by symmetry of $h$]}
\cr&=
\int \sum_{x \in \vertex(\gh)} k(\gh, x, \bp) f(\gh, \bp, x) \,d\rtd(\gh,
\bp)
\,.
\Qed
}$$

We now prove some properties that are equivalent to amenability.
Most of these are standard in the context of equivalence relations.
With appropriate modifications, these equivalences hold with a weakening of
the assumption of unimodularity.
They are essentially due to \ref b.CFW/ and \ref b.Kaim:amen/, although
(ii) seems to be new.

\procl t.many \procname{Amenability Criteria}
Let $\rtd$ be a unimodular probability measure on $\GG_*$ with $\expdeg(\rtd)
< \infty$.
The following are equivalent:
\beginitems
\itemrm{(i)} $\rtd$ is amenable;
\itemrm{(ii)} there is
a sequence of Borel functions $\lambda_n : \gtwo \to [0, 1]$
such that for all $(\gh, x, y) \in \gtwo$ and all $n$,
we have
$$
\lambda_n(\gh, x, y) =
\lambda_n(\gh, y, x)
\label e.symm
$$
and for $\rtd$-a.e.\ $(\gh, \bp)$, we have
$$
\sum_{x \in \vertex(\gh)} \lambda_n(\gh, \bp, x) = 1
\label e.prmsrae
$$
and
$$
\lim_{n \to\infty}
\sum_{x \in \vertex(\gh)}\sum_{y \sim x}
\left|\lambda_n(\gh, \bp, x) -
\lambda_n(\gh, \bp, y)\right|
= 0
\,;
\label e.asymsame
$$
\itemrm{(iii)} there is
a sequence of Borel functions $\lambda_n : \gtwo \to [0, 1]$
such that for $\rtd$-a.e.\ $(\gh, \bp)$,
$$
\sum_{x \in \vertex(\gh)} \lambda_n(\gh, \bp, x) = 1
$$
and
$$
\lim_{n \to\infty} \int
\sum_{y \sim \bp} \sum_{x \in \vertex(\gh)}
\left|\lambda_n(\gh, \bp, x) -
\lambda_n(\gh, y, x)\right|
\,d\rtd(\gh, \bp)
= 0
\,;
$$
\itemrm{(iv)}
$\rtd$ is hyperfinite, meaning that there
is a unimodular measure $\nu$ on $\GG_*$, an increasing sequence of 
Borel subsets $\marks_n \subseteq \marks$, and a Borel
function $\psi: \marks \to \marks$ such that if $G$ denotes a
network with law $\nu$ and $G_n$ the subnetwork consisting of those edges
both of whose edge marks lie in $\marks_n$, then $\psi(G)$ has law $\rtd$,
all components of $G_n$ are finite, and $\bigcup_n \marks_n = \marks$.
\enditems
\endprocl

\proof
Assume from now on that $\rtd$ is carried by networks with distinct marks.
We shall use the following construction.
Suppose that $(\nu, \marks_0) \in \fdom(\rtd)$.
Let $\kappa$ be a regular conditional probability measure of
$\rtd$ with respect to the $\sigma$-field generated by $\projection$.
By \ref l.disint/, there is a Borel symmetric 
$\lambda : \gtwo \to [0, 1]$ such that
for $\rtd$-a.e.\ $(\gh, \bp)$ and $x \in \vertex(\gh)$, we have
$$
\lambda(\gh, \bp, x)
=
\int \II{x \in K(\marks_0)}/|K(\marks_0)|
\,d\kappa_{(\gh, \bp)}
\,.
\label e.connect
$$
Clearly,
$$
\sum_{x \in \vertex(\gh)} \lambda(\gh, \bp, x) = 1
\label e.prmsr
$$
for $\rtd$-a.e.\ $(\gh, \bp)$.
For $\rtd$-a.e.\ $(\gh, \bp)$, we have
$$\eqaln{
\sum_{x \in \vertex(\gh)} \sum_{y \sim x}
|\lambda(\gh, \bp, x) &- \lambda(\gh, \bp, y)|
\cr&\le
\sum_{x \in \vertex(\gh)} \sum_{y \sim x}
\int |\II{x \in K(\marks_0)} -
\II{y \in K(\marks_0)}|/|K(\marks_0)|
\,d\kappa_{(\gh, \bp)}
\cr&\le
\int {2\ct\big(K(\marks_0), \marks_0\big) \over |K(\marks_0)|}
\,d\kappa_{(\gh, \bp)}
\,.
\label e.isobd
}$$

Now if (i) holds, then we may choose $(\rtd_n, \marks_n) \in \fdom(\rtd)$
such that
$$
\int
\sum_{n}
{\ct\big(K(\marks_n), \marks_n\big) \over |K(\marks_n)|}
\,d\rtd_{n} < \infty
\,.
$$
Let $\kappa^{(n)}$ and
$\lambda_n$ be as above (but for $\rtd_n$).
Then by \ref e.isobd/, we have
$$\eqaln{
\sum_{n} \int
\sum_{x \in \vertex(\gh)}\sum_{y \sim x}
\big|\lambda_n(\gh, \bp, x) &-
\lambda_n(\gh, \bp, y)\big|
\,d\rtd(\gh, \bp)
\cr&\le
\sum_{n} \int
\int {2\ct\big(K(\marks_n), \marks_n\big) \over |K(\marks_n)|}
\,d\kappa^{(n)}_{(\gh, \bp)}
\,d\rtd(\gh, \bp)
\cr&=
\sum_{n}
\int {2\ct\big(K(\marks_n), \marks_n\big) \over |K(\marks_n)|}
\,d\rtd_{n} < \infty
\,,
}$$
which shows that (ii) holds.

Next, suppose that (ii) holds.
The Mass-Transport Principle and \ref e.symm/
show that the integral in \ref e.prmsrae/ is the same as
$$
\int
\sum_{y \sim \bp} \sum_{x \in \vertex(\gh)}
\left|\lambda_n(\gh, \bp, x) -
\lambda_n(\gh, y, x)\right|
\,d\rtd(\gh, \bp)
\,.
$$
This gives (iii).

Next, suppose that (iii) holds.
Then we may define a sequence similar to $\lambda_n$ on
the corresponding equivalence relation (see \ref x.equivalence/),
which implies that the equivalence
relation is amenable by \ref b.Kaim:amen/, and hence hyperfinite by
a theorem of \ref b.CFW/. (Another proof of the latter theorem was sketched
by \ref b.Kaim:amen/, with more details
given by \ref b.KechrisMiller/).
Translating the definition of hyperfinite equivalence relation to rooted
networks gives (iv).

Finally, that (iv) implies (i) is an immediate consequence of Lebesgue's
Dominated Convergence Theorem, our assumption that $\expdeg(\rtd) <
\infty$, and \ref l.bdryo/.
\Qed

Now we show how to produce unimodular networks from non-unimodular ones on
amenable measures, just as we can produce invariant measures from
non-invariant ones on amenable groups.
We illustrate in the context of couplings.

\procl p.amen-average \procname{Coupling From Amenability}
Let $\rln \subseteq \marks \times \marks$ be a closed set.
If $\rtd_1, \rtd_2 \in \uni$ are amenable and $\rtd_1$ is
$\rln$-related to $\rtd_2$, then there is a unimodular $\rln$-coupling of
$\rtd_1$ to $\rtd_2$.
\endprocl

\proof
Let $\nu$ be an $\rln$-coupling of $\rtd_1$ to $\rtd_2$.
Let $\lambda_n$ be as in \ref t.many/(ii) for $(\rtd_1 + \rtd_2)/2$.
Define the measures $\nu_n$ by
$$
\nu_n(\ev B)
:=
\int
\sum_{x \in \vertex(\gh)} \lambda_n(\gh, \bp, x) \II{(\gh, x) \in \ev B}
\,d\nu(\gh, \bp)
$$
for Borel $\ev B \subseteq \GG_*$ (with mark space $\marks \times \marks$).
Then $\nu_n$ is a probability measure by \ref e.prmsrae/.
Since $\nu$ is carried by networks all of whose marks are in $\rln$, so is
$\nu_n$.
If $\ev B$ is an event that
specifies only the first coordinates of the marks, i.e., $\ev B = \ev
B_1 \times 2^\marks$ for some $\ev B_1$, then
$$\eqaln{
\nu_n(\ev B)
&=
\int
\sum_{x \in \vertex(\gh)} \lambda_n(\gh, \bp, x) \II{(\gh, x) \in \ev B}
\,d\nu(\gh, \bp)
\cr&=
\int
\sum_{x \in \vertex(\gh)} \lambda_n(\gh, \bp, x) \II{(\gh, x) \in \ev B_1}
\,d\mu_1(\gh, \bp)
\cr&=
\int
\sum_{x \in \vertex(\gh)} \lambda_n(\gh, x, \bp) \II{(\gh, \bp) \in \ev B_1}
\,d\mu_1(\gh, \bp)
\cr&=
\int
\II{(\gh, \bp) \in \ev B_1}
\,d\mu_1(\gh, \bp)
\cr&=
\mu_1(\ev B)
}$$
by the Mass-Transport Principle, \ref e.symm/, and \ref e.prmsrae/.
Likewise, if $\ev B$ is an event that
specifies only the second coordinates of the marks, then $\nu_n(\ev B) =
\mu_2(\ev B)$.
Thus, $\nu_n$ is an
$\rln$-coupling of $\rtd_1$ to $\rtd_2$.

Now by definition of $\nu_n$, we have
$$
\int f(\gh, \bp)\,d\nu_n(\gh, \bp)
=
\int
\sum_{x \in \vertex(\gh)} \lambda_n(\gh, \bp, x) f(\gh, x)
\,d\nu(\gh, \bp)
$$
for every Borel $f : \GG_* \to [0, \infty]$.
Therefore, for every Borel $h : \gtwo \to [0, 1]$ with $h(\gh, x, y) = 0$
unless $x \sim y$, we have
$$\eqaln{
\int \sum_{y \in \vertex(\gh)} h(\gh, \bp, y) \,d\nu_n(\gh, \bp)
&=
\int \sum_{x \in \vertex(\gh)} \lambda_n(\gh, \bp, x)
\sum_{y \in \vertex(\gh)} h(\gh, x, y) \,d\nu_n(\gh, \bp)
\cr&=
\int \sum_{y \in \vertex(\gh)}
\sum_{x \in \vertex(\gh)} 
h(\gh, x, y) 
\lambda_n(\gh, \bp, x)
\,d\nu(\gh, \bp)
}$$
and
$$\eqaln{
\int \sum_{y \in \vertex(\gh)} h(\gh, y, \bp) \,d\nu_n(\gh, \bp)
&=
\int \sum_{x \in \vertex(\gh)} \lambda_n(\gh, \bp, x)
\sum_{y \in \vertex(\gh)} h(\gh, y, x) \,d\nu_n(\gh, \bp)
\cr&=
\int \sum_{y \in \vertex(\gh)}
\sum_{x \in \vertex(\gh)} 
h(\gh, x, y) 
\lambda_n(\gh, \bp, y)
\,d\nu(\gh, \bp)
\,,
}$$
where, in the last step, we have interchanged $x$ and $y$.
Therefore,
$$\eqaln{
\bigg| 
\int \sum_{y \in \vertex(\gh)} h(\gh, \bp, y) \,d\nu_n(\gh, \bp)
&- 
\int \sum_{y \in \vertex(\gh)} h(\gh, y, \bp) \,d\nu_n(\gh, \bp)
\bigg|
\cr&\le
\int \sum_{y \in \vertex(\gh)}
\sum_{x \sim y}
\left|
\lambda_n(\gh, \bp, x)
-
\lambda_n(\gh, \bp, y)
\right|
\,d\nu(\gh, \bp)
\cr&=
\int \sum_{y \in \vertex(\gh)}
\sum_{x \sim y}
\left|
\lambda_n(\gh, \bp, x)
-
\lambda_n(\gh, \bp, y)
\right|
\,d\rtd_1(\gh, \bp)
\,,
}$$
which tends to 0 by \ref e.asymsame/.
Thus, any limit point of $\nu_n$ is involution invariant and, since $\rln$ is
closed, is an $\rln$-coupling of $\rtd_1$ to $\rtd_2$.
\Qed

\comment{
\proof
Let $\nu$ be an $\rln$-coupling of $\rtd_1$ to $\rtd_2$.
Let $\lambda_n$ be as in \ref t.many/ for, say, $(\mu_1+\mu_2)/2$.
Write
$$
k(\gh, \bp) :=
\sum_{x \in \vertex(\gh)} \lambda_n(\gh, x, \bp)
\,.
$$
Define the measures $\nu_n$ by
$$
\nu_n(\ev B)
:=
\int
\sum_{x \in \vertex(\gh)} \lambda_n(\gh, \bp, x) \II{(\gh, x) \in \ev B}
k(\gh, x)^{-1} \,d\nu(\gh, \bp)
$$
for Borel $\ev B \subseteq \GG_*$ (with mark space $\marks \times \marks$).
Then $\nu_n$ is a probability measure since
$$\eqaln{
\int
\sum_{x \in \vertex(\gh)} \lambda_n(\gh, \bp, x)
k(\gh, x)^{-1} \,d\nu(\gh, \bp)
&=
\int
\sum_{x \in \vertex(\gh)} \lambda_n(\gh, x, \bp)
k(\gh, \bp)^{-1} \,d\nu(\gh, \bp)
\cr&=
\int d\nu(\gh, \bp)
= 1
}$$
by the Mass-Transport Principle.
Since $\nu$ is carried by networks all of whose marks are in $\rln$, so is
$\nu_n$.
If $\ev B$ is an event that
specifies only the first coordinates of the marks, i.e., $\ev B = \ev
B_1 \times 2^\marks$ for some $\ev B_1$, then
$$\eqaln{
\nu_n(\ev B)
&=
\int
\sum_{x \in \vertex(\gh)} \lambda_n(\gh, \bp, x) \II{(\gh, x) \in \ev B}
k(\gh, x)^{-1} \,d\kappa(\gh, \bp)
\cr&=
\int
\sum_{x \in \vertex(\gh)} \lambda_n(\gh, \bp, x) \II{(\gh, x) \in \ev B_1}
k(\gh, x)^{-1} \,d\mu_1(\gh, \bp)
\cr&=
\int
\sum_{x \in \vertex(\gh)} \lambda_n(\gh, x, \bp) \II{(\gh, \bp) \in \ev B_1}
k(\gh, \bp)^{-1} \,d\mu_1(\gh, \bp)
\cr&=
\int
\II{(\gh, \bp) \in \ev B_1}
\,d\mu_1(\gh, \bp)
\cr&=
\mu_1(\ev B)
}$$
by the Mass-Transport Principle.
Likewise, if $\ev B$ is an event that
specifies only the second coordinates of the marks, then $\nu_n(\ev B) =
\mu_2(\ev B)$.
Thus, $\nu_n$ is an
$\rln$-coupling of $\rtd_1$ to $\rtd_2$.
Now by definition of $\nu_n$, we have
$$
\int f(\gh, \bp)\,d\nu_n(\gh, \bp)
=
\int
\sum_{x \in \vertex(\gh)} \lambda_n(\gh, \bp, x) f(\gh, x)
k(\gh, x)^{-1} \,d\nu(\gh, \bp)
$$
for every Borel $f : \GG_* \to [0, \infty]$.
Therefore,
\msnote{The below does not calculate in sufficient generality.}
$$\eqaln{
\bigg| \int \sum_{y \sim \bp} \II{(\gh, \bp) \in \ev B} \,d\nu_n(\gh, \bp)
&- \int \sum_{y \sim \bp} \II{(\gh, y) \in \ev B} \,d\nu_n(\gh, \bp)\bigg|
\cr&\quad=
\bigg|
\int \sum_{x \in \vertex(\gh)} \lambda_n(\gh, \bp, x)
\sum_{y \sim x} \II{(\gh, x) \in \ev B} k(\gh, x)^{-1} \,d\nu(\gh, \bp)
\cr&\qquad-
\sum_{x \in \vertex(\gh)} \lambda_n(\gh, \bp, x)
\sum_{y \sim x} \II{(\gh, y) \in \ev B} k(\gh, x)^{-1} \,d\nu(\gh, \bp)
\bigg|
\cr&\quad\le
...
}$$
\msnote{Now apply CBS and use that the square of a number in
$[0, 1]$ is at most that number.}
Thus, any limit point of $\nu_n$ is unimodular and, since $\rln$ is
closed, is an $\rln$-coupling of $\rtd_1$ to $\rtd_2$.
\Qed
}%

\procl p.recurrent \procname{Recurrence Implies Amenability}
If $\rtd \in \uni$ and simple random walk is $\rtd$-a.s.\ recurrent,
then $\rtd$ is amenable.
\endprocl

\proof
Consider the ``lazy" simple random walk that moves nowhere with probability
1/2 and otherwise moves to a random neighbor, like simple random walk.
It is recurrent by hypothesis and aperiodic by construction.
Let $\lambda_n(\gh, \bp, x)$ be the probability that lazy simple random
walk on $\gh$ starting from $\bp$ will be at $x$ at time $n$.
By \ref b.Orey/ and recurrence, the functions $\lambda_n$
satisfy property (iii) of \ref t.many/, whence $\rtd$ is amenable.
\Qed

The following is proved similarly to Remark 6.2 of \BLPSusf.

\procl p.amenforest \procname{Forests in Amenable Networks}
If $\rtd \in \uni$ is amenable and $\marks_0 \subseteq \marks$ is such that
the $\marks_0$-open subgraph $\fo$ of $(\gh, \bp)$ is a forest $\rtd$-a.s.,
then the expected degree of $\bp$ in $\fo$ is at most 2.
\endprocl

The following extends Theorem 5.3 of \BLPSgip.
In the following, we say that $\P$ is a {\bf percolation} on $\rtd$ that
{\bf gives} subgraphs $\A$ a.s.\ if there is a Borel
function $\psi: \marks \to \marks$ such that $\rtd = \P \circ \psi^{-1}$
and there is a Borel subset $\marks_0 \subseteq \marks$
such that if $G(\marks_0)$ denotes the $\marks_0$-open subnetwork of $G$,
then $\psi\big(G(\marks_0)\big) \in \A$ for $\P$-a.e.\ $G$.

\procl t.treechar
Let $\rtd \in \uni$ with $\expdeg(\rtd) < \infty$. The following are equivalent:
\beginitems
\itemrm{(i)} $\rtd$ is amenable;
\itemrm{(ii)} there is a percolation $\P$ on $\rtd$ that gives
spanning trees with at most 2 ends a.s.;
\itemrm{(iii)} there is a percolation $\P$ on $\rtd$ that gives
non-empty connected subgraphs $\omega$ that satisfy $\pc(\omega) = 1$ a.s.
\enditems
\endprocl

\proof
The proof that (i) implies (ii) is done as for Theorem 5.3 of \BLPSgip, but
uses Propositions \briefref p.amenforest/ and \briefref p.limittree/.
That (ii) implies (iii) is obvious.
The proof that (iii) implies (i) follows the first part of the proof of
Theorem 1.1 in \ref b.BLPS:crit/.
\Qed

\procl c.amentrees \procname{Amenable Trees}
A unimodular probability measure $\rtd$ on infinite
rooted trees is amenable iff $\expdeg(\rtd) = 2$ iff $\rtd$-a.s.\ $\gh$ has
1 or 2 ends.
\endprocl

\proof
Combine \ref t.treechar/ with \ref t.deg2/.
\comment{
Since the average degree in any finite tree is less than 2, we have
$\expdeg(\rtd') < 2$ for all $\rtd' \in \fdom(\rtd)$, and so $\alpha(\rtd) \le 2$.
Therefore, $\expdeg(\rtd) > 2$ implies non-amenability by \ref p.alpha/.
On the other hand, if $\expdeg(\rtd) = 2$, then by \ref t.deg2/, $\rtd$ is
concentrated on trees with at most 2 ends, whence on trees with $\pc = 1$.
Therefore, Bernoulli($p$) percolation on $\rtd$ produces only finite
components and for $p$ close to 1, produces a $\rtd' \in \fdom(\rtd)$ with
$\expdeg(\rtd')$ close to 2.
Thus, $\alpha(\rtd) = 2$ and we have amenability by \ref p.alpha/.
The last equivalence is from \ref t.deg2/.
}%
\Qed

The next result was proved for non-amenable unimodular transitive graphs in
\BLPSgip\ with a more direct proof in \ref b.BLPS:crit/. This extension is
proved similarly.
Presumably, the hypothesis that $\rtd$ is non-amenable can be replaced by the
assumption that $\pc(\rtd) < 1$. (This is a major open conjecture for
quasi-transitive graphs.)

\procl t.death \procname{Critical Percolation}
Let $\rtd$ be an extremal unimodular non-amen\-a\-ble
probability measure on $\GG_*$
with $\expdeg(\rtd) < \infty$.
There is $\P_{\pc(\rtd)}$-a.s.\ no infinite cluster.
\endprocl


This can be interpreted for finite graphs as follows.
Suppose that $\gh_n$ are finite connected graphs with bounded average
degree whose random weak limit is
extremal and non-amenable with critical value $\pc$.
Consider Bernoulli($\pc$) percolation on $\gh_n$.
Let $\alpha_n(\ell)$ be the random variable giving
the proportion of vertices of $\gh_n$ that belong to simple open paths of
length at least $\ell$.
Then $\lim_{\ell \to\infty} \lim_{n \to\infty} \alpha_n(\ell) = 0$ in
probability.

The following theorem is proved similarly to \ref b.BS:hp/, as extended by
\ref b.LP:book/, and by using \ref x.dual/.

\procl t.planarperc \procname{Planar Percolation}
Let $\rtd \in \uni$ be extremal, non-amenable,
and carried by plane graphs with one end and bounded degree.
Then $0 < \pc(\rtd) < \pu(\rtd) < 1$ and Bernoulli($\pu(\rtd)$) percolation
on $\rtd$ has a unique infinite cluster a.s.
\endprocl

The following extends Theorems 3.1, 3.2, 3.5, 3.6, and 3.10 of \ref
b.BLS:pert/ and is proved similarly.
If $\P$ is a percolation on $\rtd$ with $\P \circ \psi^{-1} = \rtd$ and
with $\marks_0 \subseteq \marks$
defining the open subgraphs, we call $\P'$ a {\bf subpercolation} on $\P$
that {\bf gives} subgraphs in $\A'$ {\bf with positive probability} if there
is a Borel
function $\psi': \marks \to \marks$ such that $\P = \P' \circ \psi'^{-1}$
and there is a Borel subset $\marks_1 \subseteq \marks$ such that if
$G(\marks_1)$ denotes the $\marks_1$-open subnetwork of $G$, then
$\P'\left[\psi\Big(\psi'\big(G(\marks_1)\big)\big(\marks_0\big)\Big) \in
\A'\right] > 0$.
For a graph $G$, define
$$
\isoe(G)
:=
\inf \Bigl\{ \frac{|\{(x, y) \st x \in K,\, y \notin K,\, (x, y) \in
\edges\}|}{|K|} \st K \subset \verts \hbox{ is finite}\Bigr\}
\,.
$$

\procl t.sub-iso \procname{Non-Amenable Subgraphs}
Let $\rtd$ be a unimodular probability measure on $\GG_*$ with finite
expected degree and $\P$ be a
percolation on $\rtd$ with open subgraph $\omega$.
\beginitems
\itemrm{(i)} If $h > 0$ and $\E[\deg_\omega \bp \mid \bp \in \omega] \ge
\alpha(\rtd) + 2h$, then there is a subpercolation $\P'$ on $\P$ that gives
a non-empty subgraph $\omega'$ with $\isoe(\omega') \ge h$ with positive
$\P'$-probability.
\itemrm{(ii)} If $\omega$ is a forest a.s., $h > 0$ and $\E[\deg_\omega \bp
\mid \bp \in \omega] \ge 2 + 2h$, then there is a subpercolation $\P'$ on
$\P$ that gives a non-empty subgraph with $\isoe \ge h$ with positive
probability.
\itemrm{(iii)} If $\rtd$ is non-amenable and extremal
and $\omega$ has exactly one
infinite cluster $\P$-a.s., then there is a subpercolation $\P'$ on $\P$
that gives a non-empty subgraph $\omega'$ with $\isoe(\omega') > 0$ $\P'$-a.s.
\itemrm{(iv)} If $\omega$ has components with at least three ends
$\P$-a.s., then there is a subpercolation $\P'$ on $\P$ that gives a
non-empty forest $\fo$ with $\isoe(\fo) > 0$ $\P'$-a.s.
\itemrm{(v)} If $\rtd$ is concentrated on subgraphs with spectral radius
less than 1 and $\omega$ has exactly one infinite cluster $\P$-a.s., then
there is a subpercolation $\P'$ on $\P$ that gives a non-empty forest $\fo$
with $\isoe(\fo) > 0$ $\P'$-a.s.
\enditems
\endprocl

The following extends a result of \ref b.Hag:rcust/ and is proved
similarly to Corollary 6.3 of \BLPSusf.

\procl p.amen-for \procname{Amenability and Boundary Conditions}
Let $\rtd$ be an amenable unimodular probability measure on $\GG_*$.
Then $\FUSF(\rtd) = \WUSF(\rtd)$ and $\FMSF(\rtd) = \WMSF(\rtd)$.
\endprocl

We may now strengthen \ref p.recurrent/, despite the fact that not every
graph is necessarily non-amenable $\rtd$-a.s.
It extends Theorem 4.3 of \ref b.BLS:pert/ and is proved similarly, using
\ref t.sub-iso/(iii) with $\P := \rtd$.

\procl t.posspeed \procname{Positive Speed on Non-Amenable Graphs}
If $\rtd \in \uni$ is non-amenable and concentrated on graphs with bounded
degree,
then the speed of simple random walk is positive $\rtd$-a.s.
\endprocl

\bsection{Examples}{s.ex}

We present here a variety of interesting examples of unimodular
measures.

\procl x.renew \procname{Renewal Processes}
Given a stationary (delayed) renewal process on $\Z$, let $\rtd$ be the law
of $(\Z, 0)$ with the graph $\Z$, some fixed mark at renewals, and some
other fixed mark elsewhere. Then $\rtd \in \uni$.
\endprocl

\procl x.halfplane \procname{Half-Plane}
Fix $d \ge 3$ and let $T$ be the $d$-regular tree.
Let $\rtd_1$ be the random weak limit of balls of growing radii in $T$.
Note that $\rtd_1$ is carried by trees with only one end.
Let $\rtd_2$ be concentrated on the fixed graph $(\Z, 0)$.
Now let $\rtd := \rtd_1 \itimes \rtd_2$.
This is a unimodular version of the half-plane $\N \times \Z$.
\endprocl

\procl x.cover
For a rooted network $(\gh, \bp)$, its {\bf universal cover} is the rooted
tree $(T, \bp) = T(\gh, \bp)$ formed as follows.
The vertices of $T$ are the finite paths in $\gh$ that start at $\bp$ and
do not backtrack.
Two such vertices are joined by an edge in $T$ if one is an extension of
the other by exactly one edge in $\gh$.
The path with no edges consisting of just the vertex $\bp$ in $\gh$ is the
root $\bp$ of $T$.
There is a natural rooted graph homomorphism $\pi : (T, \bp) \to (\gh,
\bp)$ (the cover map) that maps paths to their last point.
Marks on $T$ are defined by lifting the marks on $\gh$ via $\pi$.
It is clear that
if $\rtd \in \uni$ and $\nu$ is the law of $T(\gh, \bp)$ when $(\gh, \bp)$
has the law $\rtd$, then $\nu \in \uni$
and $\expdeg(\rtd) = \expdeg(\nu)$.
\endprocl

\procl x.cluster
Let $\P$ be a unimodular percolation on $\rtd$ that labels edges either
open or closed.
Let $\nu$ be the law of the open cluster of the root when the network is
chosen according to $\P$ conditional on the root belonging to an infinite
open cluster.
Then $\nu$ is unimodular, as a direct verification of the definition shows.
When $\rtd$ is concentrated on a fixed unimodular graph, this fact has been
widely used in the study of percolation.
\endprocl

\procl x.Voronoi \procname{Tilings}
Let $X$ be a Euclidean space or hyperbolic space (of constant curvature).
Write $\gp$ for its isometry group.
There is a Mass-Transport Principle for $X$ that says the following; see
\ref b.BS:hp/ for a proof.
Let $\rho$ be a positive Borel measure on $X \times X$ that is invariant
under the diagonal action of $\gp$.
Then there is a constant $c$ such that for all Borel $A \subset X$ of
volume $|A| > 0$, we have $\rho(A \times X) = \rho(X \times A) = c |A|$.
Suppose that $P$ is a (countable) point process in $X$ whose law is
$\gp$-invariant.
For example, Poisson point processes are $\gp$-invariant.
One often considers graphs $\gh$ that are functions $\gh = \beta(P)$, where
$\beta$ commutes with the action of $\gp$.
For a few recent examples, see
\ref b.BS:hp/, \ref b.HolPer:PP/, or \ref b.Timar:PP/.
For instance, the 1-skeleton of the Voronoi tessellation corresponding to
$P$ is such a graph.
In general, we call such measures on graphs {\bf $\gp$-equivariant factors}
of $P$.
They are necessarily $\gp$-invariant.
\par
Another way that invariant measures on graphs embedded in $X$ occur is
through (aperiodic) tilings of $X$.
Again, one can take the 1-skeleton.
An important tool for studying aperiodic tilings is a limit measure
obtained from translates of a given tiling (in the Euclidean case), or,
more generally, invariant measures on tilings with special properties; see,
e.g., \ref b.Rob/, \ref b.Radin:survey/, \ref b.Solomyak/,
\ref b.Radin:book/, or
\ref b.BowRad:densest/ for some examples.
\par
Let $\nu$ be any $\gp$-invariant probability measure on graphs embedded in
$X$.
Fix a Borel set $A \subset X$ of positive finite volume.
If $v(A) := \int |\vertex(\gh) \cap A| \,d\nu(\gh) < \infty$,
then define $\rtd$ as follows.
Choose $\gh$ with the law $\nu$ biased by $|\vertex(\gh) \cap A|$.
Then choose the root $\bp$ of $\gh$
uniformly among all vertices that belong to $A$.
The law of the resulting graph $(\gh, \bp)$ is $\rtd$.
We claim that $\rtd$ is unimodular and does not depend on $A$.
In fact, $\rtd$ is the same as the Palm measure of $(\gh, \bp)$, except
that $\rtd$ is a measure on isomorphism classes of graphs that
does not involve any geometric embedding.
To prove our claims, we first write $\rtd$ in symbols:
$$
\rtd(\A) :=
v(A)^{-1} \int \sum_{\bp \in A} \II{(\gh, \bp) \in \A} \,d\nu(\gh)
$$
for Borel $\A \subset \GG_*$.
Let $f : \gtwo \to [0, \infty]$ be Borel.
Define
$$
\rho(B \times C) :=
\int \sum_{x \in \vertex(\gh) \cap B} \sum_{y \in \vertex(\gh) \cap C}
f(\gh, x, y)\,d\nu(\gh)
$$
for Borel $B, C \subseteq X$.  Since $\nu$ is invariant, $\rho$ is
diagonally invariant.
Therefore,
$$\eqaln{
\int \sum_{x \in \vertex(\gh)} f(\gh, \bp, x) \,d\rtd(\gh, \bp)
&=
v(A)^{-1} \rho(A \times X)
=
v(A)^{-1} \rho(X \times A)
\cr&=
\int \sum_{x \in \vertex(\gh)} f(\gh, x, \bp) \,d\rtd(\gh, \bp)
\,,
}$$
which means that $\rtd$
satisfies the Mass-Transport Principle, i.e., is unimodular.
Furthermore, if we take $f(\gh, x, y) := \II{x = y}$, then we see that
$\rho(A \times X) = v(A)$, so that
there is a constant $c$ such that $v(A) = c |A|$.
Likewise, if $f(\gh, x, y) := \II{x = y, (\gh, x) \in \ev A}$, then we see
that for each $\A$, there is another constant $c_\A$ such that $v(A)
\rtd(\A) = c_\A |A|$.
It follows that $\rtd$ does not depend on $A$.
\endprocl


\procl x.dual \procname{Planar Duals}
Let $\rtd$ be a unimodular probability measure on plane graphs all of
whose faces have finitely many sides.
We are assuming that to each graph, there is a measurably associated
plane embedding.
Thus, each graph $\gh$ has a plane dual $\dual \gh$ with respect to its
embedding.
In fact, to be technically accurate in what follows, we replace the
embedding by an assignment (possibly random) of marks to the edges that
indicate the cyclic order in which they appear around a vertex
in a fixed orientation of the plane.
(For example, if a vertex $x$ has $d$ edges incident to it, then one can
let the $d$ edge marks associated to $x$ be $\{1, 2, \ldots, d\}$ in cyclic
order, with the one marked 1 chosen at random, independently of marks
elsewhere.)
Then the plane dual graph is defined entirely with respect to the
resulting network in an automorphism-equivariant way, needing no reference
to the plane.
\par
Provided a certain finiteness condition is satisfied,
there is a natural unimodular probability measure on the dual graphs,
constructed as follows.
For a face $f$, let $\deg f$ denote the
number of sides of $f$. For a vertex $x$, let $F(x) := \sum_{f \sim x}
1/\deg f$.
Assume that $Z := \int F(\bp) \,d\rtd(\gh, \bp) < \infty$.
To create a unimodular probability measure $\dual \rtd$ on the duals,
first choose $(\gh, \bp)$ with law $\rtd$ biased by $F(\bp)/Z$.
Then choose a face $f_0$ incident to $\bp$ with probability proportional to
$1/\deg f_0$.
The law of the resulting rooted graph $(\dual \gh, f_0)$ is $\dual \rtd$:
$$
\dual \rtd(\A) :=
Z^{-1} \int \sum_{f_0 \sim \bp} {1 \over \deg f_0} \II{(\dual \gh, f_0) \in
\A}
\,d\rtd(\gh, \bp)
$$
for Borel $\A \subseteq \GG_*$.
To prove that $\dual \rtd$ is indeed unimodular, let $k : \gtwo \to [0, \infty]$
be Borel. Then
$$\eqaln{
Z \int \sum_{f \in \vertex(\dual \gh)} k(\dual \gh, f_0, f)
\,d\dual \rtd(\dual \gh, f_0)
&=
\int
\sum_{f_0 \sim \bp} {1 \over \deg f_0}
\sum_{f \in \vertex(\dual \gh)} k(\dual \gh, f_0, f)
\,d\rtd(\gh, \bp)
\cr&=
\int
\sum_{f_0 \sim \bp}
\sum_{f \in \vertex(\dual \gh)}
{1 \over \deg f_0}
k(\dual \gh, f_0, f)
\sum_{x \sim f} {1 \over \deg f}
\,d\rtd(\gh, \bp)
\cr&=
\int
\sum_{x \in \vertex(\gh)}
\sum_{f_0 \sim \bp}
\sum_{f \sim x}
{1 \over \deg f_0}
{1 \over \deg f}
k(\dual \gh, f_0, f)
\,d\rtd(\gh, \bp)
\cr&=
\int
\sum_{x \in \vertex(\gh)}
\sum_{f_0 \sim x}
\sum_{f \sim \bp}
{1 \over \deg f_0}
{1 \over \deg f}
k(\dual \gh, f_0, f)
\,d\rtd(\gh, \bp)
\cr&\hskip1in\hbox{[by the Mass-Transport Principle for $\rtd$]}
\cr&=
\int
\sum_{x \in \vertex(\gh)}
\sum_{f \sim x}
\sum_{f_0 \sim \bp}
{1 \over \deg f_0}
{1 \over \deg f}
k(\dual \gh, f, f_0)
\,d\rtd(\gh, \bp)
\cr&=
Z \int \sum_{f \in \vertex(\dual \gh)} k(\dual \gh, f, f_0)
\,d\dual \rtd(\dual \gh, f_0)
\,.
}$$
Thus, $\dual \rtd$ satisfies the Mass-Transport Principle, so is unimodular.
A similar argument shows that $\dual{(\dual \rtd)} = \rtd$.
\par
Another important construction comes from combining the primal and dual
graphs into a new plane graph by adding a vertex where each edge crosses
its dual.
That is, if $\gh$ is a plane graph and $\dual \gh$ its dual, then every
edge $e \in \edge(\gh)$ intersects $\dual e \in \edge(\dual \gh)$ in one
point, $v_e$.
(These are the only intersections of $\gh$ and $\dual \gh$.)
For $e \in \edge(\gh)$, write $\widehat e$ for the pair of edges that
result from the subdivision of $e$ by $v_e$, and likewise for
$\widehat{\dual e}$.
This defines a new graph $\widehat \gh$, whose vertices are
$\verts(\gh)\cup\verts(\dual \gh)\cup\Big\{v_e\st e\in\edges(\gh)\Big\}$
and whose edges are $\bigcup_{e \in\edges(\gh)} (\widehat e \cup
\widehat{\dual e})$.
If $\expdeg(\rtd) < \infty$, then we may define a unimodular probability
measure $\widehat \rtd$ on the graphs $\widehat \gh$ from $\rtd$ as
follows.
Let $\widehat Z := 1 + (1/2) \expdeg(\rtd) + Z \le (5/2) \expdeg(\rtd) <
\infty$.
For $x \in \vertex(\gh)$, let $\widehat N(x)$ be the set consisting of $x$
itself plus the vertices of $\widehat \gh$ that correspond to edges or
faces of $\gh$ that are incident to $x$.
For $w \in \vertex(\widehat \gh)$, define
$$
\delta(w) :=
|\{x \in \vertex(\gh) \st w \in \widehat N(x)\}|^{-1}
=
\cases{
1 &if $w \in \vertex(\gh)$,\cr
1/2 &if $w = v_e$ for some $e \in \edge(\gh)$,\cr
1/\deg w &if $w \in \vertex(\dual \gh)$.\cr
}
$$
Define
$$
\widehat \rtd(\A) :=
\widehat Z^{-1} \int \sum_{w_0 \in \widehat N(\bp)} \delta(w_0)
\II{(\widehat \gh, w_0) \in \A} \,d\rtd(\gh, \bp)
\,.
$$
Note that $\widehat\rtd$ is a probability measure.
To prove that $\widehat\rtd$ is unimodular, let $k : \gtwo \to [0, \infty]$
be Borel. Then
$$\eqaln{
\widehat Z \int \sum_{w \in \vertex(\widehat \gh)} k(\widehat \gh, w_0, w)
&\,d\widehat\rtd(\widehat \gh, w_0)
=
\int
\sum_{w_0 \in \widehat N(\bp)} \delta(w_0)
\sum_{w \in \vertex(\widehat \gh)} k(\widehat \gh, w_0, w)
\,d\rtd(\gh, \bp)
\cr&=
\int
\sum_{w_0 \in \widehat N(\bp)}
\sum_{w \in \vertex(\widehat \gh)}
\delta(w_0)
k(\widehat \gh, w_0, w)
\sum_{x \st w \in \widehat N(x)} \delta(w)
\,d\rtd(\gh, \bp)
\cr&=
\int
\sum_{x \in \vertex(\gh)}
\sum_{w_0 \in \widehat N(\bp)}
\sum_{w \in \widehat N(x)}
\delta(w_0)
\delta(w)
k(\widehat \gh, w_0, w)
\,d\rtd(\gh, \bp)
\cr&=
\int
\sum_{x \in \vertex(\gh)}
\sum_{w_0 \in \widehat N(x)}
\sum_{w \in \widehat N(\bp)}
\delta(w_0)
\delta(w)
k(\widehat \gh, w_0, w)
\,d\rtd(\gh, \bp)
\cr&\hskip1in\hbox{[by the Mass-Transport Principle for $\rtd$]}
\cr&=
\int
\sum_{x \in \vertex(\gh)}
\sum_{w \in \widehat N(x)}
\sum_{w_0 \in \widehat N(\bp)}
\delta(w_0)
\delta(w)
k(\widehat \gh, w, w_0)
\,d\rtd(\gh, \bp)
\cr&=
\widehat Z \int \sum_{w \in \vertex(\widehat \gh)} k(\widehat \gh, w, w_0)
\,d\widehat \rtd(\widehat \gh, w_0)
\,.
}$$
Thus, $\widehat\rtd$ satisfies the Mass-Transport Principle, so is unimodular.

\endprocl

\procl x.PWIT \procname{Poisson Weighted Infinite Tree}
Our definition
(\ref s.notation/)
of the metric on the space $\GG_*$ of rooted graphs
refers to ``balls of radius $r$" in which distance is
graph distance, i.e., edges implicitly have length $1$.
\ref b.AS:obj/ work in the setting of graphs whose edges
have positive real lengths, so that distance becomes minimum
path length.  This setting permits one to consider 
graphs which may have infinite degree, but which are still
``locally finite" in the sense that only finitely vertices fall within
any finite radius ball.
Of course, edge lengths are a special (symmetric) case of edge marks.
An important example is the following.
Consider a regular rooted tree $T$ of infinite degree.  
Fix a continuous increasing function $\Lambda$ on $\CO{0, \infty}$ with
$\Lambda(0) = 0$ and $\lim_{t \to\infty} \Lambda(t) = \infty$.
Order the children of each vertex of $T$ via a bijection with $\Z^+$.
For each vertex $x$, consider an independent Poisson process on $\R^+$ with
mean function $\Lambda$.
Define the length of the edge joining $x$ to its $n$th child to be the
$n$th point of the Poisson process associated to $x$.
This is a unimodular random network in the extended sense of \ref b.AS:obj/.
It can be derived by taking the random weak limit of the complete graph on $n$
vertices whose edge lengths are independent with cdf 
$t \mapsto 1 - e^{-\Lambda(t)/n}$ ($t \ge 0$)
and then deleting all edges in the limit whose length
is $\infty$.
(We are working here with the mark space $[0, \infty]$.)
See \ref b.Aldous:random-assign/.
\endprocl


\procl x.stretched \procname{Edge Replacement}
Here is a general way to create unimodular random rooted graphs from existing
unimodular fixed graphs.
This is an extension of the random subdivision (or stretching) introduced
by \ref b.AdamsLyons/ and studied further in
Example 2.4.4 of \ref b.Kaim:harmonic/ and \ref b.ChenPeres/.
Let $\fin$ be
the set of isomorphism classes of finite graphs with
an ordered pair of distinct distinguished vertices.
For our construction, we may start with a fixed
unimodular quasi-transitive connected graph, $\gh$, or, more generally,
with any unimodular probability measure $\rtd$ on $\GG_*$.
In the former case, fix an orientation of the edges of $\gh$ and let $L$ be
a random field on the oriented edges of $\gh$ that is invariant under the
automorphism group of $\gh$ and takes values in $\fin$ and such that
$\big|\vertex\big(L(e)\big)\big|$ has finite mean for each edge $e$.
Replace each edge $e$ with the graph $L(e)$, where the first of the
distinguished vertices of $L(e)$ is identified with the tail of $e$ and the
second of its distinguished vertices is identified with the head of $e$.
Call the resulting random graph $H$.
There is a unimodular probability measure that is equivalent to this
measure on random graphs.
Namely, let $\mu_i$ be the law of $(H, o_i)$, where $\{o_i\}$ is a complete
section of the vertex orbits of $\gh$.
Given $H$, write
$$
A(x) := 2 + \sum_{e \sim x} \Big(\big|\vertex\big(L(e)\big)\big|-2\Big)
$$
and
$$
c := \sum_i \E\big[A(o_i)\big] \big|\Stab(o_i)\big|^{-1}
\,.
$$
Choose $o_i$ with probability $c^{-1} \E\big[A(o_i)\big]
\big|\Stab(o_i)\big|^{-1}$.
Given $o_i$, choose $(H, o_i)$ with distribution $\mu_i$.
Given this, list the non-distinguished vertices of all $L(e)$ for $e$
incident to $o_i$ as $z_1, z_2, \ldots, z_{A(o_i) - 2}$ and set $z_{A(o_i) - 1}
:= z_{A(o_i)} := o_i$.
Let $U$ be a uniform integer in $\big[1, A(o_i)\big]$.
Then $(H, z_U)$ is unimodular and, clearly, has law with respect to which
$\sum_i \mu_i$ is absolutely continuous.

Indeed, we state and prove this more generally.
Suppose that $\rtd$ is a unimodular probability measure on $\GG_*$.
Orient the edges of the rooted networks arbitrarily.
Let $\mkmp(e)$ denote the ordered pair of the
marks of the edge $e$ (ordered by the orientation of $e$).
Suppose $L : \marks^2 \to \fin$ is Borel with the property
that whenever $L(\mk_1, \mk_2) = (G, x, y)$, we also have $L(\mk_2, \mk_1) =
(G, y, x)$.
(This will ensure that the orientation of the edges will not affect the
result.)
If
$$
\int \sum_{e \sim \bp} \bigg[\Big|\vertex\Big(L\big(\mkmp(e)\big)\Big)\Big|
- 2\bigg]
\,d\rtd(\gh, \bp)
< \infty
\,,
$$
then let $\rtd'$ be the following measure.
Define
$$
A(\gh, \bp) := 2 + \sum_{e \sim \bp}
\bigg[\Big|\vertex\Big(L\big(\mkmp(e)\big)\Big)\Big| - 2\bigg]
\,.
$$
Choose $(\gh, \bp)$ with probability distribution $\rtd$ biased by $A(\gh,
\bp)$ and replace each edge $e$ by the graph $L\big(\mkmp(e)\big)$, where
the tail and head of $e$ are identified with the first and second
distinguished vertices of $L\big(\mkmp(e)\big)$, respectively; call the
resulting graph $H$.
Write $A := A(\gh, \bp)$ and
list the non-distinguished vertices of all $L\big(\mkmp(e)\big)$ for $e$
incident to $\bp$ as $z_1, z_2, \ldots, z_{A - 2}$ and set $z_{A - 1}
:= z_{A} := \bp$.
Let $U$ be a uniform integer in $[1, A]$.
Finally, let $\rtd'$ be the distribution of $(H, z_U)$.

This is unimodular by the following calculation.
Write $H(\gh)$ for the graph $H$ formed as above from the network $\gh$.
Let $\vertex_0(\mk_1, \mk_2)$ be the set of non-distinguished vertices of
the graph $L(\mk_1, \mk_2)$.
Write $z_i(\gh, \bp)$ ($1 \le i \le A(\gh, \bp)-2$) for the vertices of the
neighborhood $B(\gh, \bp) := \bigcup_{e \sim \bp} \vertex_0\big(\mkmp(e)\big)$.
Write $z_i(\gh, \bp) := \bp$ for $i = A(\gh, \bp)-1, A(\gh, \bp)$.
Put $c := \int A(\gh, \bp) \,d\rtd(\gh, \bp)$.
In order to show that $\rtd'$ is unimodular, let $f : \gtwo \to [0,
\infty]$ be Borel.
Define $\overline f : \gtwo \to [0, \infty]$ by
$$\eqaln{
\overline f(\gh, x, y) 
&:=
{1 \over c}
\sum_{z \in B(\gh, x)}
\sum_{z' \in B(\gh, y)}
f\big(H(\gh), z, z'\big)
+
{2 \over c}
\sum_{z \in B(\gh, x)}
f\big(H(\gh), z, y\big)
\cr&\qquad+
{2 \over c}
\sum_{z' \in B(\gh, y)}
f\big(H(\gh), x, z'\big)
+
{2 \over c}
f\big(H(\gh), x, y\big)
\,.
\cr
}$$
Then 
$$\eqaln{
\int \sum_{z \in \vertex(H)} &f(H, \bp, z) \,d\rtd'(H, \bp)
\cr&=
{1 \over c} \int {1 \over A(\gh, \bp)} \sum_{i=1}^{A(\gh, \bp)}
\sum_{z \in \vertex(H(\gh))} f\big(H(\gh), z_i(\gh, \bp), z\big)
A(\gh, \bp)\,d\rtd(\gh, \bp)
\cr&=
\int \sum_{x \in \vertex(\gh)} \overline f(\gh, \bp, x) \,d\rtd(\gh, \bp)
\cr&=
\int \sum_{x \in \vertex(\gh)} \overline f(\gh, x, \bp) \,d\rtd(\gh, \bp)
\cr&=
\int \sum_{z \in \vertex(H)} f(H, z, \bp) \,d\rtd'(H, \bp)
\,.
}$$
\endprocl

Our final example details the correspondence between random rooted graphs
and graphings of equivalence relations.

\procl x.equivalence
Let $\mu$ be a Borel probability measure on a topological space
$X$ and $R$ be a Borel subset of
$X^2$ that is an equivalence relation with finite or countable equivalence
classes. We call the triple $(X, \mu, R)$ a {\bf measured equivalence
relation}.
For $x \in X$, denote its $R$-equivalence class by $[x]$.
We call $R$ {\bf measure preserving} if
$$
\int_{x \in X} \sum_{y \in [x]} f(x, y) \,d\mu(x)
=
\int_{x \in X} \sum_{y \in [x]} f(y, x) \,d\mu(x)
$$
for all Borel $f : X^2 \to [0, \infty]$.
A {\bf graphing} $\Phi$ of $R$ is a Borel subset of $X^2$ such that the
smallest equivalence relation containing $\Phi$ is $R$.
A graphing $\Phi$ induces the structure of a graph on the vertex set $X$
by defining an edge between $x$ and $y$ if $(x, y) \in \Phi$ or $(y, x)
\in \Phi$.
Denote the subgraph induced on $[x]$ and rooted at $x$ by $\Phi(x)$.
Given Borel maps $\psi : X \to \marks$ and $\phi : X^2 \to \marks$,
we regard $\psi(x)$ as the mark at $x$ and $\phi(x, y)$ as the mark at $x$
of the edge from $x$ to $y$.
Thus, $\Phi(x)$ is a random rooted network. Its law (or, rather, the law of
its rooted isomorphism class) is unimodular iff $R$ is measure preserving.

Conversely, suppose that $\rtd$ is a probability measure on $\GG_*$.
Add independent uniform marks as second coordinates to the existing marks
and call the resulting measure $\nu$.
Write $\dmn \subset \GG_*$ for the set of (isomorphism classes of) rooted
networks with distinct marks.
Thus, $\nu$ is concentrated on $\dmn$.
Define $R \subset \dmn^2$ to be the set of pairs of (isomorphism classes
of) rooted networks that are
non-rooted isomorphic.
Define $\Phi \subset R$ to be the set of pairs of isomorphic rooted
networks whose roots are neighbors in the unique (non-rooted) isomorphism.
Then $(\dmn, \nu, R)$ is a measured equivalence relation with graphing
$\Phi$.
We have that $R$ is measure preserving iff $\rtd$ is unimodular.
If we define the mark map $\projection$ that forgets the second coordinate,
then $\nu$ pushes forward to $\rtd$, i.e., $\rtd = \nu \circ \projection^{-1}$.

Thus, the theory of unimodular random rooted networks has substantial
overlap with the
theory of graphed measure-preserving equivalence relations.
The largest difference between the two theories lies in the foci of
attention: We focus on probabilistic aspects of the graphing, while the
other theory focuses on ergodic aspects of the equivalence relation (and,
thus, considers all graphings of a given equivalence relation).
The origins of our work lie in two distinct areas: one is group-invariant
percolation on graphs, while the other is asymptotic analysis of finite
graphs. The origin of the study of measured equivalence relations lies in
the ergodic theory of group actions.
Some references for the latter work, showing relations to von Neumann
algebras and logic, among other things, are
\refbmulti{FM1,FM2}, \ref b.FHM/,
\ref b.CFW/,
\ref b.Zimmer:book/,
\ref b.KechrisMiller/,
\ref b.BeckerKech:book/, \ref b.Kechris:turing/,
\ref b.AdamsLyons/,
\refbmulti {Kaim:amen,Kaim:harmonic},
\ref b.Paulin/,
\refbmulti{Gaboriau:cost,Gaboriau:betti},
and \refbmulti{Furman1,Furman2}.
\endprocl

\bsection{Finite Approximation}{s.finite}

Although we do not present any theorems in this section, because of its
potential importance, we have devoted the whole section to the question of
whether
finite networks are weakly dense in
$\uni$.
Let us call random weak limits of finite networks {\bf sofic}.

\procl q.TII \procname{Finite Approximation}
Is every probability measure in $\uni$ sofic? In other words, if
$\rtd$ is a unimodular probability measure on $\GG_*$,
do there exist finite networks $G_n$ such that
$U(G_n) \cd \rtd$?
\endprocl

To appreciate why the answer is not obvious, consider the
special case of the (non-random) graph consisting of the
infinite rooted 3-regular tree.
It is true that there exist finite graphs $G_n$ that
approximate the infinite 3-regular tree in the sense of random
weak convergence; a moment's thought shows these {\it cannot}
be finite trees.
This special case is of course known
(one can use finite quotient groups of $\Z_2 * \Z_2 * \Z_2$,
random $3$-regular graphs, or expanders (\ref B.lub94/)),
but the constructions in this special case do not readily extend to
the general case.

Another known case of sofic measures is more difficult to establish.
Namely, \ref b.Bowen:periodic/ showed
that all unimodular networks on regular trees are sofic. (To deduce this from
his result, one must use the fact, easily established, that networks with
marks from a finite set are dense in $\GG_*$.)

\procl x.fdd
The general unimodular Galton-Watson
measure $\UGW$ (\ref x.AGW/) is also sofic.
To see this, consider the following random networks,
sometimes called ``fixed-degree
distribution networks" and first studied by \ref b.MolloyReed/.
Given
$\Seq{r_k}$ and $n$ vertices, give each vertex $k$ balls with probability
$r_k$, independently. Then pair the balls at random and place an edge for
each pair between the corresponding vertices. There may be one ball left
over; if so, ignore it. Let $m_0 := \sum k r_k$, which we assume is
finite.  In the limit, we get a tree where the root has degree $k$ with
probability $r_k$, each neighbor of the root, if any, has degree $k$ with
probability $k r_k/m_0$, etc. In fact, we get $\UGW$ for the
offspring distribution $k \mapsto (k+1) r_{k+1}/m_0$.  Thus, if we want
the offspring distribution $\Seq{p_k}$, we need merely start with $r_k :=
c^{-1} p_{k-1}/k$ for $k \ge 1$ and $r_0 := 0$, where $c := \sum_{k \ge 0}
p_k/(k+1)$.
\endprocl

Let us compare the
intuitions behind amenability and unimodularity.
One can define an average of any bounded function on the vertices of, say,
an amenable Cayley graph; the average will be the same for any translate of
the function.
This can also be regarded as an average of the function with respect to a
probability measure that chooses a group element uniformly at random;
however, the precise justification of this requires a measure that, though
group invariant, is only finitely additive.
Nevertheless, this invariant measure is approximated by uniform measures on
finite sets, namely, F\o{}lner sets.
By contrast, the justification that a unimodular random rooted graph
provides a uniform distribution on the vertices is via the Mass-Transport
Principle. The measure itself is, of course, countably additive; if it is
sofic, then it,
too, is approximated by uniform measures on finite sets.
The two intuitions concerning uniform measures that come from amenability
and from unimodularity agree insofar as every amenable quasi-transitive graph
is unimodular, as shown by \ref b.SoardiWoess/ and \ref b.Salvatori/.

One might think that if a sequence $\Seq{\gh_n}$ of finite graphs has a
fixed transitive graph $\gh$ as its random weak limit, then any unimodular
probability measure on networks supported by $\gh$ would be a random weak
limit of some choice of networks on the same sequence $\Seq{\gh_n}$.
This is false, however; e.g., if $\gh$ is a 3-regular tree, then almost any
choice of a sequence of growing 3-regular graphs has $\gh$ as its random weak
limit. However, there is a random independent set\ftnote{*}{This means that
no two vertices of the set are adjacent.} of density 1/2 on $\gh$ whose
law is automorphism invariant,
while the density of independent sets in random 3-regular graphs is bounded
away from 1/2 (see \ref b.FS/).
Nevertheless, if \ref q.TII/ has a positive answer, then
it is not hard to show that
there is some sequence of finite graphs that has this
property of carrying arbitrary networks.

Recall that the {\bf Cayley diagram} of a group $\gp$ generated by a finite
subset $S$ is the network $(\gh, \bp)$ with vertex set $\gp$, edge set
$\big\{(x, x s) \st x \in \gp,\ s \in S\big\}$, root $\bp$ the identity
element of $\gp$, and edge marks $s$ at the endpoint $x$ of $(x, x s)$ and
$s^{-1}$ at the endpoint $x s$ of $(x, x s)$, as in \ref r.anygroup/.
We do not mark the vertices (or mark them all the same).
\ref b.Weiss:sofic/ defined $\gp$ to be {\bf sofic} if its Cayley diagram 
is a random weak limit of finite networks with edge marks from
$S \cup S^{-1}$.
It is easy to check that this property does not depend on the generating
set $S$ chosen.
By embedding $S \cup S^{-1}$ into $\marks$, we can use a positive answer to
\ref q.TII/ to give every
Cayley diagram as a random weak limit of {\it some\/} finite networks.
Changing the marks on the finite networks to their nearest points in $S
\cup S^{-1}$ gives the kind of approximating networks desired. That is, we
would have that
every finitely generated group is sofic.
As we mentioned in the introduction, this would have plenty of
consequences.

To illustrate additional consequences of
a positive answer to
\ref q.TII/, we establish
the existence of various probability measures on sofic networks.
The idea is that if a class of probability measures is specified by a
sequence of local ``closed" conditions in such a way that there is a measure in
that class on any finite graph, then there is an automorphism-invariant
measure in that class on any sofic quasi-transitive graph.
Rather than state a general theorem to that effect, we shall
illustrate the principle by two examples.

\procl x.gibbs \procname{Invariant Markov Random Fields}
Consider networks with vertex marks $\pm 1$ and no (or constant) edge marks.
Given a finite graph $\gh$, $h \in \R$, and $\beta > 0$,
the probability measure $\nu_\gh$ on mark maps $\mkmp : \verts(\gh) \to \{-1,
+1\}$ given by
$$
\nu_\gh(\mkmp) := Z^{-1} \exp\Big\{\sum_{x \in \verts(\gh)} \beta h \mkmp(x) -
\sum_{x \sim y} \beta \mkmp(x) \mkmp(y)\Big\}
\,,
$$
where $Z$ is the normalizing constant required to make these probabilities
add to 1, is known as the
anti-ferromagnetic Ising model at inverse
temperature $\beta$ and with external field $h$ on $\gh$.
Let $\gh$ now be an infinite sofic transitive graph.
Let $\Seq{\gh_n}$ be an approximating sequence of finite graphs
and let $\nu_n := \nu_{\gh_n}$ be the corresponding probability measures.
Write $\overline \gh_n$ for the corresponding random network on $\gh_n$
and, as usual, $U(\overline \gh_n)$ for the uniformly rooted network.
Since $U(\overline \gh_n)$ is unimodular, so is any weak limit point,
$\rtd$.
By tightness, there is such a weak limit point, and it is concentrated on
networks with underlying graph $\gh$.
By \ref t.uni-vs-invar/, we may lift $\rtd$ to an
automorphism-invariant
measure $\nu = \lift_\rtd$ on networks on $\gh$.
This measure $\nu$ is a Markov random field with the required Gibbs
specification, meaning that for any finite subgraph $H$ of $\gh$, the
conditional $\nu$-distribution of $\mkmp \restrict \vertex(H)$ given $\mkmp
\restrict \bdv H$ is equal to the conditional $\nu_H$-distribution of the
same thing.
One is often interested in invariant random fields, not just any random
fields with the given Gibbs specification.
One can easily get a Markov random field with the required Gibbs
specification by taking a limit over subgraphs of $\gh$, but this will not
necessarily produce an invariant measure. In case $\gh$ is amenable, one
could take a limit of averages of the resulting measure to obtain an
invariant measure, but this will not work in the non-amenable case.
That is the point of the present construction.
A variation on this is spin glasses:
Here, for finite graphs $\gh$, the measure $\nu_\gh$ is
$$
\nu_\gh(\mkmp) := Z^{-1} \exp\Big\{\sum_{x \in \verts(\gh)} \beta h \mkmp(x) +
\sum_{x \sim y} \beta J_{x, y} \mkmp(x) \mkmp(y)\Big\}
\,,
$$
where $J_{x, y}$ are, say, independent $\pm 1$-valued random variables.
Again, one can find an invariant spin glass (coupled to the independent
interactions $J_{x, y}$) with the same parameters $h$ and $\beta$ on any
transitive sofic graph by the above method.
\endprocl

\procl x.sandpile \procname{Invariant Sandpiles}
Consider now networks with vertex marks in $\N$ and no edge marks.
Given a finite rooted 
graph $(\gh, \bp)$, a mark map $\mkmp : \vertex(\gh)
\to \N$ is called {\bf critical} (or {\bf stable and recurrent}) if for all
$x \in \vertex(\gh)$, we have $\mkmp(x) < \deg(x)$ and for all
subgraphs $W$ of $\gh \setminus \{\bp\}$, there is some $x \in \vertex(W)$
such that $\mkmp(x) \ge \deg_W(x)$; we may take $\mkmp(\bp) \equiv 0$.
In this context, one usually calls the root ``the sink".
It turns out that the set of such mark
maps form a very interesting group, called the {\bf sandpile
group} or {\bf chip-firing group} of $\gh$; see \ref b.btw/, \ref
b.Dhar:survey/, \ref b.Biggs/, and \ref b.MRZ/.
Let $\nu_{(\gh, \bp)}$ be the uniform measure on critical mark maps.
Given a transitive sofic graph $\gh$ and an approximating sequence
$\Seq{\gh_n}$ of finite graphs, let $\nu_n :=
\nu_{(\gh_n, \bp_n)}$ be the corresponding probability measures for any
fixed choice of roots $\bp_n$.
Write $(\overline \gh_n, \bp_n)$ for the corresponding random network on
$(\gh_n, \bp_n)$ and, as usual, $U(\overline \gh_n)$ for the uniformly rooted
network. (This root is unrelated to $\bp_n$.)
Since $U(\overline \gh_n)$ is unimodular, so is any weak limit point,
$\rtd$.
By tightness, there is such a weak limit point, and it is concentrated on
networks with underlying graph $\gh$.
(Probably the entire sequence $\Seq{U(\overline \gh_n)}$ in fact converges
to $\rtd$.)
By \ref t.uni-vs-invar/, we may lift $\rtd$ to an
automorphism-invariant
measure $\nu$ on networks on $\gh$.
This measure $\nu$ is supported by networks with only critical mark maps in
the sense that for all $x \in \vertex(\gh)$, we have $\mkmp(x) < \deg(x)$
and for all subgraphs $W$ of $\gh$, there is some $x \in \vertex(W)$ such
that $\mkmp(x) \ge \deg_W(x)$.
\endprocl

However, we do not know how to answer the following question.

\procl q.color \procname{Invariant Coloring}
Given a quasi-transitive infinite graph $\gh$ and a number $c$
that is 
at least the chromatic number of $\gh$, is there an $\Aut(\gh)$-invariant
probability measure on proper $c$-colorings of the vertices of $\gh$?
Let $D := \max_{x \in \vertex(\gh)} \deg_\gh x $.
A positive answer for $c \ge D + 1$ is due to Schramm (personal
communication, 1997).
If $\gh$ is sofic, then we can also obtain such a measure for $c = D$ by
using a well-known result of Brooks (see, e.g., \ref b.Bollobas/, p.~148,
Theorem V.3.)
The question is particularly interesting when $\gh$ is planar and $c = 4$.
In fact, it is then also of great interest to know whether there is a
quasi-transitive proper 4-coloring of $\gh$.
\endprocl

\medbreak
\noindent {\bf Acknowledgements.}\enspace
We are grateful to 
Oded Schramm for discussions of Theorems \briefref
t.uni-vs-invar/ and \briefref t.msf1end/.
We also thank Oded Schramm and \'Ad\'am Tim\'ar for various suggestions.
We are indebted to
G\'abor Elek for telling us about the connection of our work to
various conjectures about groups.

\def\noop#1{\relax}
\input \jobname.bbl

\filbreak
\begingroup
\eightpoint\sc
\parindent=0pt\baselineskip=10pt

Department of Statistics,
University of California Berkeley,
Berkeley, CA, 94720-3860
\emailwww{aldous@stat.berkeley.edu}
{http://www.stat.berkeley.edu/users/aldous}

Department of Mathematics,
Indiana University,
Bloomington, IN 47405-5701
\emailwww{rdlyons@indiana.edu}
{http://mypage.iu.edu/\string~rdlyons/}

\endgroup

\bye